\documentclass[11pt,a4paper]{amsart}
\usepackage{amsmath,amsfonts,amssymb,amsthm,amscd,mathrsfs}
\usepackage{stmaryrd,yhmath}
\usepackage{cases}
\usepackage[english]{babel}
\usepackage[T1]{fontenc}
\usepackage{aeguill}
\usepackage{enumerate}
\usepackage{hyperref}
\usepackage{mathrsfs}
\usepackage{tikz}

\usepackage{graphics,epsfig,color,cite,bbm}

\definecolor{light}{gray}{.9}

\definecolor{bleu}{RGB}{27,88,145}
\definecolor{mauve}{RGB}{138,20,79}

\newcommand{\mb}[1]{\mathbb{#1}}

\newcommand{\mc}[1]{\mathcal{#1}}

\newcommand{\R}[1]{\mb{R}^{#1}}

\newcommand{\lp}{\langle}
\newcommand{\rp}{\rangle}

\newcommand{\field}[1]{\mathbb{#1}} 
\newcommand{\cz}{\field{C}}

\def\N{{\mathbb N}}

\def\R{{\mathbb R}}
\def\C{{\mathbb C}}

\def\div{{\rm div}\,}
\def\Hess{{\rm Hess\,}}

\def\supp{\mathop{\rm supp} \nolimits} 
\def\Id{{\rm Id}}
\def\Ran{{\rm Ran}}
\def\sspan{{\rm Span}}

\def\and {{\rm \, and \,}}
\def\dist {{\rm \, dist \,}}

\def\Ker {{\rm \, Ker  \,}}
\def\Im {{\rm \, Im\,}}

\def\Re {{\rm \,Re\,}}

\def\dim {{\rm \, dim  \,}}

\def \rank {{\rm \, rank \,}}

\def \diag {{\rm \, diag \,}}

\newcommand {\pa}{\partial}

\makeatletter

\@addtoreset{equation}{section}
\makeatother

\renewcommand{\Im}{\operatorname{Im}}
\renewcommand{\Re}{\operatorname{Re}}
\newcommand{\bsigma}{\boldsymbol{\sigma}}

\renewcommand{\div}{\operatorname{div}}

\def\<{\langle}
\def\>{\rangle}
\newcommand{\bp}{{\it Proof. }}
\newcommand{\ep}{\hfill $\square$\\}

\def\m{\mathbf{m}}
\def\s{\mathbf{s}}

\def\u{\mathbf{u}}

\def\um{\underline\m}

\newcommand{\be}{\begin{equation}}
\newcommand{\ee}{\end{equation}}
\newcommand{\bes}{\begin{equation*}}
\newcommand{\ees}{\end{equation*}}

\numberwithin{equation}{section}
\numberwithin{figure}{section}

\def\ccc{{\mathcal C}}\def\ddd{{\mathcal D}}
\def\eee{{\mathcal E}} 
 \def\jjj{{\mathcal J}}
\def\mmm{{\mathcal M}}\def\nnn{{\mathcal N}} \def\ooo{{\mathcal O}}\def\ppp{{\mathcal P}}
\def\sss{{\mathcal S}}
\def\uuu{{\mathcal U}}\def\vvv{{\mathcal V}}
  \def\zzz{{\mathcal Z}}

\def\Dr{{\mathscr D}}
\def\Er{{\mathscr E}}\def\Gr{{\mathscr G}}

\def\Mr{{\mathscr M}}

\def\ulu{\underline{\mathcal U}}

\newtheorem{theorem}{Theorem}[section]
\newtheorem{lemma}[theorem]{Lemma}

\newtheorem{proposition}[theorem]{Proposition}
\newtheorem{definition}[theorem]{Definition}
\newtheorem{remark}[theorem]{Remark}
\newtheorem{corollary}[theorem]{Corollary}

\newtheorem{assumption} {Assumption}

\newtheorem{theorem*}{Theorem}

\begin{document}

\title[Sharp  asymptotics for non-reversible diffusion processes]{Sharp spectral asymptotics for non-reversible metastable diffusion processes}

\author{Dorian Le Peutrec}
\thanks{D. Le Peutrec : Institut Denis Poisson, Universit\'e d'Orl\'eans, Universit\'e de Tours, CNRS, Orl\'eans, France. E-mail: dorian.le-peutrec@univ-orleans.fr}

\author{Laurent Michel}
\thanks{L. Michel : Universit\'e de Bordeaux,
Institut Math\'ematiques de Bordeaux,
Talence,
France. E-mail: laurent.michel@math.u-bordeaux.fr}

\maketitle

\begin{abstract}
Let $U_h:\R^{d}\to \R^{d}$ be a smooth vector field
 and consider  the associated overdamped Langevin equation $$dX_t=-U_h(X_t)\,dt+\sqrt{2h}\,dB_t$$  in the low temperature regime $h\rightarrow 0$. In this work, we study the spectrum of the  associated 
diffusion $L=-h\Delta+U_h\cdot\nabla$ under the assumptions that
$U_h=U_{0}+h\nu$, where
the vector fields
$U_{0}:\R^{d}\to \R^{d}$ and $\nu:\R^{d}\to \R^{d}$ are independent of  $h\in(0,1]$,
and  that the dynamics
 admits $e^{-\frac Vh}dx$ as an invariant measure
  for some smooth function $V:\mathbb{R}^d\rightarrow\mathbb{R}$. Assuming additionally that  $V$ is a Morse function admitting $n_0$ local minima, we prove that there exists $\epsilon>0$ such that in the limit   $h\to 0$, $L$ admits exactly
 $n_0$ eigenvalues in the strip $\{0\leq \operatorname{Re}(z)< \epsilon\}$, which  have moreover exponentially small moduli. Under a generic assumption
on the potential barriers of the Morse function $V$, we also prove that  the asymptotic behaviors of these small eigenvalues are 
given by Eyring-Kramers type formulas.\\
\end{abstract}

\noindent
\textbf{MSC 2010:} 60J60, 35Q82, 81Q12, 35P15, 81Q20.\\
\textbf{Keywords:} Non-reversible overdamped Langevin dynamics, Metastability, Spectral theory, Semiclassical analysis, Eyring-Kramers formulas.


%
 
 \tableofcontents
 
 \section{Introduction}
 
\subsection{Setting and motivation}
\label{sub.motivation}
 
Let $d\geq 2$, $U_h:\R^d\rightarrow\R^d$ be a smooth vector field depending on a small parameter $h\in(0,1]$, and consider the associated overdamped Langevin equation
\be\label{eq:overdampedL}
dX_t\ =\ -U_h(X_{t})\,dt+\sqrt{2h}\,dB_t\,,
\ee
where $X_{t}\in \R^{d}$ and $(B_t)_{t\geq 0}$ is a standard Brownian motion in $\R^d$. 
The associated Kolmogorov (backward) and Fokker-Planck equations are then the evolution equations
\be\label{eq:KS}
\partial_t\,u+L(u)\ =\ 0\quad\text{and}\quad\partial_t\,\rho+L^{\dagger}(\rho)\ =\ 0\,,
\ee
where the elliptic differential operator 
$$L\ =\ -h\Delta+U_h\cdot\nabla$$ 
is the infinitesimal generator of the process \eqref{eq:overdampedL},
$$L^{\dagger}\ =\ -\div \,\circ\, (h\nabla +U_{h})$$ 
denotes the formal adjoint of $L$, 
and for $x\in \R^{d}$ and $t\ge0$: $u(t,x)=\mathbb E^{x}[f(X_{t})]$ is the expected
value of the observable $f(X_{t})$ when $X_{0}=x$
and  $\rho(t,\cdot)$ is the probability density 
(with respect to the Lebesgue measure on $\R^{d}$) of presence of $(X_t)_{t\geq 0}$.
In this setting, the Fokker--Planck equation, that is the second equation of \eqref{eq:KS}, is also known as the  
Kramers-Smoluchowski equation.\medskip

Throughout this paper, we assume that the vector field $U_h$ decomposes as 
$$U_h\ =\ U_0+h\nu$$
 for some real smooth vector fields $U_0$ and $\nu$ independent of $h$. Moreover, we consider the case where the above 
 overdamped Langevin
 dynamics admits a specific stationary distribution satisfying the following assumption:
\begin{assumption}\label{hyp:statmes}
There exists a smooth function $V:\R^d\rightarrow \R$ such that  $L^{\dagger }(e^{-\frac Vh})=0$
 for every $h\in(0,1]$.
\end{assumption}

A straightforward computation shows that 
Assumption~\ref{hyp:statmes}
is satisfied if and only if
the vector field $U_h=U_0+h\nu$ satisfies the following relations, where we denote 
$b:=U_0-\nabla V$,
\be\label{eq:propvectfield}
b\cdot\nabla V\ =\ 0\,,\ \ 
\div(\nu)\ =\ 0\,,\ \ \text{and}\ \ 
\div(b)\ =\ \nu\cdot\nabla V\,.
\ee
Using this decomposition, the generator $L$ writes
\be
\label{eq.L}
L_{V,b,\nu} \ :=\ L \  =\ -h\Delta + \nabla V\cdot \nabla + b_h \cdot\nabla \,,
\ee
where 
\begin{equation}
\label{eq.bh}
b_h\ :=\ b+ h\nu\ =\ U_0-\nabla V + h \nu\ =\ U_{h}-\nabla V\,.
\end{equation} 
Note moreover that the two following particular cases enter in the framework of 
Assumption~\ref{hyp:statmes}:
\begin{itemize}
\item[1.] The case where
\be
\label{eq.J-div}
b\cdot\nabla V\ =\ 0\,,\ \ \div b\ =\ 0\ \ \text{and}\ \ \nu\ =\ 0\,,
\ee
which is in particular satisfied when $\nu=0$ and $b$ is the matrix product
$b=J\,\nabla V$, where $J$ is a smooth map from $\R^{d}$
into the set of real antisymmetric matrices of size $d$ such that $\div\big (J\,\nabla V\big)=0$. For instance, this later condition holds if
$J(x)=\tilde J\circ V(x)$ for some antisymmetric matrices $\tilde J(y)$ depending smoothly on $y\in\R$. 
\item[2.] The case where 
\be
\label{eq.J-supersym}
b\ =\ J\,\nabla V\quad\text{and}\quad \nu\ =\ \Big( \,\sum_{i=1}^{d}\pa_{i} J_{ij}\,\Big)_{1\leq j\leq d}\,,
\ee
where $J$ is a smooth map from $\R^{d}$
into the set of real antisymmetric matrices of size $d$. 
\end{itemize}

 In the  case of \eqref{eq.J-supersym}, $L_{V,b,\nu}$
has in particular the following supersymmetric-type structure,
\be
\label{eq.J-supersym'}
L_{V,b,\nu} \ =\ -h \,e^{\frac Vh} \,\div \circ \big( \,e^{-\frac Vh} \big( I_{d} - J \big)\nabla \,  \big)\,,
\ee
and both cases coincide when $b_{h}$ has the form $b_{h}=b=J\,\nabla V$
for some constant antisymmetric matrix $J$. In the case of \eqref{eq.J-div}, the structure \eqref{eq.J-supersym'} fails to be true in general and we refer to 
\cite{Mi16} for more details on these questions. Let us also point out  that under Assumption \ref{hyp:statmes}, the vector 
field~$b_h$ defined in \eqref{eq.bh} is very close to the transverse vector field introduced in \cite{BoRe16} and next used in \cite{LaMaSe19}.\medskip

In this paper, we are interested in the spectral analysis
 of the operator~$L_{V,b,\nu}$ and in its connections with
 the long-time behaviour of the dynamics~\eqref{eq:overdampedL}
 when $h\to 0$. In this regime,
 the process $(X_{t})_{t\geq 0}$ solution to \eqref{eq:overdampedL}  is typically metastable,
which is characterized by a very slow return to equilibrium. We refer
especially in this connection to the related works \cite{BoRe16,LaMaSe19} dealing with the
mean transition times
between the different wells of the potential~$V$ for the process  $(X_t)_{t\geq 0}$.
Our setting is also motivated by the question of accelerating the convergence to equilibrium,
which is of  interest for computational purposes.
It is indeed known that non-gradient perturbations of 
the overdamped  gradient Langevin dynamics 
\begin{equation}
\label{eq.ODL}
dX_t\  =\  -\nabla V(X_{t})\,dt+\sqrt{2h}\,dB_t
\end{equation}
which preserve the invariant measure $e^{-\frac Vh}dx$
cannot converge slower to equilibrium than the associated gradient dynamics \eqref{eq.ODL}.
See in particular \cite{LeNiPa} on this topic, where the authors considered linear drifts
and computed the optimal rate of return to equilibrium in this case, and references therein.

\subsection{Preliminary analysis}

 In view of Assumption~\ref{hyp:statmes},  we look 
 at $L_{V,b,\nu}$ acting in
 the  natural weighted Hilbert space $L^{2}(\R^{d},m_{h})$, where
\be\label{eq:defmesure}
m_{h}(dx)\ :=\  Z^{-1}_{h}e^{-\frac{V(x)}h}dx\quad \text{and}\quad  Z_{h}\ :=\ \int_{\R^{d}} e^{-\frac{V(x)}h}dx\,.
\ee
Note that we assume here that $e^{-\frac{V}h}\in L^{1}(\R^{d})$ for every $h\in(0,1]$,
which will be a simple consequence of our further hypotheses.
In this setting, a first important consequence of \eqref{eq:propvectfield} is the following identity,
easily deduced from the relation $\div(b_{h} e^{-\frac{V}h})=0$,
$$
\forall\,u,v\,\in\, \ccc^{\infty}_{c}(\R^{d})\,,\ \ \  \langle L_{ V, b,\nu} u,v \rangle_{L^{2}(m_{h})}\ =\ \langle  u, L_{V,-b,-\nu}v \rangle_{L^{2}(m_{h})}\,.
$$
 In particular, using \eqref{eq.L},  it holds
\begin{align}
\nonumber
\Re \< L_{ V, b,\nu} u,u \>_{L^{2}(m_{h})}&\ =\ \<(-h\Delta + \nabla V\cdot \nabla)u,u\>_{L^{2}(m_{h})}\\
\label{eq:minorLV}
&\ =\ 
h \,\| \nabla u \|^{2}_{L^{2}(m_{h})}
\ \geq\  0
\end{align}
for all $u\,\in\, \ccc^{\infty}_{c}(\R^{d})$
and the operator $L_{V,b,\nu}$ acting on $\ccc^{\infty}_{c}(\R^{d})$ in $L^{2}(\R^{d},m_{h})$
is hence accretive.\medskip

Let us now introduce the following confining  assumptions at infinity on the functions $V$, $b$, and $\nu$
that we will consider in the rest of this work.

\begin{assumption}
\label{ass.confin}
There exist $C>0$ and a compact  set $K\subset\R^d$ such that
it holds
$$V\ \geq\  -C\quad\text{on \ $\R^{d}$}$$
and,
 for all $x\in\R^d\setminus K$,
\begin{equation}\label{eq:hypgenephi}
 \vert\nabla V(x)\vert\ \geq\ \frac 1 C\quad\text{and}\quad \vert \Hess V(x)\vert\ \leq\  C\vert\nabla V(x)\vert^2\,.
\end{equation}
Moreover, there exists $C>0$ such that  the vector fields $b=U_0-\nabla V$ and $\nu$ satisfy the following estimate
for all $x\in\R^d$:
\be\label{hyp:controlechamp}
 \vert b(x)\vert+\vert \nu(x)\vert\ \leq\  C\,(1+\vert\nabla V(x)\vert).
\ee
\end{assumption}

 One can show that when $V$ is bounded from below and the first  estimate of \eqref{eq:hypgenephi} is satisfied,
it also holds, for some $C>0$, $V(x)\geq  C|x|$ outside a compact set (see for example \cite[Lemma~3.14]{MeSc}).
In particular, when Assumption~\ref{ass.confin} is satisfied, then $e^{-\frac V h}\in L^1(\R^{d})$ for all $h\in(0,1]$ (which justifies the definition of $Z_h$ in  \eqref{eq:defmesure}).\medskip

%

In order to study the operator $L_{V,b,\nu}$ in $L^2(\R^d,m_h)$, it is often useful
to work with its counterpart  in the flat  space $L^2(\R^d,dx)$ by  using  the unitary
transformation
$$
\mathsf{U}
\,:\, L^2(\R^d,dx)\longrightarrow L^2(\R^d,  m_h)\,,\ \  \mathsf{U}(u)\ =\ m_h^{-\frac 12}u=\ Z_h^{\frac 12}e^{\frac V {2h}}\,u\,.
$$
Defining $\phi:=\frac V2$, we then have the unitary equivalence 
\begin{align}
\nonumber
\mathsf{U}^{*} \,h\,L_{V,b,\nu}\,\mathsf{U}&\ =\ -h^{2}\Delta+|\nabla \phi|^{2}-h\Delta \phi + \, b_h  \cdot d_{\phi,h} \\
\label{eq.unit}
 &\ =\ \Delta_{\phi} +  b_h\cdot d_{\phi}\,,
\end{align}
where 
\begin{equation}
\label{eq.d-phi}
d_{\phi}\ :=\ d_{\phi,h}\ :=\ h\nabla +\nabla \phi\ =\ h e^{-\frac \phi h} \nabla e^{\frac \phi h}
\end{equation}
 and  
$$\Delta_{\phi}\ :=\ \Delta_{\phi,h}\ :=\ -h^{2}\Delta+|\nabla \phi|^{2}-h\Delta \phi
\ =\ -h^{2}  e^{\frac \phi h} \div   e^{-\frac \phi h} d_{\phi}$$ 
denotes the usual 
semiclassical
Witten Laplacian
acting on functions. 
It is thus equivalent to study $L_{V,b,\nu}$ acting in the weighted space $L^2(\R^d,m_h)$
or
\be
\label{eq.P}
P_{\phi}\ :=\ P_{\phi,b,\nu}\ :=\ \Delta_{\phi}+b_h\cdot d_{\phi}
\ee 
acting in the flat space $L^2(\R^d,dx)$.\medskip

The Witten Laplacian $\Delta_\phi =P_{\phi,0,0}$, which is the counterpart of
the weighted Laplacian 
$$L_{V,0,0}\ =\ -h\Delta + \nabla V\cdot \nabla=h \nabla^{*}\nabla   $$ (the adjoint is considered here with respect to $m_{h}$)
acting in the flat space $L^{2}(\R^{d},dx)$,
  is moreover essentially self-adjoint on $\ccc^\infty_c(\R^n)$
(see \cite[Theorem~9.15]{Hel-spectral}). We still denote by $\Delta_\phi$ its unique self-adjoint extension and by $D(\Delta_\phi)$ the domain
of this extension.
In addition, 
it is clear that for every~$h\in(0,1]$, it holds $\Delta_{\phi} e^{-\frac \phi h}=0$ in the distribution sense. Hence,
under Assumption~\ref{ass.confin},
since $\phi=\frac V2$ satisfies the relation   \eqref{eq:hypgenephi},
it holds~$e^{-\frac \phi h}\in L^{2}(\R^{d})$ and
the essential self-adjointness of $\Delta_{\phi}$ then implies that
$e^{-\frac \phi h} \in D(\Delta_{\phi}) $  so that $0\in \Ker \Delta_{\phi}$.
It follows moreover from \eqref{eq:hypgenephi} and from  \cite[Proposition~2.2]{HeKlNi04_01} that there exists $h_0>0$ and $c_0>0$ such that
for all $h\in(0,h_0]$, it holds
$$
\sigma_{ess}(\Delta_\phi)\ \subset\  [c_0,+\infty[.
$$
Coming back to the more general operator $P_{\phi}=P_{\phi,b,\nu}$ defined in \eqref{eq.P}, or equivalently to the operator  $L_{V,b,\nu}$
 according  to the relation \eqref{eq.unit},
the following proposition gathers some of its basic properties
which specify in particular  the preceding properties of $\Delta_{\phi}$ (and their equivalents concerning 
the weighted Laplacian $L_{V,0,0}$). It will be proven in Section~\ref{sub.Proof-prop}.
\begin{proposition}\label{prop:elemPphi} 
Under Assumption~\ref{hyp:statmes},  the operator $P_\phi$ with domain 
$\ccc_{c}^{\infty}(\R^{d})$
is accretive. Moreover, assuming in addition
Assumption~\ref{ass.confin},
there exists $h_0\in(0,1]$ such that
the following hold true for every $h\in(0,h_{0}]$:
\begin{enumerate}
\item[i)] The closure of $(P_{\phi},\ccc_{c}^{\infty}(\R^{d}))$, that we still denote by $P_{\phi}$, is maximal accretive, and hence its unique maximal accretive extension.
\item[ii)] The operator $P^{*}_\phi$ is maximal accretive and $\ccc_{c}^{\infty}(\R^{d})$
 is a core for $P^{*}_\phi$.   We have moreover the inclusions 
 $$D(\Delta_\phi)\subset D(P_\phi)\cap D(P_\phi^*)
\subset  D(P_\phi)\cup D(P_\phi^*)  \subset 
\{u\in L^2(\R^{d}),\; d_\phi u\in L^2(\R^{d})\}\,,
$$ where, for 
any unbounded operator $A$,
 $D(A)$ denotes the domain
of~$A$. In addition, for ${\bf P}_{\phi}\in \{P_{\phi},P^{*}_{\phi}\}$, we have the equality
$$
\forall\,u\,\in\, D({\bf P}_{\phi})\,,\ \ \ 
\Re\<{\bf P}_{\phi} u, u\>\ =\ \Vert d_\phi u\Vert^2\,.
$$
\item[iii)]  There exists  $\Lambda_0>0$ such that, defining
 $$\Gamma_{\Lambda_0}\ :=\ \big\{\Re(z)\geq 0\ \text{and}\ \vert\Im z\vert\leq\Lambda_0\max\big(\Re(z),\sqrt{\Re(z)}\big)\big\}\ \subset\ \C\,,$$
the spectrum $\sigma(P_\phi)$ of $P_{\phi}$ is included in   $\Gamma_{\Lambda_0}$ and
$$ 
 \forall\,z\,\in\, \Gamma_{\Lambda_0}^c\cap\{\Re(z)> 0\}\,,\ 
 \Vert (P_\phi-z)^{-1}\Vert_{L^2\rightarrow L^2}\ \leq\  \frac 1{\Re (z)}\,.
 $$

 \item[iv)] There exists $c_1>0$ such that  the map $z\mapsto (P_\phi-z)^{-1}$ is meromorphic in 
 $\{\Re(z)<c_1\}$ with finite rank residues. In particular, the spectrum of $P_{\phi}$
 in $\{\Re(z)<c_1\}$ is made of isolated eigenvalues with finite algebraic multiplicities.

\item[v)]  It holds $\Ker P_{\phi} = \Ker P^{*}_{\phi} =\sspan\{e^{-\frac \phi h}\}$ and 
$0$ is an isolated eigenvalue of $ P_{\phi}$ (and then of $P^{*}_{\phi}$) with algebraic multiplicity one.
 \end{enumerate}
 \end{proposition}

 From \eqref{eq.unit} and the last item of Proposition~\ref{prop:elemPphi},  note that
 $\Ker L_{V,b,\nu}=\sspan\{1\}$ and that
$0$ is an isolated eigenvalue of $ L_{V,b,\nu}$ with algebraic multiplicity one.
Moreover, according to Proposition~\ref{prop:elemPphi} and to the Hille-Yosida theorem, 
the operators $ L_{V,b,\nu}$ and its adjoint  $ L^{*}_{V,b,\nu}$ (in $L^{2}(\R^{d},m_{h})$)  generate, for every $h>0$ small enough, 
contraction  
semigroups $(e^{-t L_{V,b,\nu}})_{t\geq 0}$
and 
$(e^{-t L^{*}_{V,b,\nu}})_{t\geq 0}$
on $L^{2}(\R^{d},m_{h})$
 which permit to solve  \eqref{eq:KS}.

\subsection{Generic Morse-type hypotheses and labelling procedure}
\label{sub.label}

In order to describe precisely, in particular by stating Eyring-Kramers type formulas,  the spectrum around $0$ of 
$L_{V,b,\nu}$ (or equivalently of $P_{\phi}$) 
  in the regime $h\to 0$, we 
 will assume from now on that $V$ is a Morse function: 
 
 \begin{assumption}\label{hyp:phiMorse}
 The function $V$ is a Morse function.
 \end{assumption}

 Under Assumption~\ref{hyp:phiMorse} and thanks to Assumption~\eqref{eq:hypgenephi}, the set $\uuu$ made  of the critical points of $V$ is finite. In the following, 
  the critical points of $V$ with index $0$ and with index $1$, that is
  its local minima  and its saddle points, will play a fundamental role, 
  and we will respectively denote by $\uuu^{(0)}$ and $\uuu^{(1)}$
  the sets made of these points.  Throughout the paper, we will moreover denote 
  $$n_0\ :=\ {\rm card } (\uuu^{(0)})\,.$$
From the pioneer work by Witten \cite{Wi82}, it is well-known that for every $h\in(0,1]$ small enough, there is a correspondence between the  small eigenvalues of $\Delta_\phi$ and the local minima of $\phi=\frac V 2$.
More precisely, we have the following result (see in particular \cite{HeSj85_01,He88_01,HeKlNi04_01} or more recently \cite{MiZw17}).

\begin{proposition}\label{prop:localspectWitten}
Assume that \eqref{eq:hypgenephi} and Assumption \ref{hyp:phiMorse} hold true. Then, there exist  $\epsilon_0>0$ and $h_0>0$ such that for every $h \in(0,h_0]$, 
$\Delta_\phi$ has precisely~$n_0$ eigenvalues (counted with multiplicity) in the interval
 $[0,\epsilon_0 h]$.
  Moreover, these eigenvalues are actually exponentially small, that is 
live in an interval $[0,Che^{-2\frac Sh}]$ for some $C,S>0$ independent of $h\in(0,h_{0}]$.
\end{proposition}


Since the operator $P_\phi=\Delta_{\phi}+b_{h}\cdot d_{\phi}$ is  not self-adjoint (when $b_{h}\neq 0$), the analysis of its spectrum is  more complicated than the one of the spectrum of $\Delta_\phi$. 
The following result states a counterpart of Proposition~\ref{prop:localspectWitten}
in this setting. In this statement and in the sequel, for any $a\in\cz$ and $r>0$,
we will denote by $D(a,r)\subset \cz$ the open disk of center $a$ and radius $r$.  
%

\begin{theorem}\label{th:small-ev-Pphi}

Assume that Assumptions~\ref{hyp:statmes} to \ref{hyp:phiMorse} hold true, and let 
$\epsilon_0>0$ be given by Proposition \ref{prop:localspectWitten}. Then, for every $\epsilon_1\in(0,\epsilon_0)$, there exists $h_0>0$ such that for all $h\in(0,h_0]$, the set 
$
\sigma(P_\phi)\cap \{\Re z < \epsilon_1h \}
$
is finite and 
consists in 
$$n_{0}\ =\ {\rm card}\big(\,\sigma(\Delta_{\phi} )\,\cap\, \{\Re z < \epsilon_0h \}\, \big)$$ eigenvalues counted with algebraic multiplicity.
%
%
Moreover, there exists $C>0$ such that for all $h\in(0,h_{0}]$,
$$\sigma(P_\phi)\cap\{\Re z < \epsilon_1h \} \ \subset\ D(0, C h^{\frac 12} e^{-\frac Sh})\,,$$
where $S$ is given by Proposition \ref{prop:localspectWitten}. 
Eventually, for every $\epsilon\in(0,\epsilon_1)$, one has, uniformly with respect to $z$,
$$
\forall\,z\,\in\, \{\Re z < \epsilon_1h \}\cap \{|z|>\epsilon h\}\,,\ \ \Vert (P_\phi-z)^{-1}\Vert_{L^2\rightarrow L^2}\ =\ \ooo(h^{-1})\,.
$$
Lastly,   all the above conclusions also
hold  for $P^{*}_\phi$.
\end{theorem}

\begin{figure}
 \begin{center}
\scalebox{0.8}{ 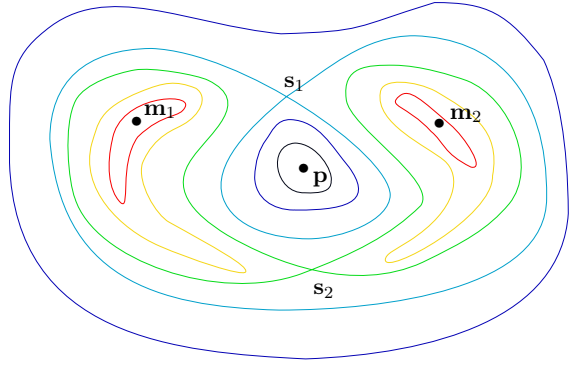}
  \end{center}
  \caption{Some level sets of a Morse function $V$ such that $V(x)\to +\infty$
  when $|x|\to+\infty$  and admitting five critical points: two local minima $\m_1$ and $\m_2$, one local maximum~${\bf p}$, and two saddle points
  $\s_1$ and $\s_2$. The point $\s_1$ is non-separating, whereas $\s_2$ is separating.}
  \label{fig:defssp}
  \end{figure}

This theorem will be proved in Section~\ref{sub-spec-origin} using Proposition~\ref{prop:localspectWitten} and a finite dimensional reduction. 
In order to give sharp asymptotics of the small eigenvalues of $P_\phi$, that is the ones in $D(0, C h^{\frac 12} e^{-\frac Sh})$, we will introduce some additional, but generic, topological assumptions on the Morse function~$V$ (see Assumption~\ref{hyp:gener} below). To this end, we first recall the general labelling of 
\cite{HeHiSj11_01} (see in particular Definition~4.1 there) generalizing the labelling of 
\cite{HeKlNi04_01} (and of \cite{BEGK,BGK}).
The main ingredient is the notion of separating saddle point, defined in Definition~\ref{def:SSP} below (see also an illustration in Figure~\ref{fig:defssp}) after the following observation.
Here and in the sequel, we define, for $a\in \mb R$,
$$\{V<a\}\ :=\     V^{-1} \big(  (-\infty, a) \big)
\quad\text{and}
\quad
\{V\leq a\}\ :=\     V^{-1} \big(  (-\infty, a] \big)\,,
$$
and $\{V>a\}$, $\{V\geq a\}$ in a similar way. The following lemma recalls the local structure of the sublevel sets of a Morse function. A proof can be found in \cite{HeKlNi04_01}.

\begin{lemma}
\label{le.SSP}
Let $z\in \R^d$ and $V:\R^d\to \R$ be a Morse function.
For any~$r>0$, we denote by $B(z,r)\subset \R^{d}$ the open ball of center
$z$ and radius $r$.
Then, for every $r>0$ small enough, $B(z,r)\cap \{V<V(z)\}$
has at least two connected components 
if and only if $z$ is a saddle point of $V$, i.e. if and only if
 $z\in \mathcal U^{(1)} $. In this case,
 $B(z,r)\cap \{V<V(z)\}$ has precisely two connected components.  
\end{lemma} 
\begin{definition}\label{def:SSP}
i) We say that  the saddle point $\s\in \uuu^{(1)}$ is a separating saddle point of $V$ if, for every $r>0$ small enough, the two connected components of $B(\s,r)\cap \{V<V(\s)\}$  (see Lemma~\ref{le.SSP})
are contained in different connected components of 
$\{V<V(\s)\}$. We will denote by $\vvv^{(1)}$ the set made of these points.\medskip

\noindent
ii) We say that $\sigma\in\R$ is a separating saddle value of $V$ if it has  the form 
$\sigma=V(\s)$ for some $\s\in\vvv^{(1)}$.\medskip

\noindent
iii) Moreover, we say that $E\subset \R^{d}$ is a critical component of $V$ if there exists 
$\sigma\in V(\mathcal V^{(1)})$ 
such that $E$ is a connected component of $\{V<\sigma\}$ satisfying
$\partial E\cap\vvv^{(1)}\neq \emptyset$. 
\end{definition}

Let us now describe the general labelling procedure of  \cite{HeHiSj11_01}. 
We will omit
details
when associating local minima and separating saddle points below, 
but  the following proposition (cf. \cite[Proposition~18]{DLLN})
can be helpful to well understand the construction.

\begin{proposition}
\label{pr.SSP}
Assume that $V$ is a Morse function with a finite number of critical points and   such that $V(x)\to +\infty$ when $|x|\to+\infty$.
Let $\lambda\in\R$ and let $\mathcal C$ be a connected component of $\{V<\lambda\}$. Then,
it holds
$$
\mathcal C \cap \vvv^{(1)} \neq \emptyset \ \ \text{iff}\ \
{\rm card}(\mathcal C\cap \mathcal U^{(0)})\geq 2\,.
$$
Let us also define
$$
\sigma\ :=\ \max_{\mathcal C\cap \vvv^{(1)}}V
$$
with the convention $\sigma:=\min_{\mathcal C}V $ when $\mathcal C\cap\vvv^{(1)}=\emptyset$.
It then holds:
\begin{enumerate}
\item[i)] For every $\mu\in(\sigma,\lambda]$, the set $\mathcal C\cap \{V<\mu\}$
is a connected component of $\{V<\mu\}$.
\item[ii)] If $\mathcal C\cap \vvv^{(1)}\neq\emptyset$, then
$\mathcal C\cap \mathcal U^{(0)}\subset \{V<\sigma\}$
and all the connected components of $\mathcal C\cap  \{V<\sigma\}$
are critical.
\end{enumerate}
\end{proposition}

Under the hypotheses of Proposition~\ref{pr.SSP}, $V(\vvv^{(1)})$ is finite. 
We moreover assume that $n_{0}\geq 2$, so that, under the hypotheses of Proposition~\ref{prop:elemPphi}  and of Theorem~\ref{th:small-ev-Pphi},
 $0$ is not the only exponentially small
eigenvalue of $P_{\phi}$ (or equivalently of $L_{V,b,\nu}$) and
$\vvv^{(1)}\neq \emptyset$ by Proposition~\ref{pr.SSP}.
We then denote the elements of $V(\vvv^{(1)})$ by $\sigma_2>\sigma_3>\ldots>\sigma_N$, where $N\geq 2$.  For convenience, we also introduce a fictive infinite saddle value $\sigma_1=+\infty$.
Starting from~$\sigma_1$, we will recursively associate to each $\sigma_i$ a finite family of local minima 
$(\m_{i,j})_j$ and a finite family of critical components $(E_{i,j})_j$ (see Definition~\ref{def:SSP}).\medskip

Let $N_{1}:=1$, $\underline{\m}=\m_{1,1}$ be a global minimum of $V$ (arbitrarily chosen if there are more than one),
and $E_{1,1}:=\R^d$. We now proceed in the following way:
\begin{enumerate}
\item[--] Let us denote, for some $N_{2}\geq 1$,  
by $E_{2,1},\dots,E_{2,N_{2}}$ the connected components of $\{V<\sigma_{2}\}$ 
which do not contain $\m_{1,1}$. 
They are all critical by the preceding proposition and we associate to each~$E_{2,j}$, where
$j\in\{1,\dots,N_{2}\}$, some global  minimum $\m_{2,j}$ of $V|_{E_{2,j}}$
(arbitrarily chosen if there are more than one).
\item[--] Let us then consider, for some $N_{3}\geq 1$,   the connected components
$E_{3,1},\dots,E_{3,N_{3}}$
of
$\{V<\sigma_{3}\}$
which do not contain the local minima of $V$ previously labelled. 
These components are also critical and included in
the $E_{2,j}\cap \{V<\sigma_{3}\}$'s, $j\in\{1,\dots,N_{2}\}$, such that $E_{2,j}\cap \{V=\sigma_{3}\}\cap \vvv^{(1)}\neq \emptyset$ (and $\sigma_{3}=\max_{E_{2,j}\cap\vvv^{(1)}}V $ for such a $j$). 
We then again associate to each $E_{3,j}$,
$j\in\{1,\dots,N_{3}\}$, some global minimum~$\m_{3,j}$ of $V|_{E_{3,j}}$. 
\item[--] We continue this process until having considered the connected components
of  $\{V<\sigma_{N}\}$, after which all the local minima of $V$ have been labelled.\medskip
\end{enumerate}

Next, we define two mappings 
$${ E}:\mathcal U^{(0)}\to \mathcal P(\R^{d}) \quad \text{and}\quad
{\bf j}: \mathcal U^{(0)}\to \mathcal P(\vvv^{(1)}\cup \{\s_{1}\})\,,$$
where, for any set $A$, $\mathcal P(A)$ denotes the power set of $A$, and
 $\s_{1}$ is a fictive saddle point such that $V(\s_{1})=\sigma_{1}=+\infty$,
as follows: for every $i\in\{1,\dots,N\}$ and $j\in\{1,\dots,N_{i}\}$,
\begin{equation}
\label{eq.E}
 { E}(\m_{i,j})\ :=\ E_{i,j}
\end{equation}
and
\begin{equation}
\label{eq.j}
  {\bf j}(\underline{\m})\ :=\ \{\s_{1}\}
\ \ \ \text{and, when $i\geq 2$,}\ \ \   
   {\bf j}(\m_{i,j})\ :=\ \partial E_{i,j}\cap\vvv^{(1)}\neq\emptyset\,.
\end{equation}
In particular, it holds $E(\underline{\m})=\R^d$ and
$$
\forall \,i\,\in\,\{1,\dots,N\}\,,\ \forall\,j\,\in\,\{1,\dots,N_{i}\}\,,\ \  \  \emptyset\neq{\bf j}(\m_{i,j})\ \subset
\ \{V=\sigma_{i}\}\,.
$$
Lastly, we define the mappings
 $$\bsigma:\uuu^{(0)}\rightarrow V(\vvv^{(1)})\cup\{\sigma_{1}\}
 \quad\text{and}\quad
S:\uuu^{(0)}\rightarrow (0,+\infty]$$ 
by
\begin{equation}
\label{eq.S}
\forall\,\m\in\uuu^{(0)}, \ \ \bsigma(\m):=V({\bf j}(\m)) \quad\text{and}\quad 
S(\m):=\bsigma(\m)- V(\m)\,,
\end{equation}
where, with a slight abuse of notation, we have identified the set $V({\bf j}(\m))$
with its unique element. Note that $S(\m)=+\infty$ if and only if $\m=\underline\m$.
An example of the preceding labelling is given in Figure~\ref{fig:association} below.\medskip

\begin{figure}[!h]
   \begin{center}
  \begin{tikzpicture}[scale=0.6]
\tikzstyle{vertex}=[draw,circle,fill=black,minimum size=5pt,inner sep=0pt]

\draw (-10,5) ..controls  (-7,-1.3).. (-5,1.25) ;
\draw (-5,1.25) ..controls  (-3.6,-0.9).. (-2,1 )  ;
\draw (-1.5,1 ) ..controls  (0,-1).. (1.9,2.6)  ;
\draw (5,5)  ..controls  (4,-3).. (2,2.6)  ;
\draw (-2,1 ) ..controls  (-1.75,1.33).. (-1.5,1)   ;

\draw (-1.74,1.2)  node[vertex,label=north: {\small{$\s_{3,2}$}}](v){}; 
\draw (-5,1.2)  node[vertex,label=north: {\small{$\s_{3,1}$}}](v){}; 
\draw (1.95,2.5)  node[vertex,label=north: {\small{$\s_2$}}](v){};
\draw (3.75,-1.35) node[vertex,label=south: {\small{$\m_{1,1}$}}](v){};
\draw (-3.5,-0.35) node[vertex,label=south: {\small{$\m_{3,1}$}}](v){};
\draw (-0.05,-0.35)  node[vertex,label=south: {\small{$\m_{3,2}$}}](v){};
\draw (-6.7,-0.35)  node[vertex,label=south: {\small{$\m_{2,1}$}}](v){};
 \draw[dashed, <->]  (-11,-2.5) -- (6,-2.5);
 \draw (-1.75,-2.2) node[]{\tiny{$ E_{1,1}=\R$}};
  \draw[dashed, <->]  (-8.7,2.5) -- (1.6,2.5);
 \draw (-3.5,2.8) node[]{\tiny{$ E_{2,1}$}};
   \draw [dashed, <->]    (-1.55,1.25)-- (1.15,1.25)   ;
 \draw (0,1.55) node[]{\tiny{$  E_{3,2}$}};
   \draw [dashed, <->]    (-4.8,1.25)-- (-2,1.25)   ;
 \draw (-3.5,1.55) node[]{\tiny{$  E_{3,1}$}};
 
  \draw  [<->]    (-13,-0.35)-- (-13,2.5)   ;
  \draw (-13.9,1) node[]{\tiny{$ S(\m_{2,1})$}};
  \draw  [densely dotted]    (-12.8,-0.35)-- (0,-0.35)   ;
  \draw  [densely dotted]    (-12.8,2.5)-- (-8.9,2.5)   ;
  
  \draw  [<->]    (-9.75,-0.35)-- (-9.75,1.25)   ;
  \draw (-10,0.5) node[]{\tiny{$ S(\m_{3,1})=S(\m_{3,2})$}};
  \draw  [densely dotted]    (-9.55,1.25)-- (-5,1.25)   ;
  
  \draw  [<-]    (6.65,-1.35)-- (6.65,5)   ;
  \draw (6.5,3) node[]{\tiny{$ S(\m_{1,1})=+\infty$}};
 \draw  [densely dotted]    (6.5,-1.35)-- (3.9,-1.35)   ;
\end{tikzpicture}

\end{center}
 \caption{A 1-D example of the preceding labelling when~$V$ admits four local minima. In this example, 
it holds $V(\m_{1,1})<V(\m_{2,1})=V(\m_{3,1})=V(\m_{3,2})$, 
  ${\bf j}(\m_{2,1})=\{\s_{2}\}$, ${\bf j}(\m_{3,1})=\{\s_{3,1},\s_{3,2}\}$, and ${\bf j}(\m_{3,2})=\{\s_{3,2}\}$.
  Note moreover that other choices of construction of the maps
  ${\bf j}$ and  $E$ are possible here since  ${\rm argmin}_{E_{2,1}}V=\{\m_{2,1},\m_{3,1},\m_{3,2}\}$.}
 \label{fig:association}
 \end{figure}

Our
generic topological assumption is the following one.
Assume that $V$ is a Morse function with a finite number
$n_{0}\geq 2$
 of critical points
   such that $V(x)\to +\infty$ when $|x|\to+\infty$,
and
let ${ E}:\mathcal U^{(0)}\to \mathcal P(\R^{d})$ and
${\bf j}: \mathcal U^{(0)}\to \mathcal P(\vvv^{(1)}\cup \{\s_{1}\})$
be the mappings defined in \eqref{eq.E} and in \eqref{eq.j}.
\begin{assumption}\label{hyp:gener} 
For every $\m\in \mathcal U^{(0)}$, the following hold true:
\begin{itemize}
\item[i)] the local minimum $\m$ is the unique global minimum of $ V|_{E(\m)}$,
\item[ii)] for all  $\m'\in\uuu^{(0)}\setminus\{\m\}$, 
${\bf j}(\m)\cap {\bf j}(\m')=\emptyset$.  
\end{itemize}
In particular, $V$ uniquely attains its global minimum, at
$\underline\m\in \mathcal U^{(0)}$.
\end{assumption}

\begin{figure}[!h]
   \begin{center}
  \begin{tikzpicture}[scale=0.6]
\tikzstyle{vertex}=[draw,circle,fill=black,minimum size=5pt,inner sep=0pt]

\draw (-10,5) ..controls  (-7,-1.3).. (-5,1.25) ;
\draw (-5,1.25) ..controls  (-3.6,-1.5).. (-2,1 )  ;
\draw (-1.5,1 ) ..controls  (0,-1).. (1.9,2.6)  ;
\draw (5,5)  ..controls  (4,-3).. (1.9,2.6)  ;
\draw (-2,1 ) ..controls  (-1.75,1.33).. (-1.5,1)   ;

\draw (-1.74,1.25)  node[vertex,label=north: {{\small $\s_{3,2}$}}](v){}; 
\draw (-5,1.25)  node[vertex,label=north: {{\small $\s_{3,1}$}}](v){}; 
\draw (1.88,2.6)  node[vertex,label=north: {{\small $\s_2$}}](v){};
\draw (3.75,-1.35) node[vertex,label=south: {{\small $\m_{1,1}$}}](v){};
\draw (-3.55,-0.8) node[vertex,label=south: {{\small $\m_{2,1}$}}](v){};
\draw (-0.05,-0.35)  node[vertex,label=south: {{\small $\m_{3,2}$}}](v){};
\draw (-6.7,-0.35)  node[vertex,label=south: {{\small $\m_{3,1}$}}](v){};
 \draw[dashed, <->]  (-11,-2.5) -- (6,-2.5);
 \draw (-1.75,-2.2) node[]{\tiny{$ E_{1,1}=\R$}};
  \draw[dashed, <->]  (-8.7,2.6) -- (1.6,2.6);
 \draw (-3.5,2.9) node[]{\tiny{$ E_{2,1}$}};

 \draw  [dashed, <->]    (2.1,2.6)-- (4.6,2.6)   ; 
\draw (3.3,2.9) node[]{\tiny{$\widehat E_{2}$}};   
   
    \draw [dashed, <->]    (-1.55,1.25)-- (1.15,1.25)   ;
 \draw (0,1.55) node[]{\tiny{$  E_{3,2}$}};
   \draw [dashed, <->]    (-8,1.25)-- (-5.3,1.25)   ;
 \draw (-6.5,1.55) node[]{\tiny{$  E_{3,1}$}};
 
 \draw  [dashed, <->]    (-4.8,1.25)-- (-1.95,1.25)   ; 
 \draw (-3.5,1.55) node[]{\tiny{$\widehat E_{3}$}};

  \draw  [<->]    (-13,-0.7)-- (-13,2.6)   ;
  \draw (-13.9,1) node[]{\tiny{$ S(\m_{2,1})$}};
  \draw  [densely dotted]    (-12.8,-0.8)-- (-3.8,-0.8)   ;
  \draw  [densely dotted]    (-12.8,2.6)-- (-9,2.6)   ;
  
  \draw  [<->]    (-9.75,-0.35)-- (-9.75,1.25)   ;
  \draw (-10,0.5) node[]{\tiny{$ S(\m_{3,1})=S(\m_{3,2})$}};
  \draw  [densely dotted]    (-9.55,-0.35)-- (0,-0.35)   ;
  \draw  [densely dotted]    (-9.55,1.25)-- (-8.2,1.25)   ;
 
  \draw  [<-]    (6.65,-1.35)-- (6.65,5)   ;
  \draw (6.5,3) node[]{\tiny{$ S(\m_{1,1})=+\infty$}};
 \draw  [densely dotted]    (6.5,-1.35)-- (3.9,-1.35)   ;

  \end{tikzpicture}

 \caption{A 1-D example 
 when~$V$ admits four local minima
 and satisfies Assumption~\ref{hyp:gener}. Here, $V(\m_{1,1})<V(\m_{2,1})<V(\m_{3,1})=V(\m_{3,2})$,
  ${\bf j}(\m_{2,1})=\{\s_{2}\}$, ${\bf j}(\m_{3,1})=\{\s_{3,1}\}$,
  and ${\bf j}(m_{3,2})=\{\s_{3,2}\}$. 
  Moreover, $\widehat E_{2}$ and
  $\widehat E_{3}$  denote respectively the sets $\widehat E(\m_{2,1})$
  and $\widehat E(\m_{3,1})=\widehat E(\m_{3,2})$
 introduced in Remark~\ref{rem:comm-2compo}.}
 \label{fig:association'}
 \end{center}
 \end{figure}
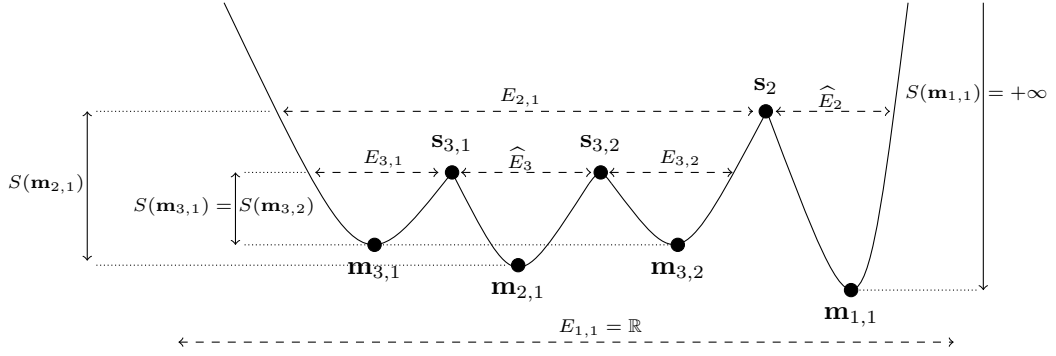

Note that the example
of Figure~\ref{fig:association} does not satisfy Assumption~\ref{hyp:gener} since
neither item i) nor ii) holds there. See also
Figure~\ref{fig:association'} below for a similar example satisfying  Assumption~\ref{hyp:gener}. 

Let us moreover underline that this assumption is  a little more general than the one
considered in the generic case in  \cite{HeKlNi04_01,HeHiSj11_01}
(see also \cite{BEGK,BGK}) where,  for instance,
 each set ${\bf j}(\m)$, $\m\in \uuu^{(0)}\setminus\{\m\}$, is assumed to only contain one element.

\begin{remark}
\label{rem:comm-2compo}
One can also show that Assumption~\ref{hyp:gener} 
implies that for every~$\m\in \uuu^{(0)}$ such that $\m\neq \underline\m$,
there is precisely one  connected component $\widehat E(\m)\neq E(\m)$ of $\{f<\bsigma(\m)\}$
such that $\overline {\widehat E(\m)}\cap \overline{E(\m)}\neq \emptyset$.
In other words, there exists a  connected component $\widehat E(\m)\neq E(\m)$ of $\{f<\bsigma(\m)\}$
such that ${\bf j}(\m)\subset \pa\widehat E(\m)$. 
Moreover, the global minimum $\m'$ of $V|_{\widehat E(\m)}$ is unique  and satisfies
$\bsigma(\m')>\bsigma(\m)$ and
$V(\m')<V(\m)$ (see examples of such sets in Figure~\ref{fig:association'}). We refer to \cite{Mi19} or \cite{DLLN} for more details on the geometry of the sublevel sets
of a Morse function.
\end{remark}

\subsection{Main results and comments} In order to state our main results, we also
need
the following lemma which  is fundamental in our analysis.
\begin{lemma}
\label{le.B}~\ For $x\in\R^{d}$, let $B(x):={\rm Jac}_{x}b$ denote the Jacobian matrix of $b=U_0-\nabla V$ at $x$,
and consider a saddle point $\s\in\mathcal U^{(1)}$.
\begin{enumerate}
\item[i)] 
The  matrix  $\Hess V(\s)+B^{*}(\s)\in\mathcal M_{d}(\R)$ 
admits precisely one negative eigenvalue $\mu=\mu(\s)$,
which has moreover geometric multiplicity one.
\item[ii)]  Denote by $\xi=\xi(\s)$ one of the two (real) unitary eigenvectors of $\Hess V(\s)+B^{*}(\s)$
associated with $\mu$.
The real symmetric matrix
$$M_{V}\ :=\ \Hess V(\s)+2|\mu|\,\xi\,\xi^{*}$$
is then positive definite and its determinant satisfies:
$$
\det M_{V}\ =\ -\det \Hess V(\s)\,.
$$
\item[iii)] Lastly, denoting by $\lambda_{1}=\lambda_{1}(\s)$ the negative eigenvalue of $\Hess V(\s)$, it holds $|\mu|\geq|\lambda_{1}|$, with equality if and only if $B^{*}(\s)\xi=0$, and
$$
\lp (\Hess V(\s))^{-1}\xi,\xi\rp \ =\ \frac1\mu\ <\ 0\,.
$$ 
\end{enumerate}
\end{lemma}

Note that the real matrix $\Hess V(\s)+B^{*}(\s)$  of Lemma~\ref{le.B}
is in general non symmetric. Let us also point out that
the statements of Lemma~\ref{le.B} already 
appeared in  the related work  \cite{LaMaSe19} 
(see in particular the beginning of Section~8 there),
and in   \cite{LaSe}, where proofs are given (see indeed Section~4.1 there).
We will nevertheless give a proof in Section~\ref{sec:preuvelemalg}
for the sake of completeness.\medskip

We can now state our main results.

\begin{theorem}\label{th.main}
Suppose that  Assumptions \ref{hyp:statmes} to \ref{hyp:gener} hold true, and let 
$\epsilon_0>0$ be given by Proposition \ref{prop:localspectWitten}.
Then, for all $\epsilon_1\in(0,\epsilon_0)$, there exists $h_0>0$ such that for all $h\in(0,h_0]$, one has,
counting the eigenvalues with algebraic multiplicity,
$$
\sigma (L_{  V,b,\nu})\,\cap\,\{{\rm Re\,}z<\epsilon_1\}\ =\ \{\lambda(\m,h),\;\m\in\uuu^{(0)}\},
$$
where, denoting by $\underline\m$ the unique absolute minimum of $V$, $\lambda(\underline{\m},h)=0$ and, for all $\m\neq\underline{\m}$, $\lambda(\m,h)$
satisfies the following Eyring-Kramers type formula:
  \begin{align}
  \label{eq.vp}
\lambda( \m, h)   &\  = \ 
\zeta(\m)    \,  e^{-\frac{S(\m)}h} 
\,\big(1+\mathcal O(\sqrt h)\big)\,,
\end{align}
where $S:\uuu^{(0)}\rightarrow (0,+\infty]$ is defined in \eqref{eq.S}
and, for every  $\m\in\uuu^{(0)}\setminus\{\underline\m\}$,
\begin{equation}
\label{eq.pref}
\zeta(\m)\ :=\ \frac{\det \Hess V(\m)^{\frac12}}{2\pi} \Big(\sum_{\s\in{\bf j}(\m)}\frac{ |\mu(\s)|}{|\det \Hess V(\s)|^{\frac12}}\Big)\,,
\end{equation}
where ${\bf j}:\mathcal U^{(0)}\to \mathcal P(\vvv^{(1)}\cup \{\s_{1}\})$  is defined in \eqref{eq.j}
and the $\mu(\s)$'s  are defined in Lemma~\ref{le.B}.\\
In addition, it holds
$$
\sigma (L_{  V,-b,-\nu})\,\cap\,\{{\rm Re\,}z<\epsilon_1\} = 
\sigma (L^{*}_{  V,b,\nu})\,\cap\,\{{\rm Re\,}z<\epsilon_1\} = 
 \{\overline{\lambda(\m,h)},\;\m\in\uuu^{(0)}\}.
$$
\end{theorem}

\begin{remark}
\label{re.double-well}
In the case where $V$ has precisely two minima $\underline\m$  and $\m$ such that
$V(\underline\m)=V(\m)$, the above result can be easily  generalized.
 In this case, using the definitions of $S$ and ${\bf j}$ given in~\eqref{eq.S} and in~\eqref{eq.j}
 (note that the choice of $\underline\m$ among the two minima of $V$ is arbitrary in this case),
 we have, counting the eigenvalues with algebraic multiplicity, for every $h>0$ small enough,
$$
\sigma(L_{  V,b,\nu})\,\cap\,\{{\rm Re\,}z<\epsilon_1\}\ =\ \{0,\lambda(\m,h)\}\,,
$$
where 
  \begin{equation*}
\lambda( \m,h)   \  =    \  
\zeta(\m) \,  e^{-\frac{S(\m)}{h}} 
\,\big(1+\mathcal O(\sqrt h)\big)
\end{equation*}
with
$$
\zeta(\m)\ =\   \frac{\det \Hess V(\m)^{\frac12}+\det \Hess V(\underline\m)^{\frac12}}{2\pi} \Big(\sum_{\s\in{\bf j}(\m)}\frac{ |\mu(\s)|}{|\det \Hess V(\s)|^{\frac12}}\Big)\,.
$$
Moreover, since $\sigma(L_{  V,b,\nu})=\overline{\sigma(L_{  V,b,\nu})}$, the eigenvalue $\lambda(\m,h)$ is real.
\end{remark}

Let us make a few comments on the above theorem.\medskip 

First, observe that if we assume that $U_{h}=\nabla V$, that is if $b_h=0$ (see \eqref{eq.bh}), we obtain the precise asymptotics of the small eigenvalues of $L_{V,0,0}$
 (or equivalently of $\Delta_{\phi}$ after multiplication by $\frac1h$, see \eqref{eq.unit}) 
 and hence recover the results already proved
 in this reversible setting in
 \cite{BGK,HeKlNi04_01}
(see also \cite{MeSc} for an extension to logarithmic Sobolev inequalities). 
In this case, for every saddle point $\s$ appearing in \eqref{eq.pref},
the real number $\mu(\s)$ is indeed the negative eigenvalue of $\Hess V(\s)$ according
to the  first item of Lemma~\ref{le.B}.
Let us also point out that under the hypotheses made in \cite{BGK,HeKlNi04_01}, 
the set ${\bf j}(\m)$ actually contains one unique element for every $\m\in\uuu^{(0)}\setminus\{\underline\m\}$.
Moreover, our analysis permits in this case to recover that
the error term $\ooo(\sqrt h)$ is actually of order $\ooo(h)$, as proven in \cite{HeKlNi04_01}. However, it does not permit
to prove that this $\ooo(h)$ actually admits a full asymptotic expansion in $h$
as proven in \cite{HeKlNi04_01}.\medskip

To the best of our knowledge, the above theorem is the first result giving sharp asymptotics of the small eigenvalues of the generator 
 $L_{V,b,\nu}$ in the non-reversible case. 
 Similar results were obtained by H\'erau-Hitrik-Sj\"ostrand for the Kramers-Fokker-Planck (KFP) equation in \cite{HeHiSj11_01}. Compared to our framework, they deal with non-self-adjoint and non-elliptic operators, which makes the analysis more complicated. However, the KFP equation enjoys several symmetries which are crucial in their analysis. First of all, the KFP operator has a 
 supersymmetric structure (for a non-symmetric skew-product $\<.,.\>_{\text{KFP}}$) which  permits to write the interaction matrix associated with the small eigenvalues as a square   $M=A^{*}A$, where the adjoint $A^*$ is taken with respect to  $\<.,.\>_{\text{KFP}}$. 
 Using this square structure, the authors can then follow the strategy of \cite{HeKlNi04_01} to construct accurate approximations of the matrices  $A$ and $A^*$. However, since $\<.,.\>_{\text{KFP}}$ is not a scalar product, they cannot identify the squares of the singular values of $A$ 
 with the eigenvalues of $M$.  This difficulty  is solved by using an extra symmetry (the PT-symmetry), which permits to modify the skew-product 
 $\<.,.\>_{KFP}$ into a new product $\<.,.\>_{KFPS}$, which is a scalar product when restricted 
 to the  ``small spectral subspace", and  for which the identity $M=A^*A$ remains true with an adjoint taken with respect to $\<.,.\>_{KFPS}$.
This permits  to conclude as in  \cite{HeKlNi04_01},  using in particular the Fan inequalities to estimate the singular values of $A$.
 
  In the present case, none of these two symmetries are available in general ($L_{V,b,\nu}$,
  or equivalently $P_\phi$, enjoys  however a supersymmetric structure when 
  $b$ and $\nu$ satisfy the relation \eqref{eq.J-supersym}, see indeed \eqref{eq.J-supersym'}
  or
  Remark \ref{rem:susy} below in this connection).
  We then developed an alternative approach based on the construction of very accurate quasimodes
  and partly inspired by   \cite{DiLe} (see also the related constructions made in  \cite{BEGK,LaMaSe19,LeNe}). 
  This permits  the construction of the interaction matrix $M$ as above.
  However, since we cannot write $M=A^{*}A$ and use the Fan inequalities
  as in \cite{HeKlNi04_01,HeHiSj11_01} (and e.g. in \cite{HeNi06,Lep10,Mi19,DLLN,LeNe}),
 we have to compute  directly the  eigenvalues of $M$. 
 To this end, we use  crucially  the Schur complement method. 
 This leads to Theorem~\ref{th:specgradedmat} in appendix, which permits to replace
 the use of the Fan inequalities  to perform the final analysis in our setting. We believe that these two arguments are quite general and may be used in other contexts.\medskip

Though it is generic, one may ask if Assumption~\ref{hyp:gener} is necessary to get Eyring-Kramers type formulas as in 
Theorem~\ref{th.main}. In the reversible setting, the full general (Morse) case was recently  treated by the second author in \cite{Mi19}, but applying the methods developed there to our non-reversible setting was not straightforward and we decided to postpone this analysis to future works.
 Let us point out in this connection
 that in the general (Morse) case,
 some tunneling effect
 between the characteristic wells of $V$ defined by the mapping $E$ (see \eqref{eq.E})  mixes
 their corresponding prefactors,
 see  indeed Remark~\ref{re.double-well}, or  \cite{Mi19} for more intricate situations in the reversible setting.\medskip

Note that Theorem~\ref{th.main} does not state that the operator $L_{  V,b,\nu}$ is diagonalizable when restricted to the spectral subspace associated with its small eigenvalues. Indeed, since $L_{  V,b,\nu}$ is not self-adjoint, we cannot exclude the existence of Jordan's blocks.
We cannot neither  exclude the existence of  non-real eigenvalues, but 
the spectrum of $L_{  V,b,\nu}$ is obviously  stable by complex conjugation  since $L_{  V,b,\nu}$ is a partial differential operator with real coefficients.
 However, 
 in the case where for every $\m\in\uuu^{(0)}\setminus\{\underline\m\}$,
 the prefactors $\zeta(\m')$ defined in \eqref{eq.pref}
  are all distinct for $\m'\in S^{-1}(S(\m))$, the $\lambda(\m,h)$'s, $\m\in\uuu^{(0)}$,  are then real eigenvalues of multiplicity one of
 $L_{  V,b,\nu}$,
   and its restriction to its small spectral subspace is diagonalizable.\medskip

Coming back to the contraction semigroups
$(e^{-t L_{V,b,\nu}})_{t\geq 0}$
and 
$(e^{-t L^{*}_{V,b,\nu}})_{t\geq 0}$
on $L^{2}(\R^{d},m_{h})$
 introduced just after
Proposition~\ref{prop:elemPphi}, 
Theorem~\ref{th.main} has the following consequences
on the rate of convergence to equilibrium for the process~\eqref{eq:overdampedL}.
 
 \begin{theorem}
 \label{co.main}
Assume that the hypotheses of Theorem~\ref{th.main}  hold
and let
$\m^{*}\in \uuu^{(0)}\setminus\{\underline\m\}$ be such that
\begin{equation}
\label{eq.m*}
S(\m^{*})\ =\ \max_{\m\in\uuu^{(0)}} S(\m)
\quad\text{and}\quad  \zeta(\m^{*})\ =\ \min_{\m\in S^{-1}(S(\m^{*}))}\zeta(\m)\,,
\end{equation}
where the prefactors $\zeta(\m)$'s, $\m \in  \uuu^{(0)}\setminus\{\underline\m\}$, are defined
 in \eqref{eq.pref}, and $S:\uuu^{(0)}\rightarrow (0,+\infty]$ is defined in \eqref{eq.S}.
 Let us then define, 
 for any $h>0$,
 $$
 \lambda(h)\ :=\ \zeta(\m^{*})\,e^{-\frac{S(\m^{*})}h}\,.
 $$
Then, there exist $h_{0}>0$ and $C>0$
such that for every $h\in(0,h_{0}]$,
it   holds
\begin{equation}
\label{eq.CV1}
 \forall\,t\,\geq\,0\,,\ \ \ \|\,e^{-tL_{V,b,\nu}} - \Pi_{0}\,\|_{L^2(m_h)\rightarrow L^2(m_h)}\ 
 \leq\ C\,e^{-\lambda(h)(1-C\sqrt h)t}\,,
\end{equation}
 where
$\Pi_{0}$ denotes the orthogonal projector on $\Ker L_{V,b,\nu} =\sspan\{1\}$:
 $$
 \forall\,u\,\in\, L^{2}(m_{h})\,,\ \ \  \Pi_{0} u \ =\  \lp u,1 \rp_{L^{2}(m_{h})}   \ =\   \int_{\R^{d}} u \,dm_{h}  \,.
 $$ 
Assume moreover that $(X_{t})_{t\geq 0}$ is solution to \eqref{eq:overdampedL}
and that the probability distribution $\varrho_{0}$ of $X_{0}$ admits a density $\mu_{0}\in L^{2}(\R^{d},m_{h})$ with respect to the probability measure $m_{h}$.
Then, for every $t\geq 0$, the probability distribution~$\varrho_{t}$ of  $X_{t}$ admits   the density $\mu_{t}=e^{-tL^{*}_{V,b,\nu}}\mu_{0}\in L^{2}(\R^{d},m_{h})$ with respect to $m_{h}$, and for every $h\in(0,h_{0}]$,
it holds
\begin{equation}
\label{eq.CV2}
 \forall\,t\,\geq\,0\,,\ \ \ \|\,\varrho_{t}-\nu_{h}\,\|_{TV}
 \ \leq\  C\,\|\mu_{0}-1\|_{L^{2}(m_{h})}\,e^{-\lambda(h)(1-C\sqrt h)t}\,,
\end{equation}
where $\|\cdot \|_{TV}$ denotes the total variation distance.\\
Finally, when there exists one unique $\m^{*}$
satisfying \eqref{eq.m*}, the eigenvalue
$\lambda(\m^{*},h)$ associated with $\m^{*}$ (see \eqref{eq.vp}) is real and simple,
 and 
 the estimates \eqref{eq.CV1} and \eqref{eq.CV2}
 remain valid if one replaces $\lambda(h)(1-C\sqrt h)$ by $\lambda(\m^{*},h)$
 in the exponential terms.
\end{theorem}

Theorems~\ref{th.main} and \ref{co.main} 
describe the metastable behaviour of the dynamics~\eqref{eq:overdampedL}  from
a spectral perspective.\medskip 

Concerning the question of accelerating the convergence 
to equilibrium mentioned at the end of Section~\ref{sub.motivation},
 the exponential rate of convergence to equilibrium 
appearing in
the estimates \eqref{eq.CV1} and \eqref{eq.CV2}
 is generically strictly larger than the optimal  rate for the associated
gradient dynamics~\eqref{eq.ODL}.
To be more precise, let us assume, as in  the last part of   the statement of
Theorem~\ref{co.main}, that   there exists one unique $\m^{*}$
satisfying \eqref{eq.m*}.
The exponential rate of return to equilibrium appearing in 
\eqref{eq.CV1} and \eqref{eq.CV2} is then given  
by the spectral gap $\lambda(\m^{*},h)$ of $L_{V,b,\nu}$. Moreover, denoting by~$\lambda^{\!\nabla}(\m^{*},h)$
the spectral gap of the  generator $L_{V,0,0}$ of the associated gradient dynamics \eqref{eq.ODL},
that is
the optimal rate of return to equilibrium in the gradient setting,
it follows from Theorem~\ref{th.main} and item~iii) in Lemma~\ref{le.B} that, as soon as $B^{*}(\s^{*})\neq 0$ for at least one $\s^{*}\in{\bf j}(\m^*)$, the ratio of the rates~$\frac{\lambda(\m^{*},h)}{\lambda^{\!\nabla}(\m^{*},h)}$ converges to some constant $c>1$
when $h\to 0$.

 In addition, it is not difficult to see that playing with $b_{h}$,
one can make $\lim\limits_{h\to0}\frac{\lambda(\m^{*},h)}{\lambda^{\!\nabla}(\m^{*},h)}$ arbitrarily big. 
Taking for example
$b_{h}=b=aJ\,\nabla V$ around~$\s^{*}$
for $a\in\R$ and some  constant  antisymmetric and invertible matrix $J$, 
it holds
 $$\lim\limits_{a\to\infty}\lim\limits_{h\to0}\frac{\lambda(\m^{*},h)}{\lambda^{\!\nabla}(\m^{*},h)}=+\infty\,.$$
 
Nevertheless, making this limit too big will deteriorate the constant $C$ appearing in 
\eqref{eq.CV1} and \eqref{eq.CV2}, as well as the interval $(0,h_{0}] \ni h$ for which these estimates
remain relevant. A more interesting problem is the computation of the optimal rate when $h_0$ is small but fixed, that is when the preceding~$J$ 
has a constant size (see \cite{LeNiPa} in the case of linear drifts). We did not make the whole computation, but a partial one seems to indicate that the 
optimal (or at least a reasonable) choice for $J$ is given when it sends the unstable direction of $\Hess V( \s^* )$ onto one of its stable directions
corresponding to a maximal eigenvalue.
\medskip

A closely related point of view to ours is to study the
 mean transition times
between the different wells of the potential $V$ 
 for the process $(X_t)_{t\geq 0}$ solution to \eqref{eq:overdampedL}. 
In the non-reversible case,  this question  
 has been studied recently e.g. in \cite{BoRe16,LaMaSe19}, to which we also refer for more details
 and references on this subject. 
 
 In \cite{BoRe16}, an Eyring-Kramers  type
 formula (for the mean transition times)
  is derived from  formal computations relying on the study of the appropriate quasi-potential. 
  In the case of a  double-well  potential $V$ and under the assumption
  that $U_{h}=\nabla V+b$ (that is that $\nu=0$, see \eqref{eq.bh})
 for some vector field $b$ only satisfying 
 $b\cdot\nabla V=0$ (that is without assuming $\div b=0$ as we do when $\nu=0$, see~\eqref{eq:propvectfield}),
 the authors derived formula (5.65),  where, in comparison with
 a formula such as (the inverse of) \eqref{eq.vp} in Theorem~\ref{th.main},   appears in the prefactor some extra term 
measuring the non-Gibbsianness of their situation.
In this general setting, 
the measure $m_{h}$ is indeed invariant for the dynamics if and only if $\div b=0$,
and this
extra term 
involves the integral of the  function $F:= \div( b)$ along 
the so-called instanton trajectory.
Under the additional assumption that $m_{h}$ is  invariant (that is that $F=0$), 
this extra term equals~$1$, which leads to
the formula (5.66) in  \cite{BoRe16},
which is similar to (the inverse of)  \eqref{eq.vp} in Theorem~\ref{th.main} (see more precisely Corollary~\ref{co.Landim} below, which clarifies the relation between 
eigenvalues of $L_{V,b,\nu}$ and mean transition  times). In the present paper, we restrict ourselves to the Gibbsian case, so that our formulas do not contain any extra prefactor as discussed above. It would be of great interest to study the general case
of a drift of the form $\nabla V +b$, where $b\cdot\nabla V=0$ but without assuming $\div b=0$,   by mixing our approach and quasi-potential constructions.
 
  In  \cite{LaMaSe19}, the authors use a potential theoretic approach to prove  an Eyring-Kramers 
 type formula similar to the formula (5.66) of \cite{BoRe16}  
 in the case of a  double-well  potential $V$, when $b$ and $\nu$ satisfy
 the relation \eqref{eq.J-supersym} in such a way that $L_{V,b,\nu}$ has  the form~\eqref{eq.J-supersym'}. Though the mathematical objects considered in \cite{LaMaSe19} and in the present paper are not the same, these two works 
 share some similarities. Nevertheless, we would like to emphasize that our approach permits to go beyond the supersymmetric assumption \eqref{eq.J-supersym} and to treat the case of multiple-well potentials.\medskip

To be more precise on the connections between the present paper and \cite{LaMaSe19} (and also \cite{BoRe16}), let us conclude this introduction with the 
corollary below which combines the results given by Theorem~\ref{th.main} when $V$ is a double-well  potential
 and  \cite[Theorem~5.2 and Remarks~5.3 and~5.6]{LaMaSe19}.
 This  result generalizes in particular, in this non-reversible double-well setting, the results obtained in the reversible case
 in \cite{BEGK,BGK} on the relations between the small eigenvalues of $L_{V,b,\nu}$
and the mean transition times of \eqref{eq:overdampedL} when $b=\nu=0$.

\begin{corollary}
 \label{co.Landim} 
 Assume that the hypotheses of Theorem~\ref{th.main}  hold with moreover
 $$
 \lim_{|x|\to+\infty} \frac{x}{|x|}\cdot \nabla V(x)\ =\ +\infty\quad\text{and}\quad
 \lim_{|x|\to+\infty}  |\nabla V(x)|-2\Delta V(x)\ =\ +\infty\,,
 $$ 
 and that
 $V$ admits precisely two local minima $\underline\m$ and $\m$
 such that $V(\underline\m)<V(\m)$ (it then holds $\vvv^{(1)}={\bf j}(\m)$).
  Assume in addition
 that $b$ and $\nu$ satisfy the relation \eqref{eq.J-supersym},
and hence that $b=J\,\nabla V$ for some  smooth map $J$  from $\R^{d}$
into the set of real antisymmetric matrices of size $d$,
and that $J$ is uniformly bounded on $\R^{d}$.\\
Let $\mc O(\underline\m)$ be a smooth open connected set
containing $\underline\m$ such that $\overline{\mc O(\underline\m)}\subset\{V<\bsigma(\m)\}$.  
Let then $(X_{t})_{t\geq 0}$ be the  solution to \eqref{eq:overdampedL}
such that $X_{0}=\m$
and let
$$
\tau_{\mc O(\underline\m)}\ :=\ \inf\{t\geq0\,,\ X_{t}\in \mc O(\underline\m) \}
$$
 be the first hitting time of $\mc O(\underline\m)$. 
 The expectation of $\tau_{\mc O(\underline\m)}$ and 
  the non-zero small eigenvalue $\lambda(\m,h)$ of $L_{V,b,\nu}$ 
 are then related by the following formula in the limit 
$h\to 0$: 
 $$
\mathbb{E}(\tau_{\mc O(\underline\m)})\ =\ \frac{1}{\lambda(\m,h)}\,\Big(1+\ooo\big(\sqrt{h |\ln h|^{3}}\big)\Big)\,. $$ 
\end{corollary} 

Let us mention here that the hypotheses of Corollary~\ref{co.Landim} 
are simply the minimal hypotheses  permitting to apply
at the same time
Theorem~\ref{th.main} and 
\cite[Theorem~5.2]{LaMaSe19} in its refinement specified in
\cite[Remark~5.6]{LaMaSe19}.\\

\noindent 
\textbf{Acknowledgements.}  The authors thank the anonymous referees for their remarks who permitted to improve the quality of the paper. Both authors are members of the ANR project QuAMProcs 19-CE40-0010-01.\\

\section{General spectral estimates}

\subsection{Proof of Proposition~\ref{prop:elemPphi}}
\label{sub.Proof-prop}

\noindent
The unbounded operator $(P_\phi,\ccc^\infty_c(\R^d))$ is accretive, since, according to 
\eqref{eq:minorLV}, one has: 
\be\label{eq:accretPphi}
\forall\,u\,\in \,\ccc^\infty_c(\R^d)\,,\ \ \Re\<P_\phi u, u\>=\<\Delta_\phi u, u\>=\Vert d_\phi u\Vert^2\geq 0\,.
\ee
In order to prove that its closure  is maximal accretive, it then suffices to show
that $\Ran(P_\phi+1)$ is dense in $L^2(\R^{d})$   
(see for example \cite[Theorem~13.14]{Hel-spectral}). The proof of this fact is rather standard 
but we give it for the sake of completeness (see in particular   the proof 
of 
 \cite[Proposition~5.5]{HeNi05} for a similar proof). Suppose that 
$f\in L^{2}(\R^{d})$ is orthogonal to $\Ran(P_\phi+1)$. 
It then holds $(P_\phi^*+1)f=0$ in the distribution sense and,
since $P_\phi$ is real, one can assume that $f$ is real. In particular, since $P_\phi^*=\Delta_\phi-b_h\cdot d_\phi$ is elliptic with smooth coefficients, 
$f$ belongs to $\ccc^\infty(\R^{d})$. Thus, for every $\zeta \in \ccc^\infty_c(\R^d,\R)$, one has
\begin{align*}
h^{2}\<\nabla (\zeta f),\nabla(\zeta f)\>
+ \int \zeta^{2}(|\nabla \phi|^{2}-&h\Delta \phi+1)f^{2}  = 
\<(P_\phi^*+1)\zeta f,\zeta f\> \\
& = h^{2}\!\!\!\int |\nabla\zeta|^{2}f^{2}-h\int (b_{h}\cdot d\zeta)\zeta f^{2} .
\end{align*}
Take now $\zeta $ such that $0\leq\zeta\leq 1$, $\zeta =1$ on  $B(0,1)$ and $\supp\zeta \subset B(0,2)$,
and define $\zeta_{k}:=\zeta(\frac\cdot k)$ for $k\in\N^{*}$. 
According to \eqref{hyp:controlechamp} and to the above relation, there exists $C>0$ such that
 for every $k\in \N^{*}$, it holds
\begin{align*}
 \int \zeta_{k}^{2}(|\nabla \phi|^{2}-h\Delta \phi+1)f^{2} 
&\ \leq\ C\frac{h^{2}}{k^{2}} \|f\|^{2} + C\frac hk \|f\|\,\|(1+|\nabla \phi|) \zeta_{k} f \|\\
&\ \leq\ C(1+\frac{1}{2\varepsilon})\frac{h^{2}}{k^{2}} \|f\|^{2} +  \frac\varepsilon2 C\|(1+|\nabla \phi|) \zeta_{k} f \|^{2}\,,
\end{align*}
where $\varepsilon>0$ is arbitrary. Choosing $\varepsilon=\frac{1}{2C}$ and using \eqref{eq:hypgenephi}, it follows that
for every $h>0$ small enough, it holds
$$
\frac14\|\zeta_{k}f\|^{2} \ \leq\ \int \zeta_{k}^{2}(\frac12|\nabla \phi|^{2}-h\Delta \phi+\frac12)f^{2} \ \leq\ \frac 43C(1+\frac{1}{2\varepsilon})\frac{h^{2}}{k^{2}} \|f\|^{2}\,,
$$ 
which implies, taking the limit $k\to+\infty$,
 that $f=0$. Hence,
the closure of
$P_\phi$, that we still denote by $P_{\phi}$, is maximal accretive.
Note moreover, that \eqref{eq:accretPphi} implies that
$
D(P_\phi) \subset \{u\in L^2(\R^{d}),\; d_\phi u\in L^2(\R^{d})\}$
and that
$\Re\<P_\phi u, u\>=\Vert d_\phi u\Vert^2$
for every $u\in D(P_\phi)$.\medskip

%
%

\noindent
Let us now prove that $D(\Delta_\phi)\subset D(P_{\phi})$,
which amounts to show
that for every $u\in D(\Delta_\phi)$, there exists a sequence $(u_{n})_{n\in\N}$
of $\ccc^\infty_c(\R^d)$ such that $u_{n}\to u$ in $L^{2}(\R^{d})$ and $(P_{\phi}u_{n})_{n\in\N}$
is a Cauchy sequence. Since $(\Delta_\phi,\ccc^\infty_c(\R^d))$ is essentially self-adjoint,
for any such $u$, there exists a sequence $(u_{n})_{n\in\N}$
in $\ccc^\infty_c(\R^d)$ such that $u_{n}\to u$ in $L^{2}(\R^{d})$ and $(\Delta_{\phi}u_{n})_{n\in\N}$
is a Cauchy sequence, and it thus suffices to show that $(b_h\cdot d_\phi u_{n})_{n\in\N}$
is also a Cauchy sequence. For this purpose, we introduce the exterior derivative $d$ acting from $0$-forms into $1$-forms and the 
twisted semiclassical derivative $d_\phi=e^{-\phi/h}\circ hd\circ e^{\phi/h}$. Note that 
the notation $d_\phi$ has actually already been defined  in \eqref{eq.d-phi} with a different meaning;
we are thus making here
a slight abuse of notation, by identifying the exterior derivative~$d$ acting on functions with $\nabla$.
Thanks to \eqref{eq:hypgenephi} and to \eqref{hyp:controlechamp}, there exists $C>0$ such that
for every $h>0$ small enough and every $u\in \ccc^\infty_c(\R^d)$, one has  
\begin{equation*}
\begin{split}
\Vert b_h\cdot d_\phi u\Vert^2\leq\int\vert b_h\vert^2\vert d_\phi u\vert^2
&\leq C\<\vert\nabla\phi\vert^2 d_\phi u,d_\phi u\>+C\Vert d_\phi u\Vert^2\\
&\leq 2C\< \Delta_\phi^{(1)}\,d_\phi u,d_\phi u\>+2C\Vert d_\phi u\Vert^2\,,
\end{split}
\end{equation*}
where  $\Delta_\phi^{(1)}$ denotes the Witten Laplacian acting on $1$-forms, that is
$$
\Delta_\phi^{(1)}=\Delta_\phi^{(0)}\otimes \Id +2h\Hess\phi
=(-h^2\Delta+\vert\nabla\phi\vert^2-h\Delta\phi)\otimes \Id + 2h\Hess \phi.
$$
Combined  with the intertwining relation $\Delta_\phi^{(1)}d_\phi=d_\phi \Delta_\phi^{(0)}$, we get
\be\label{eq:estimbdphi}
\Vert b_h\cdot d_\phi u\Vert^2
\ \leq\ 2C\big(\Vert \Delta_\phi^{(0)} u\Vert^2
+\Vert d_\phi u\Vert^2\big)
\ \leq\ 2C\Vert \Delta_\phi^{(0)} u\Vert\big(\Vert \Delta_\phi^{(0)} u\Vert+\Vert u\Vert\big)
\ee
for every $u\in \ccc^\infty_c(\R^d)$. This implies that for any Cauchy sequence  
$(u_{n})_{n\in\N}$ in $L^{2}(\R^{d})$ such that $(\Delta_{\phi}u_{n})_{n\in\N}$ is a Cauchy sequence,
$(b_h\cdot d_\phi u_{n})_{n\in\N}$
is  also a Cauchy sequence, and thus that $D(\Delta_\phi)\subset D(P_{\phi})$.\medskip

\noindent
The statement about $P_\phi^*$ is then a straightforward consequence of the above analysis.
Indeed, since
$P_\phi^*=\Delta_\phi-b_{h}\cdot d_{\phi}$ on $\ccc^\infty_c(\R^d)$, 
the above arguments 
imply that the closure of $(P_\phi^*,\ccc^\infty_c(\R^d))$ is maximal accretive and that its domain contains
$D(\Delta_{\phi})$. Moreover, $P_\phi^*$   is maximal accretive since $P_{\phi}$ is, and 
hence coincides with the closure of  $(P_\phi^*,\ccc^\infty_c(\R^d))$.\medskip

\noindent
Let us now prove the statement on the spectrum of $P_\phi$. 
Throughout, we will denote $\C_+=\{\Re(z)\geq 0\}$.
 It follows from \eqref{eq:hypgenephi} and from \eqref{hyp:controlechamp} that 
for every $u\in \ccc^\infty_c(\R^d)$, it holds, for some $C>0$ and every $h>0$ small enough,
\begin{equation}
\label{eq.C}
\begin{split}
\vert\<b_h\cdot d_\phi u,u\>\vert\leq \Vert d_\phi u\Vert \Vert b_h u\Vert&\leq C(\Vert d_\phi u\Vert^2+\Vert u\Vert \Vert d_\phi u\Vert)\,.
\end{split}
\end{equation}
Let us set $\Lambda_0=5C$ for some $C\geq1$ satisfying~\eqref{eq.C}, and let $z\in\C_+$ be such that $\vert\Im(z)\vert\geq \Lambda_0\max(\Re(z),\sqrt{\Re(z)})$. 
Suppose first that $ \Re(z) \Vert  u\Vert^2\geq\frac 12 \Vert d_\phi u\Vert^2$. Then, thanks to the  estimate~\eqref{eq.C}, we have
\begin{equation*}
\begin{split}
\vert \<(b_h\cdot d_\phi-i\Im(z))u,u\>\vert
&\geq \Big(\,\vert\Im(z)\vert-C\big(2\Re(z)+\sqrt{2\Re(z)}\big)\,\Big)\Vert u\Vert^2\\
&\geq C\max\big(\Re(z),\sqrt{\Re(z)}\big) \Vert u\Vert^2\geq  \Re(z) \Vert u\Vert^2.
\end{split}
\end{equation*}
Since $\vert \<(b_h\cdot d_\phi-i\Im(z))u,u\>\vert\leq \vert\<(P_\phi-z)u,u\>\vert$, this implies that
\be\label{eq:estresolvsect}
\vert\<(P_\phi-z)u,u\>\vert\geq  \Re(z) \Vert u\Vert^2.
\ee
Suppose now that $ \Re(z) \Vert  u\Vert^2\leq\frac 12 \Vert d_\phi u\Vert^2$. One then directly obtains
$$
\vert\<(P_\phi-z)u,u\>\vert\geq \<(\Delta_\phi-\Re(z))u,u\>\geq \Re(z)\Vert u\Vert^2\,,
$$
which, combined with \eqref{eq:estresolvsect}, implies that 
%
\be\label{eq:estresolvsect3}
\Vert(P_\phi-z)u\Vert\geq  \Re(z) \Vert u\Vert
\ee
for every $z\in \C_+\setminus\Gamma_{\Lambda_0}$ and $u\in \ccc^\infty_c(\R^d)$.
Since $P_\phi$ is closed,  it follows
that $P_\phi-z$ is injective with closed range, and hence semi-Fredholm, for every  $z\in \C_+\setminus\Gamma_{\Lambda_0}$
such that $\Re(z)\neq 0$. Assume now for a while that the fourth item in Proposition~\ref{prop:elemPphi},
which is proved  independently just below,
is satisfied, and let $\lambda\in \R$ be such that $i\lambda \in \sigma(P_\phi) $.
By assumption, $i\lambda $ is  then an eigenvalue of $P_\phi$ 
and there exists some $u\in D(P_\phi)\setminus\{0\}$ such that $P_{\phi}u=i \lambda  u$.
In particular,
it holds
$$
0\ =\  \Re\<P_\phi u, u\> \ =\ \Vert d_\phi u\Vert^2\ =\ h^{2}\Vert e^{-\frac \phi h} \nabla (e^{\frac \phi h} u)\Vert^2\,,
$$
which implies 
$u\in \sspan\{e^{-\frac \phi h}\}$ and then $\lambda=0$.
This shows that  $\sigma(P_\phi)\cap i\R\subset\{0\}$
and thus, $P_{\phi}$ being maximal accretive, that $\sigma(P_\phi)\cap \{\Re(z)\leq 0\}\subset\{0\}$.
It follows that $P_\phi-z$ is semi-Fredholm for every  $z\in \C\setminus\Gamma_{\Lambda_0}$,
and has index~$0$ on $\{\Re z\leq 0\}\setminus\{0\}$.
But the open set $\C\setminus\Gamma_{\Lambda_0}$ being connected, the index of $P_\phi-z$
 is constant, and then equal to $0$, on $\C\setminus\Gamma_{\Lambda_0}$ (see \cite[Theorem~5.17 in Chap.~4]{Kato}). 
 Hence, $P_\phi-z$ being injective on $\C\setminus\Gamma_{\Lambda_0}$, it
  is invertible from $D(P_\phi)$ onto $L^2(\R^d)$ on $\C\setminus\Gamma_{\Lambda_0}$ and
the resolvent estimate stated in Proposition~\ref{prop:elemPphi} becomes a direct consequence of \eqref{eq:estresolvsect3}.\medskip

\noindent
Let us now prove the fourth item of Proposition~\ref{prop:elemPphi}. Thanks to \eqref{eq:hypgenephi}, there exist $c>0$ and $R>0$ such that 
$$\forall \vert x\vert \geq R,\;\vert\nabla\phi(x)\vert^2\geq c.$$
Take $c_1\in(0,c)$ and let $W$ be a nonnegative smooth function such that $\supp(W)\subset B(0,R)$ and 
$W(x)+\vert\nabla\phi(x)\vert^2\geq \frac {c+c_1}2$ for all $x\in\R^d$. There exists consequently $h_0>0$ such that for all $h\in(0,h_0]$, one has
$$
\tilde W\ :=\  W+\vert\nabla\phi\vert^2-h\Delta\phi\ \geq\  c_1
$$
on $\R^d$. Introduce the operator 
$$\tilde P_\phi=P_\phi+W=-h^2\Delta+\tilde W+b_hd_\phi$$
with domain $D(P_\phi)$. Since $P_\phi$ is maximal accretive
and $W\in \ccc_{c}^{\infty}(\R^{d},\R^{+})$,  $\tilde P_\phi$
is also maximal accretive (see for example \cite[Theorem~13.25]{Hel-spectral}).
Moreover, 
for every $u\in \ccc_c^\infty(\R^d)$ and then for every $u\in D(P_\phi)$, one has
$$
\Re\<\tilde P_\phi u,u\>=\<(-h^2\Delta+\tilde W)u,u\>\geq c_1\Vert u\Vert^2\,,
$$
which implies as above that  for every $z$ in $\{\Re(z)<c_1\}$, $\tilde P_\phi-z$ is invertible from $D(P_\phi)$ onto $L^2(\R^d)$.
Hence, for every $z$ in $\{\Re(z)<c_1\}$, we can write
$$
P_\phi-z\ =\ \tilde  P_\phi-z -W\ =\ (\Id-W(\tilde  P_\phi-z)^{-1})(\tilde  P_\phi-z).
$$
Of course, $z\mapsto (\tilde  P_\phi-z)^{-1}$ is holomorphic on $\{\Re z<c_1\}$ and thanks to the analytic Fredholm theorem, it then suffices to prove that 
$$K(z)\ :=\ W(\tilde  P_\phi-z)^{-1}\,:\,L^{2}(\R^{d})\to L^{2}(\R^{d})$$ is compact  for every $z$ in $\{\Re(z)<c_1\}$. This follows
 from the compactness of the embedding $H^1_R\subset L^2(\R^{d})$ and
from the fact that for every $z\in \{\Re z<c_1\}$,  $K(z)$ acts continuously from $L^{2}(\R^{d})$
into $ H_R^1$, where
$$ H_R^1\ :=\ \{u\in H^1(\R^d),\;\supp(u)\subset B(0,R)\}\,.$$
Indeed, for any  $z$ in $\{\Re(z)<c_1\}$, the operator
$d_\phi(\tilde P_\phi-z)^{-1}: L^2(\R^{d})\to L^{2}(\R^{d})$ is continuous 
thanks to \eqref{eq:accretPphi} and hence, since $W$ is smooth and  supported in $B(0,R)$,  $K(z):L^2(\R^{d})\to H^1_{R}$
is also continuous.\medskip

\noindent
To conclude, it remains to prove the last statement of  Proposition~\ref{prop:elemPphi}.
To this end, note first that $P_{\phi} e^{-\frac \phi h} =0$ according to \eqref{eq.P}
and let us recall that, according to \eqref{eq:hypgenephi},
 $ e^{-\frac \phi h}\in D(\Delta_{\phi})\subset D(P_{\phi})$. Thus, $\sspan\{e^{-\frac \phi h}\}\subset \Ker P_{\phi} $ and $0$ is an eigenvalue of $P_{\phi}$. It has moreover finite algebraic multiplicity
according to the preceding analysis.
Conversely, the relation
$$
\forall\,u\,\in\, D(P_\phi)\,,\ \ \ 
\Re\< P_{\phi} u, u\>\ =\ \Vert d_\phi u\Vert^2\ =\ h^{2}\Vert e^{-\frac \phi h} \nabla (e^{\frac \phi h} u)\Vert^2
$$
leads to $ \Ker P_{\phi} \subset \sspan\{e^{-\frac \phi h}\}$
and the same arguments also show that $\Ker P^{*}_{\phi}=\sspan\{e^{-\frac \phi h}\}$.  This implies that
 $0$ is an eigenvalue of $P_{\phi}$ with algebraic multiplicity one. Indeed, 
if it was not the case, there would exist $u\in D(P_{\phi})$ such that 
$u\notin \Ker P_{\phi}$ and $P_{\phi}u= e^{-\frac \phi h}$, and hence such that
$$
0\ <\ \lp P_{\phi}u, e^{-\frac \phi h}\rp\ =\ \lp u,  P^{*}_{\phi}e^{-\frac \phi h}\rp\ =\ 0\,.
$$

\subsection{Spectral analysis near the origin}
\label{sub-spec-origin}
Let us denote by $(e_{k}^{W})_{k\geq 1}$ the eigenfunctions of 
$\Delta_\phi$ associated with the non-decreasing sequence of eigenvalues $(\lambda_{k}^W)_{k\geq 1}$.
Let $\epsilon_0$ and $h_0>0$ be given by Proposition \ref{prop:localspectWitten}.  We recall that 
for every $h\in(0,h_{0}]$, it holds
$$
{\rm card}\big(\sigma(\Delta_\phi)\cap \{\Re z < \epsilon_0 h \}\big)\ =\ n_0\,,$$
where $n_0$ is the number of local minima of $\phi$.
We define
$$
\begin{array}{rcl}
R_-\,:\,\C^{n_0}&\longrightarrow& L^2(\R^d)\\
(\alpha_{k})&\longmapsto& \sum_{k=1}^{n_0}\alpha_{k} e_{k}^{W}
\end{array}
$$
and $R_+:=R_-^*$, i.e.
$$
\begin{array}{rcl}
R_+\,:\,L^2(\R^d)&\longrightarrow&\C^{n_0}\\
u&\longmapsto& (\<u,e_{k}^W\>)_{k=1,\ldots,n_0}\,.
\end{array}
$$
Note in particular the relations
\be
\label{eq.Pi}
R_+R_-\ =\ \Id_{\C^{n_0}}\quad\text{and}\quad R_-R_+\ =\ \Pi\,,
\ee 
where $\Pi$ denotes the orthogonal projection onto $\Ran( R_-)=\sspan\big(e_{k}^{W}, k\in\{1,\dots,n_{0}\}\big)$.
We also define the spectral projector 
$$
\hat \Pi\ :=\ 1-\Pi\,.
$$
For $z\in\C$, let us then consider on the Hilbert space $\hat E:=\Ran (\hat\Pi)$  the following unbounded operator
which will be useful in the rest of this section: 
\be
\label{eq.hatP}
\hat P_{\phi,z}\ :=\ \hat\Pi(P_\phi-z)\hat\Pi \quad\text{with domain}\quad D(\hat P_{\phi,z})\ :=\ \hat\Pi (D(P_\phi))\,.
\ee
 Hence  $D(\hat P_{\phi,z})$ is dense in $\hat E$ and, since $\Ran\,\Pi\subset D(\Delta_\phi)\subset D(P_\phi)$, it holds
$\hat\Pi (D(P_\phi))\subset D(P_\phi)$ and $\hat P_{\phi,z}$ is well and densely defined.

\begin{lemma}
\label{le.hatP}
Let $\epsilon_0$ and $h_0>0$ be given by Proposition~\ref{prop:localspectWitten}.
Then, for every $h\in(0,h_0]$, 
 the operator $\hat P_{\phi,z}:D(\hat P_{\phi,z})\to \hat E$ 
defined in \eqref{eq.hatP} 
 is invertible on $\{\Re z < \epsilon_0h \}$.
 Moreover, for any $\epsilon_{1}\in(0,\epsilon_{0})$
 it holds:
 $$
 \forall\,z\,\in \,\{\Re z < \epsilon_1h \}\,,\ \   \|\hat P_{\phi,z}^{-1}\|_{\hat E\to \hat E}\ =\ \mathcal O(h^{-1})\,,
 $$
 uniformly with respect to $z$.
\end{lemma}

\noindent
\bp
We begin by the following  observation:
 the unbounded operator  
 $$\hat\Pi(P^{*}_\phi-z)\hat\Pi \quad\text{with domain}\quad
\hat\Pi (D(P^{*}_\phi))\ \subset\  D(P^{*}_\phi)$$
is well and densely defined on $\hat E$, and satisfies moreover 
$$\hat\Pi(P^{*}_\phi-z)\hat\Pi \ =\ \hat P^{*}_{\phi,z}\,.$$ 
Indeed, the relation
$
\<\hat\Pi(P_\phi-z)\hat\Pi v, w\>
 = 
\< v, \hat\Pi(P^{*}_\phi-z)\hat \Pi w\>
$,
valid
for every $v\in D(P_\phi)$ and $w\in D(P^{*}_\phi) $,
implies that $\hat\Pi(P^{*}_\phi-z)\hat\Pi \subset \hat P^{*}_{\phi,z}$.
Moreover, for every $v\in D(P_{\phi})$ and $w\in D(\hat P_{\phi,z}^*)$, one has
\begin{align*}
   \<(P_{\phi}-z) v,w\>&\ =\ \<(P_{\phi}-z)\Pi v,w\>+\<(P_{\phi}-z)\hat \Pi v,w\>\\
   &\ = \ \<(P_{\phi}-z)\Pi v,w\>+\<\hat \Pi v,\hat P_{\phi,z}^*w\>.
\end{align*}
Since $P_{\phi}\Pi $ is continuous, $\Pi$ being continuous with finite rank, one has   $\vert \<P_{\phi}\Pi v,w\>\vert\leq C\Vert v\Vert\Vert w\Vert$ for some $C>0$ independent of $(v,w)$, which implies that $w \in D(P_{\phi}^*)$. Hence
   $D(\hat P_{\phi,z}^*)\subset D(P_{\phi}^*)$ 
 and   since $\Ran(\Pi)\subset D(\Delta_\phi)\subset D(P_\phi^*)$, this implies $\hat\Pi (P_{\phi}^*-z)\hat\Pi = \hat P_{\phi,z}^*$.\medskip

\noindent
Let now  consider $z$ in $\{\Re z < \epsilon_{0}h\}$
and
let us prove  that $\hat P_{\phi,z}$
 is invertible from $D(\hat P_{\phi,z})$ onto $\hat E$.
 First, according to Proposition~\ref{prop:localspectWitten},
  we have for every $u\in D(\Delta_{\phi})$,
 \begin{align}\nonumber
\Re\<(P_\phi-z)\hat \Pi u,\hat \Pi u\>&\ =\ \<(\Delta_\phi-\Re(z))\hat \Pi u,\hat \Pi u\>\\
\label{eq:corcivPphi}
&\ \geq\  (\epsilon_0h-\Re z)\Vert \hat \Pi u\Vert^2\,,
\end{align}
and the inequality \eqref{eq:corcivPphi} is also true when $u\in D( P_{\phi})$.
 Indeed, for any $u\in D( P_{\phi})$, there exists a sequence
 $(u_{n})_{n\in\N}$ in $D(\Delta_{\phi})$ such that $u_{n}\to u$ and $P_{\phi}u_{n}\to P_{\phi}u$
 in $L^{2}(\R^{d})$. Hence $\hat \Pi u_{n}\to \hat \Pi u$ and, since $P_{\phi}\Pi$ is continuous,
 it also holds
$P_{\phi}\hat\Pi u_{n}\to P_{\phi}\hat\Pi u$. 
In particular,
it  follows that
 $\hat P_{\phi,z}$ is injective.
 Note that a similar analysis  shows that
 $\hat P_{\phi,z}^*$
is also  injective.
 \\
Second,
  let us show that 
  $\hat P_{\phi,z}$ is closed, which will in particular
  imply that  $\Ran(\hat P_{\phi,z})$  is closed according to
 \eqref{eq:corcivPphi}. 
   For shortness, we denote $\hat P=\hat P_{\phi,z}$ and $P=P_\phi$. Suppose that $(u_n)_{n\in\N}$ is a sequence in 
  $D(\hat P)\subset D(P)$ such that $u_n\rightarrow u$ and $\hat P u_n\rightarrow v$ in $\hat E$. 
  Since 
  $\Ran\,\Pi\subset D(\Delta_\phi)\subset D(P^*)$, it holds
   \begin{align*}
 \Pi P u_{n}\ =\ \sum_{k=1}^{n_{0}}\<Pu_{n},e_{k}^{W} \> e_{k}^{W} \ =\ 
\sum_{k=1}^{n_{0}}\<u_{n},P^{*}e_{k}^{W} \> e_{k}^{W}
  \  \underset{n\to +\infty}{\longrightarrow}\ 
 \sum_{k=1}^{n_{0}}\<u,P^{*}e_{k}^{W} \> e_{k}^{W}\,,
\end{align*}
and thus
   $(P--z)u_n=\hat P u_n+\Pi (P-z)u_n$
   converges. Since $P$ is closed, this implies that
 $u\in D(P)\cap \Ran \hat\Pi=\hat\Pi(D(P))$ and that 
    $$(P-z) u\ =\ v +g \quad\text{with} \ \ g\in \Ran\,\Pi\,.$$ 
    Multiplying this relation by $\hat \Pi$, we get 
   $v=\hat P u$, which proves that $\hat P $ is closed.\\ 
   To prove that $\hat P$ is invertible from $D(\hat P)$ onto $\hat E$, it is thus enough to prove that $\Ran(\hat P)$ is dense in $\hat E$. Let  then $v\in \hat E$ be such that 
   $\<\hat P u,v\>=0$ for all $u\in D(\hat P)$. Then $v\in D(\hat P^*)$ and  $\hat P^* v=0$.
   By injectivity of $ \hat P^*$, it thus holds $v=0$, which
 proves the invertibility of $\hat P:D(\hat P_{\phi,z})\to\hat E$.\medskip

\noindent    
The relation \eqref{eq:corcivPphi} then implies that  for all $z\in\{\Re z\leq \epsilon_1 h\}$, one has 
$$
 \Re\<(P_\phi -z)\hat \Pi u,\hat \Pi u\>\geq \delta h\Vert\hat\Pi u\Vert^2
$$
 with $\delta=\epsilon_0-\epsilon_1>0$.
 Hence, for the operator norm on $\hat E \subset L^2(\R^{d})$, one has
$$
\hat P_{\phi,z}^{-1}=\ooo(h^{-1})\,,
$$
 uniformly with respect to 
$z\in\{\Re z < \epsilon_1 h \}$.
\ep

\noindent
 For  $z\in\C$, we now consider the Grushin operator $\ppp_\phi(z):D(P_\phi)\times \C^{n_0}\rightarrow L^2(\R^d)\times \C^{n_0}$
defined by
\be
\ppp_\phi(z)=
\left(
\begin{array}{cc}
P_\phi-z&R_-\\
R_+&0
\end{array}
\right).
\ee

\begin{lemma}\label{prop:invert-grushin} 
Let $\epsilon_0$ and $h_0>0$ be given by Proposition~\ref{prop:localspectWitten}.
Then, the operator $\ppp_\phi(z)$ is invertible on $\{\Re z < \epsilon_0h \}$.
More precisely, for every $z\in \{\Re z < \epsilon_0h \}$, $(u,u_-)\in D(P_\phi)\times \C^{n_0}$
and $(f,y)\in L^2(\R^d)\times \C^{n_0}$,
it holds 
$$\ppp_\phi(z)(u,u_{-})\ =\ (f,y)$$ if and only if
$$
(u,u_{-})\ =\ \big(\,R_{-}y+v\,,\, R_{+}f-R_{+}(P_{\phi}-z)R_{-}y - R_{+}P_{\phi}v\,\big)\,,
$$
where 
$$
v\ :=\ \hat P_{\phi,z}^{-1}\hat\Pi f-\hat P_{\phi,z}^{-1}\hat\Pi P_\phi R_-y
\ \in\ \hat\Pi (D(P_\phi))\,.
$$
\end{lemma}


\noindent
\bp 
Let $(f,y)\in L^2(\R^d)\times \C^{n_0}$ and assume  that $(u,u_-)\in D(P_\phi)\times \C^{n_0}$ satisfies
\be\label{eq:grushin00}
\left\{
\begin{array}{c}
(P_\phi-z)u+R_-u_-=f\\
R_+u=y.\phantom{*******}
\end{array}
\right.
\ee
Applying $R_{+}$ to the first equation and $R_-$ to the second one, we get, according to \eqref{eq.Pi}: 
$$u_{-}=R_{+}f-R_{+}(P_{\phi}-z)u \quad\text{and}\quad
u=R_-y+v\,,$$ with $v\in\Ran\, \hat\Pi \cap D(P_{\phi})=\hat\Pi (D(P_\phi))$ solution to
$$
(P_\phi-z)R_-y+(P_\phi-z)v+R_-u_-=f.
$$
Then, applying  $\hat\Pi$ to the latter equation,  we get, using $\hat\Pi R_-=0$,
\be\label{eq:grushin1}
\hat\Pi(P_\phi-z)\hat\Pi v\ =\ \hat\Pi f-\hat\Pi(P_\phi-z)R_-y-
\hat\Pi R_{-}u_{-}
\ =\ \hat\Pi f-\hat\Pi P_\phi R_- y\,.
\ee
Conversely, note  that if $v\in \Ran \,\hat\Pi \cap D(P_{\phi})$ is solution to \eqref{eq:grushin1},
then according to \eqref{eq.Pi}, 
$$\big(\,u= R_-y+v\,,\,u_{-}=R_{+}f-R_{+}(P_{\phi}-z)(R_-y+v) \,\big)
\ \in\ D(P_\phi)\times \C^{n_0}$$
is solution to \eqref{eq:grushin00}.\medskip

\noindent
Hence, the statement of Lemma~\ref{prop:invert-grushin} simply follows from
Lemma~\ref{le.hatP} which implies that, for every $z\in \{\Re z < \epsilon_0h \}$,
$$v\ =\ \hat P_{\phi,z}^{-1}\hat\Pi f-\hat P_{\phi,z}^{-1}\hat\Pi P_\phi R_-y\ \in\ \hat\Pi (D(P_\phi))
$$
is the unique solution to 
\eqref{eq:grushin1}.
\ep

\noindent
\emph{Proof of Theorem~\ref{th:small-ev-Pphi}}. Let $\epsilon_0$ and $h_0$  be as in 
Lemmata~\ref{le.hatP} and \ref{prop:invert-grushin}, and take $\epsilon_{1}\in(0,\epsilon_{0})$. For 
$z\in \{\Re z < \epsilon_0h \}$, let  $\eee_\phi(z)=\ppp_\phi(z)^{-1}$. 
According to Lemma~\ref{prop:invert-grushin}, 
it thus holds
$$
\eee_\phi(z)=
\left(
\begin{array}{cc}
E(z)&E_+(z)\\
E_-(z)&E_{-+}(z)
\end{array}
\right)\,,
$$
where $E,E_-,E_+,E_{-+}$ are holomorphic in $\{\Re z < \epsilon_0h \}$
and satisfy the following  formulas:
\be\label{eq:grushin40}
E_+(z)\ =\ R_--\hat P_{\phi,z}^{-1}\hat\Pi P_\phi R_-\,,\ \ E_-(z)\ =\ R_+-R_+P_\phi \hat P_{\phi,z}^{-1}\hat\Pi \,,
\ee
\be\label{eq:grushin4}
E_{-+}(z)\ =\ -R_+(P_\phi-z)R_-+R_+P_\phi  \hat P_{\phi,z}^{-1}\hat\Pi P_\phi R_-
\ee
and 
\be\label{eq:grushin5}
E(z)\ =\ \hat P_{\phi,z}^{-1}\hat\Pi\,.
\ee
Moreover, $P_\phi-z$ is invertible if and only if $E_{-+}(z)$ is, in which case
it holds
\be\label{eq:grushin3}
(P_\phi-z)^{-1}=E(z)-E_+(z)E_{-+}(z)^{-1}E_-(z).
\ee
We refer in particular to \cite{SjZw} for more details in this connection.\medskip

\noindent
We now want to use these formulas to compute the number of poles of 
$(P_\phi-z)^{-1}$. Thanks to  \eqref{eq:estimbdphi}, one has, for some $C>0$ and all  $k\in\{1,\ldots, n_0\}$,
$$
\Vert b_h\cdot d_\phi e_{k}^W\Vert
\leq C \big(\Vert\Delta_\phi e_{k}^W\Vert+\Vert d_\phi e^W_k\Vert\big)
\leq C(\lambda_{k}^W+\sqrt{ \lambda_{k}^W})\,.
$$
Using the bound  $\lambda_{k}^W\leq Ch e^{-2\frac Sh}$ given by Proposition~\ref{prop:localspectWitten}, this yields the existence of some $C>0$ such that
for every  $k\in\{1,\ldots, n_0\}$,
\be\label{eq:QMPphi}
\Vert b_h\cdot d_\phi e_{k}^W\Vert\leq C \sqrt h e^{-\frac Sh}
\quad\text{and}
\quad
\Vert P_\phi e_{k}^W\Vert\leq C \sqrt h e^{-\frac Sh}\,.
\ee
This shows that 
$R_+\Delta_\phi R_-=\ooo(h e^{-2\frac Sh})$ and
$R_+ b\cdot d_\phi R_-=\ooo(\sqrt he^{-\frac Sh})$. 
Hence, for all $z\in\C$, it holds
\be\label{eq:grush}
\begin{split}
R_+(P_\phi-z)R_-&=R_+P_\phi R_--z\,\Id_{\C^{n_0}}\\
&=-z\,\Id_{\C^{n_0}}+\ooo(\sqrt he^{-\frac Sh}).
\end{split}
\ee
On the other hand, 
we deduce from
 \eqref{eq:QMPphi}
 and from the related relation
$$
\<P_\phi u,e_k^W\>=\<u,\Delta_\phi e_k^W\>-\< u,b_{h}\cdot d_\phi e_k^W\>
=\ooo(he^{-2\frac Sh}+\sqrt he^{-\frac Sh})\Vert u\Vert,
$$
valid for any $u\in D(P_{\phi})$ and $k\in\{1,\ldots, n_0\}$,
that
\begin{equation}
\label{eq.norm-PR}
P_\phi R_- = \ooo(h^{\frac 12} e^{-\frac Sh})\quad
\text{and}\quad
R_+P_\phi=\ooo( h^{\frac12}e^{-\frac Sh})\,.
\end{equation}
Moreover, 
we know from Lemma~\ref{le.hatP}
 that, uniformly on $\{\Re z<\epsilon_1 h\}$, it holds $\hat P_{\phi,z}^{-1}=\ooo(h^{-1})$. Therefore,
 injecting this estimate and
 \eqref{eq:grush}, \eqref{eq.norm-PR} into 
  \eqref{eq:grushin4} and \eqref{eq:grushin40}, we obtain respectively,
  uniformly on $\{\Re z<\epsilon_1 h\}$,
 \be\label{eq:estimE-+}
E_{-+}(z)\ =\ z\,\Id_{\C^{n_0}}+\ooo( h^{\frac12}e^{-\frac Sh})
\ee
and 
  \be\label{eq:estimE+/E-}
 E_{+}(z)\ =\ R_{-}+\ooo(h^{-\frac 12}e^{-\frac Sh})
\quad \text{and}
\quad
E_{-}(z)\ =\ R_{+}+\ooo(h^{-\frac 12}e^{-\frac Sh})\,.
\ee
According to \eqref{eq:estimE-+}, 
$E_{-+}(z)$ is then invertible when 
$z \in\{\Re z < \epsilon_1h \}$ satisfies $|z|\geq C h^{\frac12}e^{-\frac Sh} $ for $C$ large enough
and  the spectrum of $P_{\phi}$ in $\{\Re z < \epsilon_1h \}$ is then of order $\mathcal O( h^{\frac12}e^{-\frac Sh})$.
Moreover, 
for $|z|= \frac{\epsilon_{1}}{2}h $, it holds
\be\label{eq:expandE-+'}
E_{-+}(z)=z\big(\Id_{\C^{n_0}}+\ooo(h^{-\frac12}e^{-\frac Sh})\big)
\ee
and
injecting \eqref{eq:expandE-+'} and \eqref{eq:estimE+/E-}
into \eqref{eq:grushin3} shows that 
$$
(P_\phi-z)^{-1}=E(z)- \frac 1z\big(\Pi+ \ooo(h^{-\frac12}e^{-\frac Sh})\big)\,.$$
Thus, the spectral projector on the open disk $D(0,\frac{\epsilon_{1}}{2}h)$ satisfies
$$
\Pi_{D(0,\frac{\epsilon_{1}}{2}h)}:=-\frac{1}{2\pi i}\int_{\pa D(0, \frac{\epsilon_{1}}{2}h)}(P_\phi-z)^{-1} dz
= \Pi +\ooo(h^{-\frac12}e^{-\frac Sh})\,,
$$
where we recall that $\Pi$ is a projector of rank $n_{0}$.
This implies that for  every $h>0$ small enough,
the rank of $\Pi_{D(0,\frac{\epsilon_{1}}{2}h)}$, which is 
the number of eigenvalues of $P_{\phi}$ in $D(0,\frac{\epsilon_{1}}{2}h)$
counted with algebraic multiplicity,
is precisely $n_0$.\medskip

\noindent
In order to achieve the proof of Theorem~\ref{th:small-ev-Pphi}, it just remains to prove the resolvent estimate stated there. 
On the one hand, it follows easily from  \eqref{eq:grushin5}, \eqref{eq:estimE+/E-}, and 
Lemma~\ref{le.hatP} that 
$$E(z)=\ooo(h^{-1})\,,\ E_-(z)=\ooo(1)\,,\ \ \text{and}\ \ E_+(z)=\ooo(1)\,,$$
uniformly with respect to  $z\in \{\Re z < \epsilon_1h \}$.
On the other hand, taking $\epsilon\in(0,\epsilon_1)$, it follows from \eqref{eq:estimE-+} that 
$E_{-+}^{-1}(z)=\ooo(h^{-1})$, uniformly with respect to  $z\in \{\Re z < \epsilon_1h \}\cap\{|z|>\epsilon h\}$. Plugging all these estimates into \eqref{eq:grushin3}, we obtain the announced result.\medskip

\noindent
Eventually, since 
$\sigma(P_\phi^*)=\overline{\sigma(P_\phi)}$ and, for all $z\notin \sigma(P_\phi)$, $\Vert (P^{*}_\phi-\overline z)^{-1}\Vert=\Vert(P_\phi- z)^{-1}\Vert,$
it follows
easily that the conclusions of Theorem \ref{th:small-ev-Pphi} also hold true for $P_\phi^*$.
\ep
 
\section{Geometric preparation}\label{sec:preuvelemalg}

\noindent
Let us begin this section by  observing
that the identity $b\cdot\nabla V=0$ arising from \eqref{eq:propvectfield} implies that $\uuu\subset \{x\in\R^{d},\ b(x)=0\}$, where
we recall that $\uuu$ denotes the set of critical points of the Morse function $V$, as it can be easily  proved using a Taylor expansion. 
Moreover, we have the following 
\begin{lemma}\label{lem:structure-b}Suppose that  Assumptions \ref{hyp:statmes} and \ref{hyp:phiMorse} hold true and let $\u\in\uuu$ be a critical point of $V$. Then, there exists a smooth map $J_\u:\R^{d}\rightarrow\mmm_d(\R)$ such that  $J_\u(\u)$ is antisymmetric and
$b(x)=J_\u(x)\nabla V(x)$ for all $x$ in some neighborhood of $\u$.
Moreover, it holds
$$
J_\u(\u)\ =\ B(\u)\Hess V(\u)^{-1}\,,
$$
where $B(\u)={\rm Jac}_\u b$ is the Jacobian matrix of $b$ at $\u$.
\end{lemma}

\noindent
\bp Let $\u\in\uuu$ that we assume to be $0$
to lighten the notation. Thanks to the Taylor formula, there exists  a  smooth map 
$G:\R^{d}\rightarrow \mmm_d(\R)$ such that $b(x)=G(x) x$ for all $x\in \R^{d}$ and $G(0)={\rm Jac}_{0}\,b$. The same construction works for $\nabla V$ and 
denoting by $\sss_d$ the set of symmetric matrices, there exists
 a smooth map $A:\R^{d}\rightarrow\sss_d$ such that $\nabla V(x)=A(x)x$
 for all $x\in \R^{d}$
 and $A(0)=\Hess V(0)$. The equation 
 $\<b(x),\nabla V(x)\>=0$ for all $x\in\R^{d}$ then yields
 $\<G(x)x,A(x)x\>=0$
 and hence, since $A(x)$ is symmetric,
 $
 \<A(x)G(x)x,x\>=0
 $
 for all  $x\in\R^{d}$. Expanding $A(x)G(x)$ in powers of $x$, this implies that 
 $$
\forall \,x\,\in\,\R^{d}\,,\ \ \  \<A(0)G(0)x,x\>\ =\ 0\,.
 $$
Hence, the matrix $A(0)G(0)$ is antisymmetric. Since $A(0)$ is symmetric and invertible (since $V$ is a Morse function), this implies that $G(0)A(0)^{-1}$ is antisymmetric.
Moreover, $A(x)$ is then also invertible in a neighborhood $\mc V$ of $0$
and we can thus define 
 $J_{0}(x)=G(x)A(x)^{-1}$ on $\mc V$. One then has 
 $$
 J_{0}(x)\nabla V(x)\ =\ G(x)A(x)^{-1}A(x)x\ =\ b(x)
 $$
 for all $x\in\mc V$
 and 
 $
 J_{0}(0)=G(0)A(0)^{-1}
 $
 is antisymmetric thanks to the above analysis.
\ep

\begin{remark}\label{rem:susy}
%
It is not clear from the above proof that 
the relation $b\cdot \nabla V=0$ implies the existence of a 
smooth
map $J:\R^{d}\rightarrow\mmm_d(\R)$ with antisymmetric
matrices values such that $b=J\,\nabla V$.
However, 
it follows from  \eqref{eq:propvectfield} that
for such a map $J$, the
vector fields of the form  $b_{h}=J\,\nabla V+h\nu$
 enter in our framework as soon as
 \be 
\label{eq.J-anti}
\div \nu\ =\ 0\quad\text{and}\quad
  \big(\sum_{i=1}^{d}\pa_{i} J_{ij}\big)_{j=1,\dots,d}\cdot \nabla V \ =\  \nu\cdot \nabla V\,.
  \ee
This is for instance the case when 
$
\nu = \big(\sum_{i=1}^{d}\pa_{i} J_{ij}\big)_{j=1,\dots,d}
$, 
which is in particular satisfied when $J$ appears to be constant.
Moreover, when  $
\nu = \big(\sum_{i=1}^{d}\pa_{i} J_{ij}\big)_{j=1,\dots,d}
$,
$L_{V,b,\nu}$ (or equivalently $P_{\phi}$)
admits a supersymmetric structure according to (see indeed \eqref{eq.J-supersym'})
$$
L_{V,b,\nu}  \ =\ -h \,e^{\frac Vh} \,\div \circ \big( \,e^{-\frac Vh} \big( I_{d} - J \big)\nabla \,  \big)
\ =\ h \,\nabla^{*} \big( I_{d} - J \big)\nabla\,,$$
where the adjoint is considered with respect to $m_{h}$ 
(or equivalently 
$$
P_{\phi}\ =\ \Delta_{\phi}+ b_{h}\cdot d_{\phi}\ =\ d_\phi^*\big(I_{d}- J\big)d_\phi\,,
$$
where the adjoint is now considered with respect to the Lebesgue measure).
Using this structure, we may follow the general approach of \cite{HeHiSj11_01} to analyse the spectrum of $P_\phi$. 
Nevertheless, the operator $P_\phi$ still does not have any PT-symmetry and following this approach would again require 
to replace the use of the Fan inequalities by the one of
Theorem~\ref{th:specgradedmat} in the final part of the analysis. We believe
that this approach may yield complete asymptotic expansions of the small eigenvalues of $P_\phi$ (or $L_{V,b,\nu}$)
in this  setting.\medskip

\noindent
However, when $J$ has antisymmetric matrices values and  \eqref{eq.J-anti} holds but $\nu \neq \big(\sum_{i=1}^{d}\pa_{i} J_{ij}\big)_{j=1,\dots,d}$,
 the operator $P_\phi$ is not supersymmetric anymore
 (see \cite{Mi16} for related results). 
\end{remark}

\noindent
We are now in position to prove Lemma \ref{le.B}.
Throughout the rest of this section, we denote 
$$-\mu_{1}<0<\mu_{2}\leq\cdots\leq \mu_{d}$$
the eigenvalues of $\Hess V(\s)$ counted with multiplicity. 
For shortness, we will denote 
$$B\ =\ B(\s)\ =\ {\rm Jac}_{\s}b
\quad\text{and}
\quad
J\ =\ J(\s)\ =\ B(\s)(\Hess V(\s))^{-1}\,.
$$ 
We recall from Lemma~\ref{lem:structure-b} that $J$ is antisymmetric.\medskip

\noindent
{\bf Step 1 :} Let us first prove that $\det(\Hess V(\s)+B^{*})<0$. Since the matrix 
$\Hess V(\s)+B^{*}$ is real, it thus admits at least one negative eigenvalue.\\[0.3cm]
Since $\Hess V(\s)$ is real and symmetric, there exists $P\in \mathcal M_{d}(\R) $ such that 
$$P^{*}=P^{-1}\quad
\text{and}\quad
 \Hess V(\s)= P \,D\,P^{-1}\,,$$
 where $D:={\rm Diag}(-\mu_{1},\mu_{2},\dots,\mu_{d})$.
 It then holds:
\begin{equation}
\label{eq.decomp-B}
 \Hess V(\s)+B^{*}  \ =\ \Hess V(\s)\,(I_{d}-J  )\ =\ 
 P\,D\,(I_{d}-P^{-1}\,J  \,P)\,P^{-1}\,. 
\end{equation}
 Since $(P^{-1}\,J  \,P)^*=-P^{-1}\,J   \,P$,
 there exist moreover $p\in \{0,\dots,\lfloor{\frac d2}\rfloor\}$, $\eta_{1},\dots,\eta_{p}>0$,
 and $Q \in \mathcal M_{d}(\R) $ satisfying $Q^{*}=Q^{-1}$ such that
 $$
Q^{-1}\,P^{-1}\,J  \,P\,Q\ =\  \begin{bmatrix}
A_{1} && (0) &\\
& \ddots & &&\\
 (0) & & A_{p} &&\\
& & &   (0) &
\end{bmatrix}
 $$
where, for every  $k\in\{1,\dots,p\}$,
$$
A_{k}=\begin{bmatrix} 0 & -\eta_{k}\\
\eta_{k} & 0 \end{bmatrix}\,.
$$
 Here, the rank of the matrix  $J$  is $2p$ and its nonzero eigenvalues
 are the $\pm i \eta_{k}$, $k\in \{1,\dots,p\}$. Therefore,
 it holds
\begin{equation}
\label{eq.spec-decomp}
Q^{-1}\,(I_{d}-P^{-1}\,J  \,P)\,Q\ =\ \begin{bmatrix}
B_{1} && (0) &\\
& \ddots & &&\\
 (0) & & B_{p} &&\\
& & &   I_{d-2p} &
\end{bmatrix}
\end{equation}
where,  for every  $k\in\{1,\dots,p\}$,
$$
B_{k}=\begin{bmatrix} 1 & \eta_{k}\\
-\eta_{k} & 1 \end{bmatrix}\,.
$$
We then deduce from \eqref{eq.decomp-B} and \eqref{eq.spec-decomp} that
$$
\det(\Hess V(\s)+B^{*})\ =\ -(\Pi_{k=1}^{d}\mu_{k})\,(\Pi_{k=1}^{p}(1+\eta_{k}^{2}))\ <\ 0\,,
$$
which concludes this first step.\\[0.2cm]
{\bf Step 2 :} Let us denote by $\mu$ a negative eigenvalue of $\Hess V(\s)+B^{*}$
and let us show that $\mu$ is its only negative eigenvalue 
and has geometric multiplicity one.\\[0.3cm]
Assume first by contradiction that $\mu$ has geometric multiplicity two and denote by $\xi_{1},\xi_2$
two  associated unitary eigenvectors such that $\lp\xi_{1},\xi_{2}\rp=0$. Let us also define
$\xi'_{i}:=P^{-1}\xi_{i}$ for $i\in\{1,2\}$ so that $\xi'_{1}$ and $\xi'_{2}$ are orthogonal and unitary. 
According to \eqref{eq.decomp-B}, it holds moreover for $i\in\{1,2\}$,
\begin{equation*}
\label{eq.D-1}
D\,(I_{d}-P^{-1}\,J  \,P)\,\xi'_{i}\ =\ \mu\,\xi_{i}'\quad\text{and hence}
\quad D^{-1}\xi_{i}'\ =\ \frac{1}{\mu}(I_{d}-P^{-1}\,J  \,P)\,\xi'_{i}\,.
\end{equation*}
In particular, since $(P^{-1}\,J  \,P)^*=-P^{-1}\,J   \,P$, it holds for every $(a,b)\in\R^{2}$ satisfying $a^{2}+b^{2}=1$:
\begin{equation*}
\label{eq.Min-Max}
\lp D^{-1}(a \xi'_{1} + b\xi'_{2}),a \xi'_{1} + b\xi'_{2}\rp
\ =\ \frac{1}{\mu}\,.
\end{equation*}
Applying the Max-Min principle to the symmetric matrix $D^{-1}$, this shows that the second eigenvalue $\mu_2(D^{-1})$ of the matrix $D^{-1}$
satisfies $\mu_{2}(D^{-1})\leq \frac{1}{\mu}<0$, contradicting 
$D^{-1}={\rm Diag}(-\frac1{\mu_{1}},\frac1{\mu_{2}},\dots,\frac1{\mu_{d}})$.\\[0.1cm]
Hence the negative eigenvalue $\mu$ has geometric multiplicity one and we have to show that it is the only negative eigenvalue
of $\Hess V(\s)+B^{*}$. We reason again by contradiction, assuming that
$\Hess V(\s)+B^{*}$ admits another negative eigenvalue that we denote by $\eta$. Note in particular
that it follows from the relation (see indeed \eqref{eq.decomp-B})
$$
 \Hess V(\s)\,(I_{d}+J  )\ =\ \Hess V(\s)\big(\Hess V(\s)+B^{*}\big)^{*}(\Hess V(\s))^{-1}
$$
that $\eta$ is also an eigenvalue of $\Hess V(\s)-B^{*}(\s)=\Hess V(\s)\,(I_{d}+J  )$.
Denote now by $\xi_{1}$ a unitary eigenvector of $\Hess V(\s)+B^{*}$
associated with $\mu$ and by $\xi_{2}$ a unitary eigenvector of $\Hess V(\s)-B^{*}$
associated with $\eta$. 
Defining again $\xi'_{i}:=P^{-1}\xi_{i}$ for $i\in\{1,2\}$,
we have thus
$$
D^{-1}\xi_{1}'\ =\ \frac{1}{\mu}(I_{d}-P^{-1}\,J  \,P)\,\xi'_{1}
\quad\text{and}\quad
D^{-1}\xi_{2}'\ =\ \frac{1}{\eta}(I_{d}+P^{-1}\,J  \,P)\,\xi'_{2}\,.
$$
It follows that
$$
\lp D^{-1}\xi_{1}',\xi'_{2}\rp\ =\ 0\ ,\ \ \  \lp D^{-1}\xi_{1}',\xi'_{1}\rp\ =\ \frac1\mu
\ \ \ \text{and}\ \ \ \lp D^{-1}\xi_{2}',\xi'_{2}\rp\ =\ \frac1\eta\,.
$$
The vectors $\xi_{1}'$ and $\xi_{2}'$ are in particular linearly independent
and  it holds for some positive constant $c$ and every $(a,b)\in\R^2\setminus\{(0,0)\}$,
$$
\lp D^{-1}(a\xi_{1}'+b\xi'_{2}),a\xi_{1}'+b\xi'_{2}\rp\ =\ \frac{a^{2}}\mu+
\frac{b^{2}}\eta\ \leq\ -c\,\|a\xi_{1}'+b\xi'_{2}\|^{2}
$$
Applying again the Max-Min principle to the symmetric matrix $D^{-1}$ leads to 
$\mu_{2}(D^{-1})\leq -c<0$ and hence to
a contradiction.
This concludes the proof of the second step.
\\[0.2cm]
{\bf Step 3 :} Let us now prove the relation
\begin{equation}
\label{eq.det}
\det \big(\,\Hess V(\s)+2|\mu|\,\xi\,\xi^{*}\,\big)\ =\ -\det\Hess V(\s)\,,
\end{equation}
which is equivalent to
\begin{equation}
\label{eq.det'}
\det\big(\,I_{d}+2\,|\mu|D^{-1}\,\xi'\,\xi'^{*}\,\big)\ =\ -1\,,
\end{equation}
where $\xi$ denotes a unitary eigenvector of
$\Hess V(\s)+B^{*}$  associated with $\mu$ and  $\xi':=P^{-1}\xi$.
To this end, note first that  it obviously holds
\begin{equation}
\label{eq.xi-perp}
\forall\,x\,\in\,(\xi')^{\perp}\,,\ \ 
\big(\,I_{d}+2\,|\mu|D^{-1}\,\xi'\,\xi'^{*}\,\big)\,x\ =\ x\,.
\end{equation} 
Moreover, since  $D^{-1}\xi' = \frac{1}{\mu}(I_{d}-P^{-1}\,J  \,P)\xi'$, it also holds
\begin{align}
\nonumber
\big(\,I_{d}+2\,|\mu|D^{-1}\,\xi'\,\xi'^{*}\,\big)\,\xi'&\ =\ 
\xi'+2\,|\mu|D^{-1}\,\xi'\\
\label{eq.xi-perp'}
&\ =\ -\xi'+2P^{-1}\,J  \,P\xi'\,.
\end{align}
Since $P^{-1}\,J  \,P\xi'$ belongs to $(\xi')^{\perp}$, we deduce 
\eqref{eq.det'} and then \eqref{eq.det} from \eqref{eq.xi-perp}
and \eqref{eq.xi-perp'}.\\[0.2cm]
{\bf Step 4 :} To conclude the proof of the second item of Lemma~\ref{le.B}, it only remains to show that the real symmetric matrix
$M_{V} := \Hess V(\s)+2|\mu|\,\xi\,\xi^{*}$ is positive definite,
where we recall that $\xi$
denotes a unitary eigenvector of
$\Hess V(\s)+B^{*}$  associated with $\mu$.
This is an easy consequence of the Max-Min principle and of 
the relation $\det M_{V} = -\det D>0 $ obtained in the previous step. We have
indeed, defining again $\xi':=P^{-1}\xi$, 
\begin{align*}
\nonumber
\forall\,x\,\in\, \big((1,0,\dots,0)^{*}\big)^{\perp}\ ,\ \ 
\lp (D+2|\mu|\,\xi'\,\xi'^{*})x,x\rp&\ =\ \lp Dx,x\rp+2|\mu|\,\lp \xi,x\rp^{2}\\
&\ \geq\ \mu_{2}\,\|x\|^{2}\,,
\end{align*}
which implies that the second eigenvalue of $D+2|\mu|\,\xi'\,\xi'^{*}$,
that is the second eigenvalue of $M_{V}$, is greater than or equal to $\mu_{2}$, and hence positive. 
The first eigenvalue of $M_{V}$ is then positive according to $\det M_{V}>0$. This concludes
this step of the proof.\\[0.2cm]
{\bf Step 5 :} We now prove the third item of Lemma~\ref{le.B}. Since
$\Hess V(\s)(I_{d}-J)\xi = \mu \xi$ and $J^{*}=-J$, 
it first holds
\begin{equation}
\label{eq.Hess-1}
(\Hess V(\s))^{-1}\,\xi\ =\ \frac1\mu\,(I_{d}-J)\xi
\ \ \text{and then}\ \ 
\lp(\Hess V(\s))^{-1}\,\xi,\xi\rp\ =\ \frac1\mu\,,
\end{equation}
which proves the second part of the third item of Lemma~\ref{le.B}.
Defining again $\xi':=P^{-1}\xi$, this also means
$$
-\frac1{\mu_{1}}+
\sum_{k=2}^{d}(\frac1{\mu_{k}}+\frac1{\mu_{1}})\xi'^{2}_{k}
\ =\ -\frac1{\mu_{1}}\xi'^{2}_{1}+\sum_{k=2}^{d}\frac1{\mu_{k}}\xi'^{2}_{k}
\ =\ \lp D^{-1}\,\xi',\xi'\rp\ =\ \frac1\mu\,.
$$
This implies that $\frac1\mu\geq -\frac1{\lambda_{1}}$, i.e. that
$|\mu|\geq \mu_{1}$, with equality if and only if $\xi'=\pm(1,0,\dots,0)^{*}$, that is if and only if
$\xi$ is a unitary eigenvector of $(\Hess V(\s))^{-1}$ associated with $-\frac1{\mu_{1}}$,
which is equivalent to the relation
$J\xi=0$ by \eqref{eq.Hess-1}, and hence to $B^{*}\xi=0$ since $J=-(\Hess V(\s))^{-1}B^{*}$.

\section{Spectral analysis in the case of Morse functions}
\subsection{Construction of accurate quasimodes}

In the following, we assume that Assumption~\ref{hyp:gener} is satisfied.
Let us then consider some arbitrary  $\m\in \mathcal U^{(0)}\setminus\{ \underline{\m}\}$, that is,
according to Assumption~\ref{hyp:gener}, a local minimum of $V$ which is not the global 
minimum $\underline\m$ of $V$. 
According to  the labelling procedure
of Section~\ref{sub.label}
 leading to the definitions 
\eqref{eq.E}--\eqref{eq.S}, it holds in particular $\m=\m_{i,j}$
and
 $\bsigma(\m)=\sigma_i$
for some $i\in\{2,\dots,N\}$ and $j\in\{1,\dots,N_{i}\}$. 
For every $\s \in {\bf j}(\m)$ and $\rho, \delta >0$, where we recall that the mapping
${\bf j}: \mathcal U^{(0)}\to \mathcal P(\vvv^{(1)}\cup \{\s_{1}\})$ has been defined in \eqref{eq.j}
and that $V(\s)=\bsigma(\m)$, we define the set
$$
     \mc B_{\s,  \rho,\delta}   \  :=   \    \left\{  V\leq \bsigma(\m)+\delta \right\}\cap\left\{x\in \R^{d}\,,\ |\xi(\s)  \cdot (x-\s)| \leq \rho   \right \}
$$     
and the set $\mc C_{\s,  \rho,\delta}$ by:
\be
\label{eq.C-s}
\mc C_{\s,  \rho,\delta}\quad \text{is the connected component of $ \mc B_{\s,  \rho,\delta}$ containing $\s$}\,,
\ee
where $\xi(\s)  $ has been defined in Lemma~\ref{le.B}.
We recall that $\xi(\s)  $ is an unitary 
 eigenvector of the matrix $\Hess V(\s)+B^{*}(\s)$ associated with its only negative eigenvalue
$\mu(\s)$ which has geometric multiplicity one. 
Let us also define
\be
\label{eq.E-m}
E_{\m,  \rho, \delta}\  :=\  \  \big( E_{-}(\m)\cap \{V< \bsigma(\m) + \delta\}\big)    \setminus  \cup_{\s\in {\bf j}(\m)}\mc C_{  \s,\rho,\delta}\,,
\ee
where 
\be
\label{eq.E-}
E_{-}(\m)\ \ \ \text{is the connected component of $\{V<\sigma_{i-1}\}$ 
containing $\m$}.
\ee
According to Assumption~\ref{hyp:gener} and Remark~\ref{rem:comm-2compo}, we recall that
there is precisely one connected component $\widehat E(\m)\neq E(\m)$ of $\{V<\bsigma(\m)\}$
such that $\overline{E(\m)}\cap \overline{\widehat E(\m)}\neq \emptyset$ (see examples in Figure~\ref{fig:association'}). Moreover, it holds
${\bf j}(\m)=\pa \widehat E(\m) \cap\pa E(\m)$ and the global minimum
$\hat \m$ of $V|_{\widehat E(\m)}$ satisfies $\bsigma(\hat\m)>\bsigma(\m)$ and $V(\hat \m)< V(\m)$ (see in this connection \cite{Mi19},  where the notation $\widehat E(\m)$
is introduced for an arbitrary Morse function).\medskip

\noindent
 According
 to the geometry of the Morse function $V$ around 
$\pa E(\m)$ and to Lemma~\ref{le.B}, we have then the following result.
 \begin{lemma}
 \label{le.sep-m0-m1}
Assume  that Assumption~\ref{hyp:gener} is satisfied and
let $\m\in \mathcal U^{(0)}\setminus\{\underline{\m}\}$,
 $\s\in {\bf j}(\m)$, and $\xi(\s)  $
be some unitary 
 eigenvector of the matrix $\Hess V(\s)+B^{*}(\s)$ associated with its unique negative eigenvalue (see Lemma~\ref{le.B}). 
Then, there exists a neigborhood $\mathcal O$ of $\s$ such that:
$$
\forall\, x\in \, \mathcal O\setminus\{\s\}, \ \  \left( x-\s \in \xi(\s)^{\perp}\,\Longrightarrow
V(x)>V(\s)\right).
$$
It follows that
there exist $\rho_0, \delta_0>0$ sufficiently small such that 
for all $\rho\in(0,\rho_{0}]$ and $\delta\in(0,\delta_{0}]$,
 the set $E_{\m,  3\rho, 3\delta}$ defined in \eqref{eq.E-m} has exactly 
two connected components, $E_{  \m,3\rho_0, 3\delta_0}^{+}$
and $E^{-}_{ \m, 3\rho, 3\delta}$,  containing respectively $\m$ and  $\hat\m$.
\end{lemma}

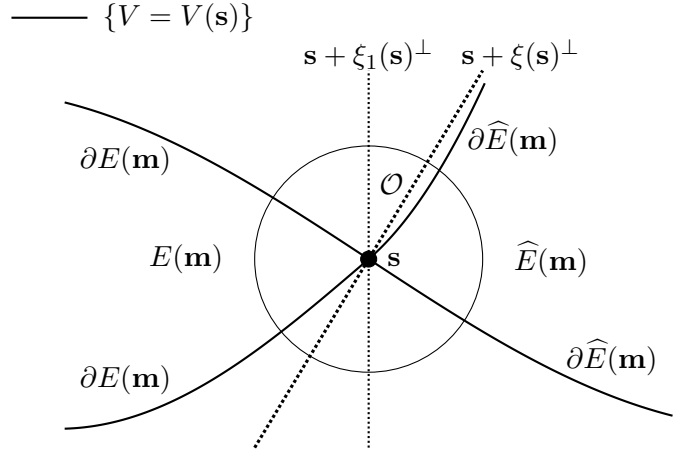
\begin{figure}[h!]
\begin{center}
\begin{tikzpicture}
\tikzstyle{vertex}=[draw,circle,fill=black,minimum size=6pt,inner sep=0pt]
\draw[thick] (-4.7,3.2)-- (-3.7,3.2);
\draw (-2.5, 3.2) node[]{$\{V=V(\s)\}$};
\draw[thick, domain=-4:0, samples=200] plot ({\x}, { 9/4*sin(\x/2.6 r)}) ;
\draw[thick, domain=0:1.53, samples=200] plot ({\x}, { 9/(8*2.6)*((\x+1)*(\x+1)-1)}) ;
\draw[thick, domain=-4:4, samples=200] plot ({\x}, { -9/4*sin(\x/3.4 r)}) ;
\draw[densely dotted, thick] (0,-2.5)--(0,2.5);
\draw[densely dotted, very thick] (-1.5,-2.5)--(1.5,2.5);
 \draw (0,2.7) node[]{$\s + \xi_{1}(\s)^{\perp}$};
 \draw (2,2.7) node[]{$\s + \xi(\s)^{\perp}$};
  \draw (0,0)  node[vertex,label=east: {$\s$}](v){};
      \draw (-3.2,1.3) node[]{$\pa  E(\m)$};
         \draw (-3.2,-1.6) node[]{$\pa  E(\m)$};
  \draw (-2.4,0) node[]{$ E(\m)$};
      \draw (3.2,-1.3) node[]{$\pa \widehat  E(\m)$};
         \draw (1.9,1.6) node[]{$\pa  \widehat E(\m)$};
  \draw (2.4,0) node[]{$ \widehat E(\m)$};
   \draw (0,0) circle(1.5);
   \draw (0.3 , 1 ) node[]{$\mathcal O$};
 \end{tikzpicture}
\caption{Representation of the Morse function $V$ near $\s\in{\bf j}(\m)$. 
Here, $\xi_{1}(\s)$ denotes an eigenvector of $\Hess V(\s)$ associated with its negative eigenvalue
and $B^{*}(\s)\xi(\s)\neq0$. Note that
according to the last item in Lemma~\ref{le.B},
$\s + \xi_{1}(\s)^{\perp}$ and $\s + \xi(\s)^{\perp}$ coincide if and only if
$B^{*}(\s)\xi(\s)=0$.}
 \label{fig:ssp_cutoff}
 \end{center}
\end{figure}

\begin{proof} For shortness, we denote $\xi=\xi(\s)$.
By a continuity argument, note that to prove the first part of Lemma~\ref{le.sep-m0-m1},
it is sufficient to prove that the linear hyperplane $\xi^{\perp}$
does not meet the cone $\{X\in\R^{n}\,;\ \lp \Hess V(\s) X,X\rp\leq0\}$
outside the origin.
The second part of the lemma  then simply  follows from the observation that
the set $\mc C_{\s,  \rho,\delta}$ defined in \eqref{eq.C-s}
is thus an arbitrary small neighborhood of $\s$ when $ \rho,\delta>0$ tend to $0$.
\\[0.1cm]
When $d\geq 3$,  it is then enough to show that
for any  column vector $X\in\R^{d}\setminus\{0\}$ such that $\lp \Hess V(\s) X,X\rp=0$,
it holds ${\rm Span\,}X\oplus \xi^{\perp}=\R^{d}$, i.e. $\lp X,\xi\rp\neq0$.
Indeed, when $d\geq 3$, any linear hyperplane meets $\{X\in\R^{n}\,;\ \lp \Hess V(\s) X,X\rp>0\}$
and then meets $\{X\in\R^{d}\setminus\{0\}\,;\ \lp \Hess V(\s) X,X\rp=0\}$ if and only if it meets
$\{X\in\R^{d}\setminus\{0\}\,;\ \lp \Hess V(\s) X,X\rp\leq0\}$.
Let us then consider $X\in\R^{d}\setminus\{0\}$ such that $\lp \Hess V(\s) X,X\rp=0$
and let us prove that $\lp X,\xi\rp\neq0$.
To show this, let us work in orthonormal coordinates  of $\R^{d}$
where 
$\Hess V(\s)$ is diagonal, i.e. where $\Hess V(\s)={\rm Diag}(-\mu_{1},\mu_{2},\dots,\mu_{d})$.
It then follows from $\lp \Hess V(\s) X,X\rp=0$ and from the third item of Lemma~\ref{le.B} that
\begin{equation*}
\mu_{1}X_{1}^{2}\ =\ \sum_{k=2}^{d}\mu_{k}X_{k}^{2}\quad
\text{and}\quad\frac{1}{\mu_{1}}\xi_{1}^{2}\ >\ \sum_{k=2}^{d}\frac{1}{\mu_{k}}\xi_{k}^{2}\ \geq\ 0\,.
\end{equation*}
It holds in particular $X_{1}\neq 0$ and thus, by multiplying the two above relations,
\begin{align*}
|\xi_{1}\,X_{1}|&\ >\ \big(\sum_{k=2}^{d}\frac{1}{\mu_{k}}\xi_{k}^{2}\big)^{\frac12}\big(\sum_{k=2}^{d}\mu_{k}X_{k}^{2}\big)^{\frac12}
\ \geq\ |\sum_{k=2}^{d}\xi_{k}\,X_{k}|\,,
\end{align*}
the last inequality resulting from the Cauchy-Schwarz inequality.
The relation $\lp X,\xi\rp\neq0$ follows.\\[0.1cm]
When $d=2$, the situation is slightly different since for any hyperplane $H$,
either $H\setminus\{0\}\subset \{X\in\R^{2}\setminus\{0\}\,;\ \lp \Hess V(\s) X,X\rp\leq0\}$
or $H\setminus\{0\}\subset \{X\in\R^{2}\setminus\{0\}\,;\ \lp \Hess V(\s) X,X\rp>0\}$.
Take again orthonormal coordinates where $\Hess V(\s)={\rm Diag}(-\mu_{1},\mu_{2})$.
We have then only to prove that the vector $\xi':=(-\xi_{2},\xi_{1})^{*}$, which spans $\xi^{\perp}$, satisfies 
$$
-\mu_{1}\xi^{2}_{2}+\mu_{2}\xi^{2}_{1}\ =\ \lp \Hess V(\s)\xi',\xi'  \rp\ >\ 0\,.
$$
This is obviously satisfied since equivalent to 
$$
0\ >\ \frac1{\mu_{2}}\xi^{2}_{2}-\frac1{\mu_{2}}\xi^{2}_{1}\ =\ \lp (\Hess V(\s))^{-1}\xi,\xi  \rp\,,
$$
which holds true thanks to iii) of Lemma \ref{le.B}.
This concludes the proof of Lemma~\ref{le.sep-m0-m1}.
\end{proof}

\begin{figure}
 \begin{center}
\scalebox{0.8}{ 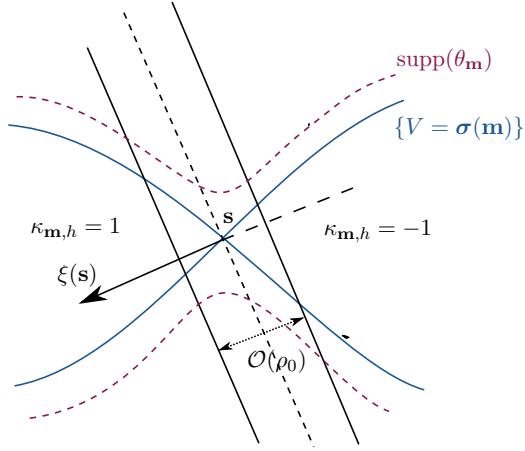}
  \end{center}
  \caption{The support of the function $\kappa_{\m,h}$}
  \label{fig:ssp_cutoff'}
  \end{figure}

\noindent
Let us now define, for every $h\in(0,1]$ and for every $\rho_{0},\delta_{0}>0$ small enough, the function $ \kappa_{\m,h}$   on the sublevel set  $E_-(\m)\cap\{V< \bsigma(\m) + 3\delta_0\}$ (see \eqref{eq.E-}) as follows:
\begin{enumerate}
\item[1.] On  the disjoint open sets $E^{+}_{  \m,3\rho_0,3\delta_0} $ and $E^{-}_{  \m,3\rho_0,3\delta_0} $
introduced in Lemma~\ref{le.sep-m0-m1},
\be
\label{de.kappa1}     \kappa_{\m,h}(x)     \    :=    \     \begin{cases}      +1      &    \ \ \text{ for }      x\in
E^{+}_{  \m,3\rho_0,3\delta_0}         \\
-1       &  \ \     \text{ for }        x\in
E^{-}_{  \m,3\rho_0,3\delta_0}\end{cases}\,.
\ee
\item[2.] For every $\s\in {\bf j}(\m)$ and $x\in \mc C_{  \s,3\rho_0,3\delta_0}  \cap \{V< \bsigma(\m) + 3\delta_0\}$
(see \eqref{eq.C-s}),
\begin{equation}
\label{de.kappa2} 
 \kappa_{\m,h}(x)\ :=\ C^{-1}_{\s,h}  \int_0^{\xi(\s)\cdot(x-\s)}   \,  \chi(\rho_0^{-1} \eta) \, e^{- \frac{  |\mu(\s)| \eta^2}{2h}}   \, d\eta  
 \,,    
\end{equation}
where 
the orientation of $\xi(\s)$ is chosen in such a way that there exists a  neighborhood
$\mathcal O$ of $\s$  such that 
$E(\m)\cap \mathcal O$
is included in the half-space $\{\xi(\s)  \cdot (x-\s)>0\}$
(see Lemma~\ref{le.sep-m0-m1} and Figures~\ref{fig:ssp_cutoff} and~\ref{fig:ssp_cutoff'}),
$\chi\in C^\infty(\mb R;[0,1])$ is even and satisfies $\chi\equiv 1$ on $[-1,1]$, $\chi(\eta)= 0$ 
 for $|\eta|\geq 2$, and
 \begin{equation*}
\label{de.kappa3} 
     C_{\s,h}     \    :=     \     \frac 12    \  \int_{-\infty}^{+\infty} \chi(\rho_0^{-1} \eta) \, e^{- \frac{ |\mu(\s)| \eta^2}{2h}}   \, d\eta    \,. 
 \end{equation*}
 \end{enumerate}
 
 \noindent
Note in particular that  
\be
\label{de.kappa4} 
 \exists \gamma>0   \  \text{ s.t. }      \   \    C^{-1}_{\s,h}      \ =\    \sqrt{  \frac{2|\mu(\s)|}{\pi h}}   \    \left(  1 +   \mc O(e^{-\frac \gamma h}) \right)  \,.  
\ee
Note also that  for every $\rho_{0},\delta_{0}>0$ small enough, thanks to
the definitions \eqref{de.kappa1} and \eqref{de.kappa2},
and since the sets $E^{+}_{  \m,3\rho_0,3\delta_0} $, $E^{-}_{  \m,3\rho_0,3\delta_0} $,
and $\mc C_{  \s,3\rho_0,3\delta_0}$'s, $\s\in {\bf j}(\m)$, are two by two disjoint 
(see  Lemma~\ref{le.sep-m0-m1}), $\kappa_{\m,h}$ is well defined and is $\ccc^\infty$ on $E_-(\m)\cap\{V< \bsigma(\m) + 3\delta_0\}$.\medskip

\noindent
Consider now  a smooth function $\theta_{\m}$ 
such that 
\be
\label{de.kappa5}      \theta_\m(x)    \    :=   \          \begin{cases}      1   &    
\text{ for }      x \in \{V\leq \bsigma(\m) + \frac32\delta_0\} \cap E_{-}(\m)    \\
0      &      \text{ for }     x\in \mb R^d \setminus\big(  \{V< \bsigma(\m) + 2\delta_0\}   \cap E_{-}(\m)\big)    \\
      \end{cases}  \,.
\ee
The function $\theta_\m\kappa_{\m,h}$ then belongs to $\ccc_{c}^{\infty}(\R^{d};[-1,1])$ and 
$$\supp \theta_\m\kappa_{\m,h}\ \subset\  E_{-}(\m)\cap
\{V< \bsigma(\m) + 2\delta_0\}\,.$$

\begin{definition}
\label{de.quasi} 
For any $\m\in\uuu^{(0)}$ let us define the function $\psi_{\m,h}$
 by
$$ 
\psi_{\m,h}(x)     \  :=\ \theta_{\m}(x)        \, \big(\kappa_{\m,h}(x)+1\big)   
$$
when $\m\neq \underline{\m}$ and, when $\m=\underline{\m}$,
$\psi_{\m,h}(x)      := 1$.
We then define,  for any $\m\in\uuu^{(0)}$, the quasimode $\varphi_{\m,h}$
by 
$$
\varphi_{\m,h}(x)\ :=\ \frac{\psi_{\m,h}(x)}{\|\psi_{\m,h}\|_{L^2(m_h)}}\,.
$$
\end{definition}

\noindent Note that, for every $h\in(0,1]$, 
it holds $L_{V,b,\nu} \varphi_{\underline\m,h}=0$  and for every $\m\in \uuu^{(0)}\setminus\{\underline{\m}\}$, the quasimodes $\psi_{\m,h}$ and $\varphi_{\m,h}$
belong to $\ccc_c^\infty(\mb R^d;\R^{+})$ with supports included in $E_{-}(\m)\cap
\{V< \bsigma(\m) + 2\delta_0\}$. 
We have more precisely the following lemma resulting from the previous construction. 

\begin{lemma}\label{lem:orthogQM}
Assume  that Assumption~\ref{hyp:gener} is satisfied.
For every $\m\in \uuu^{(0)}$ and every small $\epsilon>0$ fixed,   there exist $\rho_{0},\delta_{0}>0$ small enough such that  for every $h\in(0,1]$, one has:
\begin{itemize}
\item[i)] It holds $$\supp \psi_{\m,h}\ \subset\  \overline{E(\m)}+B(0,\epsilon)\,.$$
\item[ii)] When $\m\neq \underline\m$, there exists a  neighborhood $\mathcal O_{\rho_{0},\delta_{0}}$ of $\overline{E(\m)}$
such that:
$$
\mathcal O_{\rho_{0},\delta_{0}}\setminus \cup_{\s\in {\bf j}(\m)}\mc C_{  \s,3\rho_{0},3\delta_{0}}\ \subset\ 
\{\theta_{\m}\,\kappa_{\m,h}=1\}\,.
$$ 
In particular, it holds 
$${\rm argmin}_{\supp \psi_{\m,h}} V\ = \ 
{\rm argmin}_{\{\theta_{\m}\,\kappa_{\m,h}=1\}} V  \ =\  {\rm argmin}_{E(\m)} V\ =\ \{\m\}\,.$$
\item[iii)] When $\m\neq \underline\m$, it holds:
$$\forall\,x\,\in\, \supp \nabla\psi_{\m,h}\,,\ \left(  V(x)< \bsigma(\m)+\frac32\delta_{0}\ \Longrightarrow\ \ x\,\in\, \cup_{\s\in {\bf j}(\m)}\mc C_{  \s,3\rho_{0},3\delta_{0}}\right).$$
\end{itemize}
Let moreover  $\m'$ belong to  $\uuu^{(0)}$ with $\m\neq\m'$. The following then hold true for every $\rho_{0},\delta_{0}>0$ small enough and every $h\in(0,1]$:
\begin{itemize}
\item[iv)] if $\bsigma(\m)=\bsigma(\m')$, then $\supp(\psi_{\m,h})\cap \supp(\psi_{\m',h})=\emptyset$,
\item[v)] if $\bsigma(\m)>\bsigma(\m')$, then 
\begin{itemize}
\item[--] either $\supp(\psi_{\m,h})\cap \supp(\psi_{\m',h})=\emptyset$,
\item[--] or $\psi_{\m,h}=2$ on $\supp(\psi_{\m',h})$ and $V(\m')>V(\m)$.
\end{itemize}
\end{itemize}
\end{lemma}

 \noindent
\bp The first part of Lemma~\ref{lem:orthogQM} follows from Assumption~\ref{hyp:gener} and from the construction of 
the quasimodes $\varphi_{\m,h}$ defined in Definition~\ref{de.quasi} 
for $\m\in\mathcal U^{(0)}$, see indeed 
\eqref{de.kappa1}, \eqref{de.kappa2}, and \eqref{de.kappa5}. Let us then prove the second part of Lemma~\ref{lem:orthogQM}.\medskip

\noindent
When  $\bsigma(\m)=\bsigma(\m')$ and $\m\neq \m'$, note first that $\m$ and $\m'$
 differ from $\underline{\m}$ since $\bsigma(\m)=+\infty$ if and only if $\m=\underline{\m}$. When moreover $\m'\notin E_{-}(\m)$,
it holds $E_{-}(\m)\neq E_{-}(\m')$ and hence
$E_{-}(\m)\cap E_{-}(\m') =\emptyset$, implying 
$\supp(\psi_{\m,h})\cap \supp(\psi_{\m',h})=\emptyset$.
In the case when  $\m'\in E_{-}(\m)$, the statement of Lemma~\ref{lem:orthogQM}
follows from 
 ii) of Assumption~\ref{hyp:gener} and of Remark~\ref{rem:comm-2compo},
 which indeed imply that $\overline{E(\m)}\cap \overline{E(\m')}=\emptyset$
 (see the first item of Lemma~\ref{lem:orthogQM}).\medskip
 
 \noindent
When $\bsigma(\m)>\bsigma(\m')$ and  $\m'\notin E(\m)$,
it holds $\overline{E(\m)}\cap \overline{E(\m')}=\emptyset$, and again,
according to the first item of Lemma~\ref{lem:orthogQM}, it holds
$\supp(\psi_{\m,h})\cap \supp(\psi_{\m',h})=\emptyset$ for every $\rho_{0},\delta_{0}>0$ small enough.
Lastly, when $\bsigma(\m)>\bsigma(\m')$ and  $\m'\in E(\m)$, it holds
$\overline{E(\m')}\subset E_{-}(\m')\subset E(\m)$
and then, according to the second item of Lemma~\ref{lem:orthogQM},  $\psi_{\m,h}=2$ on $\supp(\psi_{\m',h})$
for every $\rho_{0},\delta_{0}>0$ small enough.
Besides,
the relation $V(\m')>V(\m)$ follows
from $\m'\in E(\m)$ and from the first item of Assumption~\ref{hyp:gener}.
\ep

\subsection{Quasimodal estimates}

We write in the sequel   $a\eqsim b$ and  $a\lesssim b$ to mean, in the limit $h\to 0,$ equality/inequality up to a multiplicative factor $1 + \mc O(h)$. 
Moreover, we define for shortness, for any critical point $\mathbf{u}$ of $V$:
\begin{equation*}
\label{de.Di} 
D_{\mathbf{u}}     \   :=       \       \sqrt {\vert \det \Hess V(\mathbf{u})\vert}        \   >   \     0.
\end{equation*}


%
%

\begin{proposition}
\label{pr.L2}
Assume that
Assumption~\ref{hyp:gener} is satisfied
and consider the families $\big(\psi_{\m,h}, \m\in\uuu^{(0)}\big)$
and $\big(\varphi_{\m,h}, \m\in\uuu^{(0)}\big)$
of Definition~\ref{de.quasi}.
Then,  for every $\m\in \mathcal U^{(0)}\setminus\{\underline\m\}$ and
$\rho_0,\delta_0>0$ small enough,
it holds in the limit $h\to 0$:
\be \label{L2}
\|\psi_{\m,h}  \|^{2}_{L^2(m_h)}    \ \eqsim\ 
4     \, \frac{D_{\underline \m}}{D_{\m}}  \,  e^{- \frac{V(\m) - V(\underline \m)}h}\,.
     \ee
Moreover,
 there exists $C>0$ such that for every $\m,\m'\in \uuu^{(0)}$, it holds in the limit $h\to 0$:
\be\label{eq:orthogQM}
\<\varphi_{\m,h},\varphi_{\m',h}\>\ =\ \delta_{\m,\m'}+\ooo(e^{-\frac c h}).
\ee
\end{proposition}

\noindent
\bp To prove the relation \eqref{L2},  write, according to Definition~\ref{de.quasi}, 
\begin{gather*}
\|\psi_{\m,h}  \|^{2}_{L^2(m_h)}    \    =    \      Z^{-1}_h   \,    \int     \big(\theta_\m (\kappa_{\m,h}+1)\big)^2  \, e^{-\frac{V(x)}{h}} \, dx  \,,       
\end{gather*}
where
$Z_h$ is the normalizing constant defined by \eqref{eq:defmesure}. Hence, 
according to Lemma~\ref{lem:orthogQM} and standard tail estimates and Laplace asymptotics, we get,
in the limit $h\to 0$,
\begin{gather*}
\label{eq.Zbeta}
     Z_h        \  \eqsim    \  
(2\pi h )^{\frac d2}  \,    D_{\um}^{-1}      \, e^{- \frac{V(\um)}{h}} 
   \end{gather*}
as well as
 \begin{align}
 \nonumber
     \int     \big(\theta_{\m} (\kappa_{\m,h}+1)\big)^2  \, e^{-\frac{V(x)}{h}} \, dx     
        &\ \eqsim\ 
   4\, (2\pi h )^{\frac d2}  \, 
      D_{\m}^{-1}
      \, e^{- \frac{V(\m)}{h}} \,.
\end{align}
The estimate \eqref{L2} 
then follows easily.\medskip

\noindent
Let us now prove the relation \eqref{eq:orthogQM}.
According to Definition~\ref{de.quasi}, 
note first that
$\<\varphi_{\m,h},\varphi_{\m,h}\>=1$ for every $\m\in\uuu^{(0)}$.
Moreover, when $\m,\m'\in\uuu^{(0)}$ and $\m\neq\m'$,
it follows from Lemma \ref{lem:orthogQM} that, up to switching $\m$ and $\m'$, we are in one of the two following cases:
\begin{itemize}
\item[--] either $\supp(\varphi_{\m,h})\cap\supp(\varphi_{\m',h})=\emptyset$, and then
$$
\<\varphi_{\m,h},\varphi_{\m',h}\>\ =\ 0\,,
$$
\item[--] 
or $\psi_{\m,h}=2$ on $\supp(\psi_{\m',h})$ and $V(\m')>V(\m)$, and then, using the preceding estimates,
\begin{align*}
\<\varphi_{\m,h},\varphi_{\m',h}\>&\ =\ \frac{2}{\|\psi_{\m,h}  \|_{L^2(m_h)} }\int_{\supp 
\psi_{\m',h}} \frac{\psi_{\m',h} }{\|\psi_{\m',h}  \|_{L^2(m_h)}}\,\frac{e^{-\frac{V(x)}{h}}}{Z_{h}} \, dx\\
&\ =\ \frac{1}{\|\psi_{\m,h}  \|_{L^2(m_h)}\|\psi_{\m',h}  \|_{L^2(m_h)} }\mathcal O\big(e^{-\frac{V(\m')-V(\underline\m)}{h}}\big)\ =\ \mathcal O\big(e^{-\frac{C}{h}}\big)\,,
\end{align*}
where $C=\frac{V(\m')-V(\m)}{2}>0$.
\end{itemize}
This leads to \eqref{eq:orthogQM}.
\ep
\
\\

 \begin{proposition} 
\label{pr.quad-form}  For every $\m\in\ulu^{(0)}$ and  $\rho_0,\delta_0>0$  small enough, it holds
in the limit $h\to 0$:
\begin{equation}\label{eq:estim:dirich}
\lp L_{V,b,\nu}\psi_{\m,h} , \psi_{\m,h}   \rp_{L^2(m_h)}   \    \eqsim    \         \sum_{\s \in {\bf j}(\m)} \frac{2 |\mu(\s)|}{\pi}   \,    \frac{ D_{\um}}{D_{\s} }  \,  e^{-\frac{V(\s)-V(\um)}{h}}       
    \end{equation}
    and then 
    \begin{equation}\label{eq:estim:dirichstar}
\lp L_{V,b,\nu}\varphi_{\m,h} , \varphi_{\m,h}   \rp_{L^2(m_h)}   \    \eqsim    \ 
     \sum_{\s \in {\bf j}(\m)} \frac{|\mu(\s)|}{2\pi}   \,    \frac{ D_{\m}}{D_{\s} }  \,  e^{-\frac{V(\s)-V(\m)}{h}} .
    \end{equation}
 \end{proposition}

\noindent
\bp
Note first that thanks to \eqref{eq:propvectfield}, one has $\div(b_h m_h)=0$ and hence:
$$
\forall\, u\,\in\,\ccc_c^\infty(\R^d;\R)\,,\ \ \<b_h\cdot\nabla u,u\>_{L^{2}(m_h)}\ =\ -\frac 12\int u^2\div(b_h m_h) dx\ =\ 0\,.
$$
Using this relation together with \eqref{eq.L}, \eqref{de.kappa1}--\eqref{de.kappa5}, Definition~\ref{de.quasi}, and Lemma~\ref{lem:orthogQM}, we get, in the limit $h\to 0$,
\begin{align}
\nonumber
\lp L_{V,b,\nu}\psi_{\m,h} , \psi_{\m,h}   \rp_{L^2(m_h)}  &\ =\ 
\lp(-h\Delta + \nabla V\cdot \nabla)\psi_{\m,h} , \psi_{\m,h}   \rp_{L^2(m_h)}
\\
\nonumber
&\ =\ Z_{h}^{-1}h\int |\nabla\big(\theta_{\m} (\kappa_{\m,h}+1)\big)|^2   \  e^{-\frac V h}   \, dx\\
\nonumber
&\ =\ Z_{h}^{-1}h\int \theta_{\m}^{2}|\nabla \kappa_{\m,h}|^2   \  e^{-\frac V h}   \, dx
\,+\,Z_{h}^{-1}\mathcal O(e^{-\frac{\bsigma(\m)+\delta_0}{h}})\\
\nonumber
&\ =\ 
Z^{-1}_{h}\mathcal O(e^{-\frac{\bsigma(\m)+\delta_0}{h}})\\
\label{eq.Dir1}
+Z^{-1}_{h}\sum_{\s\in {\bf j}(\m)} &C_{\s,h}^{-2}h\int_{\mc C_{  \s,3\rho_0,3\delta_0}} \!\!\!\!\!\!\!\theta_{\m}^{2}(x)\chi^{2}(\rho_0^{-1} \xi\cdot(x-\s))
e^{ - \frac{|\mu| (\xi\cdot(x-\s))^{2}}h}   \  e^{-\frac V h}   \, dx\,,
\end{align}
where for short we denote $\xi=\xi(\s)$ and $\mu=\mu(\s)$.
From the second item in Lemma~\ref{le.B} and the Taylor expansion of $V+|\mu| \lp \xi,\cdot-\s\rp^{2}$ around $\s \in{\bf j}(\m)$,
\begin{align*}
V(x)+|\mu| (\xi\cdot(x-\s))^{2}\ =\ V(\s)&+\frac12\lp\Hess V(\s)\,(x-\s),x-\s\rp 
\\&+|\mu| \lp\xi\xi^{*}
(x-\s),x-\s\rp+\mathcal O(|x-\s|^{3})\,,
\end{align*}
it is clear  that for $\rho_{0}$ and $\delta_{0}$ small enough,
$V+|\mu| \lp \xi,\cdot-\s\rp^{2}$ uniquely attains its minimal value in $\mc C_{\s,  3\rho_0,3\delta_0}$
at $\s$ since:
$$
\nabla \big(\,V+|\mu| \lp \xi,\cdot-\s\rp^{2}\,\big)(\s)=0\quad
\text{and}\quad \Hess\big(\,V+|\mu| \lp \xi,\cdot-\s\rp^{2}\,\big)(\s)=M_{V}\,.$$

 \noindent
Moreover, using again the second item in Lemma~\ref{le.B} and a standard Laplace method, it holds in the limit $h\to 0$, for every $\s\in{\bf j}(\m)$,
\begin{align}
\nonumber
\!\!\!\!\! C_{\s,h}^{-2}\int_{\mc C_{ \s, 3\rho_0,3\delta_0}} \!\!\!\!\!\!\!\!\!\!\theta_{\m}^{2}\chi^{2}(\rho_0^{-1} \lp\xi,\cdot-\s\rp) \, e^{ - \frac{|\mu| \lp\xi,\cdot-\s\rp^2}h}   \  e^{-\frac V h}   \, dx
&\ \eqsim\  \frac{(2\pi h )^{\frac d2}}{C_{\s,h}^{2}D_{\s}}\, e^{- \frac {V(\s)} h}\\
\label{eq.Dir2}
&\ \eqsim\ \frac{2\,(2\pi h )^{\frac d2}\,|\mu|}{\pi\,h\,D_{\s}}\, e^{- \frac{V(\s)} h}\,,
\end{align}
where we also used \eqref{de.kappa4}  at the last line.
The statement of Proposition~\ref{pr.quad-form} then follows from \eqref{eq.Dir1} and \eqref{eq.Dir2}, using also
$
Z_{h} \eqsim (2\pi h )^{\frac d2}  D_{\um}^{-1}           \, e^{- \frac{V(\um)}{h}}
$.
\ep

 \begin{proposition}\label{pr.quad-1-form} Let $\m\in\ulu^{(0)}$. 
For $\rho_{0}$ and $\delta_{0}$ sufficiently small, it holds in the limit $h\to 0$:
\begin{equation}\label{quasimodal2}
   \| L_{V,b,\nu}\psi_{\m,h} \|^{2}_{L^2(m_h)}    \    =   \ 
  \lp L_{V,b,\nu}\psi_{\m,h} , \psi_{\m,h}   \rp_{L^2(m_h)} \, \mathcal O(h)     \,.   
\end{equation}
and
\begin{equation}\label{quasimodal2star}
\| L^{*}_{V,b,\nu}\psi_{\m,h} \|^{2}_{L^2(m_h)}       \    =   \ 
  \lp L_{V,b,\nu}\psi_{\m,h} , \psi_{\m,h}   \rp_{L^2(m_h)} \, \mathcal O(1)     \,.   
\end{equation}
\end{proposition}

\noindent
\bp Let $\s\in{\bf j}(\m)$ and denote for short $\xi=\xi(\s)$ and $\mu=\mu(\s)$.  We first recall  the Taylor expansion of $V+|\mu| \lp \xi,\cdot-\s\rp^{2} $ around $\s$, 
\begin{align*}
V(x)+|\mu| (\xi\cdot(x-\s))^{2}
&\ =\ V(\s)+\frac12\lp M_{V}\,(x-\s),x-\s\rp +\mathcal O(|x-\s|^{3})\,,
\end{align*}
which implies, according to the second item of Lemma~\ref{le.B},
 that for $\rho_{0}$ and $\delta_{0}$ small enough:
 \begin{itemize}
 \item[--] $\nabla \big(\,V+|\mu| \lp \xi,\cdot-\s\rp^{2}\,\big)(\s)=0$,
 \item[--]
$V+|\mu| \lp \xi,\cdot-\s\rp^{2}$ uniquely attains its minimal value in $\mc C_{\s,  3\rho_0,3\delta_0}$
at $\s$.
\end{itemize} 

\noindent
Note now that according to \eqref{eq.L}, it holds
$$
L_{V,b,\nu}\,\psi_{\m,h}\ =\  \theta_{\m}\,L_{V,h}\,\kappa_{\m,h}
+
\big(1+\kappa_{\m,h}\big)\,L_{V,h}\, \theta_{\m} -2h\nabla\kappa_{\m,h}\cdot\nabla\theta_{\m}\,,
$$
with on  $\mc C_{  \s,3\rho_0,3\delta_0}$, for every $\s\in{\bf j}(\m)$,  according to  \eqref{de.kappa2},
\begin{align*}
L_{V,b,\nu}\,\kappa_{\m,h}&\ =\ -h\Delta  \kappa_{\m,h}+ \nabla V\cdot \nabla \kappa_{\m,h} + 
b_h  \cdot\nabla \kappa_{\m,h}\\
&\ =\ C_{\s,h}^{-1}\chi(\rho_{0}^{-1}\lp\xi,\cdot-\s\rp) \, e^{ - \frac{|\mu| \lp\xi,\cdot-\s\rp^2}{2h}}\,\big(\,\nabla V\cdot \xi+b_h  \cdot \xi
+|\mu|\lp\xi,\cdot-\s\rp\,\big)\\
&\qquad -h\,C_{\s,h}^{-1}\,\div\big(\,\chi(\rho_{0}^{-1}\lp\xi,\cdot-\s\rp) \,\xi\,\big) \, e^{- \frac{  |\mu| \lp\xi,\cdot-\s\rp^2}{2h}}\,,
\end{align*}
where we recall that $b_h=b+h\nu$.
It then follows from \eqref{de.kappa1}--\eqref{de.kappa5}  that
in the limit $h\to 0$,
\begin{align*}
\| L_{V,b,\nu}\psi_{\m,h} \|^{2}_{L^2(m_h)} 
&\ =\ 
\sum_{\s\in {\bf j}(\m)}
\| {\bf 1}_{\mc C_{ \s, 3\rho_0,3\delta_0}} L_{V,b,\nu}\psi_{\m,h} \|^{2}_{L^2(m_h)} 
+\,\frac{\mathcal O(e^{- \frac{\bsigma(\m)+\delta_{0}}h   })}{Z_{h}}\\
&\ =\ 
\sum_{\s\in {\bf j}(\m)}\frac{C_{\s,h}^{-2}}{Z_{h}}\int_{\mc C_{ \s, 3\rho_0,3\delta_0}}   \chi^{2}(\rho_{0}^{-1}\xi\cdot(x-\s)) \, e^{- \frac{V+ |\mu| (\xi\cdot(x-\s))^2}h}\\
&\qquad\qquad\times\big(\,\nabla V\cdot \xi+b \cdot \xi
+|\mu|\xi\cdot(x-\s)\,+h\nu\cdot\xi\big)^{2}\,dx\\
& +\  \frac{\mathcal O(e^{-\frac{\bsigma(\m)+c}h})}{Z_{h}}
\end{align*}
for some real constant $c\in(0,\delta_{0})$.
Moreover,
using $b(\s)=0$ and  the first item of Lemma~\ref{le.B},
the Taylor expansion of $\nabla V+b  $  around $\s$
satisfies
\begin{align*}
(\nabla V+b  )\cdot \xi
+|\mu|\xi\cdot(x-\s)&\ =\ \lp(\Hess V(\s)+B  ) (x-\s),\xi \rp +|\mu|\xi\cdot(x-\s)\\&\phantom{**************}+\mathcal O((x-\s)^{2})\\
&\ =\ \mu \xi\cdot(x-\s)+|\mu|\xi\cdot(x-\s)+\mathcal O((x-\s)^{2})\\
&\ =\ \mathcal O((x-\s)^{2})\,.
\end{align*}
It then follows from Proposition~\ref{pr.quad-form}, standard tail estimates, and Laplace asymptotics, that
in the limit $h\to 0$,
\begin{align*}
\| L_{V,b,\nu}\psi_{\m,h} \|^{2}_{L^2(m_h)}
&\ =\ \sum_{\s\in {\bf j}(\m)}\frac{C_{\s,h}^{-2}}{Z_{h}}
\int_{\mc C_{  \s,3\rho_0,3\delta_0}} \!\!\!\!\!\!\mathcal O((x-\s)^{4}+h^2) \, e^{-\frac{V+ |\mu| (\xi\cdot(x-\s))^2}h}\,dx\\
&\phantom{******}+\  \frac{\mathcal O(e^{-\frac {\bsigma(\m)+c}h})}{Z_{h}}\\
&\ =\ \lp L_{V,b,\nu}\psi_{\m,h} , \psi_{\m,h}   \rp_{L^2(m_h)} \, \mathcal O(h)\,,
\end{align*}
which proves \eqref{quasimodal2}.\medskip

\noindent
To prove \eqref{quasimodal2star}, we observe that since $L_{V,b,\nu}^*=L_{V,-b,-\nu}$, the same computation as above shows that 
in the limit $h\to 0$,
\begin{align*}
\| L^{*}_{V,b,\nu}\psi_{\m,h} \|^{2}_{L^2(m_h)}
&\ =\ 
\sum_{\s\in {\bf j}(\m)}\frac{C_{\s,h}^{-2}}{Z_{h}}\int_{\mc C_{ \s, 3\rho_0,3\delta_0}}   \chi^{2}(\rho_{0}^{-1}\xi\cdot(x-\s)) \, e^{- \frac{V+ |\mu| (\xi\cdot(x-\s))^2}h}\\
&\times\big(\,\nabla V\cdot \xi-b  \cdot \xi
+|\mu|\xi\cdot(x-\s)-h\nu\cdot\xi\,\big)^{2}\,dx\\
& \phantom{**************}+\  \frac{\mathcal O(e^{-\frac{\bsigma(\m)+c}h})}{Z_{h}}.
\end{align*}
However, contrary to the preceding case, one has here only 
$$
\nabla V\cdot \xi-b  \cdot \xi
+|\mu|\xi\cdot(x-\s)=\ooo(x-\s)\,,
$$
which implies, in the limit $h\to 0$,
\begin{align*}
\| L^{*}_{V,b,\nu}\psi_{\m,h} \|^{2}_{L^2(m_h)}
&\ =\ \sum_{\s\in {\bf j}(\m)}\frac{C_{\s,h}^{-2}}{Z_{h}}
\int_{\mc C_{  \s,3\rho_0,3\delta_0}} \!\!\!\!\!\!\mathcal O((x-\s)^{2}+h^2) \, e^{-\frac{V+ |\mu| (\xi\cdot(x-\s))^2}h}\,dx\\
&\phantom{******}+\  \frac{\mathcal O(e^{-\frac {V(\s)+c}h})}{Z_{h}}\\
&\ =\ \lp L_{V,b,\nu}\psi_{\m,h} , \psi_{\m,h}   \rp_{L^2(m_h)} \ooo(1)\,,
\end{align*}
which is exactly \eqref{quasimodal2star}.
\ep

\subsection{Proof of Theorem~\ref{th.main}}

Throughout this section, we denote for shortness
$$\lp\cdot,\cdot \rp=\lp\cdot,\cdot \rp_{L^{2}(m_{h})}\,,\ \ \ 
\|\cdot \|=\|\cdot \|_{L^{2}(m_{h})}\,,\ \ \ 
L_{V,b,\nu}=L_V\,,$$
and we label the local minima $\m_1,\ldots,\m_{n_0}$ of $V$ in so that
$(S(\m_{j}))_{j\in\{1,\dots, n_{0}\}}$ is non-increasing (see \eqref{eq.S}):
$$S(\m_1)=+\infty \ \ \text{and, for all $j \in\{2,\dots,n_{0}\}$,}\ \  S(\m_{j+1})\leq S(\m_j)<+\infty\,.$$
For all $j\in\{1,\ldots,n_0\}$, we will also denote  for shortness
$$
S_j\ :=\ S(\m_j)\,,\  \varphi_j\ :=\ \varphi_{\m_j,h}\,,\ \ \text{and}\ \ 
\tilde\lambda_j(h)\ :=\ \<L_V\varphi_j,\varphi_j\>\,.$$ 
From Proposition~\ref{pr.quad-form}, one knows that for all $j\in\{2,\ldots,n_0\}$, one has 
\be\label{eq:asympt-tildelambdaj}
\tilde\lambda_j(h)\ =\  \sum_{\s \in {\bf j}(\m_{j})}\frac{|\mu(\s)|}{2\pi}   \,    \frac{ D_{\m_j}}{D_{\s} }  \,  e^{-\frac{S_j}h}\big(1+\ooo(h)\big).
\ee
Moreover, since $(S_j)_{j\in\{1,\dots,n_{0}\}}$ is non-increasing, we deduce from this estimate that there exists $h_0>0$ and $C>0$ such that for all $h\in(0,h_0]$
and all $i,j\in\{1,\ldots,n_0\}$, one has
\be\label{eq:inclambdaj}
i\leq j\ \Longrightarrow\  \lambda_i(h)\leq C \lambda_j(h).
\ee
\ 

\noindent
The two following lemmata are straightforward consequence of the previous analysis.

\begin{lemma}\label{lem:orthogphik}
For every $j,k \in\{1,\ldots,n_0\}$ and $h\in(0,1]$, one has 
$$
\<L_V \varphi_j,\varphi_k\>\ =\ \delta_{jk}\,\tilde\lambda_j(h)\,.
$$
\end{lemma}

\noindent
\bp
When $j=k$, the statement if obvious. When $j\neq k$, then it follows from Lemma \ref{lem:orthogQM} that we are in one of the three following cases:

\begin{itemize}
\item[--] either $\supp(\varphi_j)\cap\supp(\varphi_k)=\emptyset$ and the conclusion is obvious,
\item[--] either there exists $c_{h}>0$ such that  $\varphi_j=c_{h}$ on $\supp(\varphi_k)$ and
$$
\<L_V \varphi_j,\varphi_k\>\ =\ \<L_V (c_{h}),\varphi_k\>\ =\ 0\,,
$$
\item[--] 
 or there exists $c_{h}>0$ such that  $\varphi_k=c_{h}$ on $\supp(\varphi_j)$ and
$$
\<L_V \varphi_j,\varphi_k\>\ =\ \< \varphi_j,L_V^*\varphi_k\>\ =\ \< \varphi_j,L_V^*(c_{h})\>\ =\ 0\,.
$$
\end{itemize}
\ep

 \begin{lemma}\label{prop:quad1formadj} 
For $\rho_{0},\delta_{0}$ sufficiently small and every $j\in\{1,\ldots,n_0\}$, it holds in the limit $h\to 0$,
\begin{equation}\label{quasimodal2bis}
 \Vert L_{V}\, \varphi_j\Vert =
 \mathcal O(\sqrt{h\tilde\lambda_j(h)})     \,.   
\end{equation}
and 
\begin{equation}\label{quasimodal2bis-star}
 \Vert L_{V}^*\, \varphi_j\Vert =
 \mathcal O(\sqrt{\tilde\lambda_j(h)})     \,.   
\end{equation}

\end{lemma}

\noindent
\bp
This is a simple rewriting of Proposition \ref{pr.quad-1-form}, using the fact that for every $\m\in \mathcal U^{(0)}$
and $h\in(0,1]$, $\varphi_{\m,h}=\frac{\psi_{\m,h}}{\Vert\psi_{\m,h}\Vert}$.
\ep

\noindent
We now introduce, for every $h>0$ small enough, the spectral projector $\Pi_{h}$ associated with the $n_{0}$ smallest eigenvalues of $L_V$ as described in Theorem~\ref{th:small-ev-Pphi}. Let then $\epsilon_0$ be given by Theorem~\ref{th:small-ev-Pphi}. According to Theorem~\ref{th:small-ev-Pphi},  
 for every $h>0$ small enough, $\Pi_{h}$ satisfies 
\be
\label{eq.def-Pi-h}
\Pi_h\ :\ =\frac 1{2i\pi}\int_{z\in \partial D(0,\frac {\epsilon_0}2)}(z-L_V)^{-1}dz
\ee
and in particular:
\be
\label{eq.Pi-h}
\Pi_h\ =\ \ooo(1)\,.
\ee
\begin{lemma}\label{lem:estprojQM}
For all $j\in\{1,\ldots,n_0\}$, we have, in the limit $h\to 0$, 
\be\label{eq:estprojQM1}
\Vert (1-\Pi_h)\varphi_j\Vert=\ooo(\sqrt{h\tilde\lambda_j(h)})
\ee
and
\be\label{eq:estprojQM1star}
\Vert (1-\Pi_h^*)\varphi_j\Vert=\ooo(\sqrt{\tilde\lambda_j(h)})
\ee
\end{lemma}

\noindent
\bp
Thanks to the resolvent identity, one has
\begin{equation*}
\begin{split}
(1-\Pi_h)\varphi_j&=\frac 1{2i\pi}\int_{z\in \partial D(0,\frac {\epsilon_0}2)}(z^{-1}-(z-L_V)^{-1})\varphi_j\,dz\\
&=\frac {-1}{2i\pi}\int_{z\in \partial D(0,\frac {\epsilon_0}2)}z^{-1}(z-L_V)^{-1}L_V\varphi_j\,dz.
\end{split}
\end{equation*}
Moreover, it follows from Theorem \ref{th:small-ev-Pphi} and from \eqref{eq.unit} that for any $z\in \partial D(0,\frac {\epsilon_0}2)$,
$$
 \Vert (z-L_V)^{-1}\Vert_{L^2(m_h)\rightarrow L^2(m_h)}\ =\ \ooo(1).
$$
Combined with \eqref{quasimodal2bis},  this proves \eqref{eq:estprojQM1}. On the other hand, one has similarly
$$
 (1-\Pi_h^*)\varphi_j\ =\ \frac {-1}{2i\pi}\int_{z\in \partial D(0,\frac {\epsilon_0}2)}z^{-1}(z-L_V^*)^{-1}L_V^*\varphi_j\,dz
$$
and $\Vert (z-L_V^*)^{-1}\Vert_{L^2(m_h)\rightarrow L^2(m_h)}=\ooo(1)$. Then, \eqref{eq:estprojQM1star} follows immediately from 
\eqref{quasimodal2bis-star}.
\ep

\begin{proposition}\label{lem:formquadproj}
For every $j\in\{1,\dots,n_{0}\}$ and $h>0$ small enough, let us define
$v_j:=\Pi_h\varphi_j$. Then, there exists $c>0$ such that for all $j,k\in\{1,\ldots,n_0\}$, one has
in the limit $h\to 0$,
\be\label{eq:orthogproj}
\<v_j,v_k\>\ =\ \delta_{jk}+\ooo(e^{-\frac ch})
\ee
and
\be\label{eq:formquadproj}
\<L_Vv_j,v_k\>\ =\ \delta_{jk}\tilde \lambda_j(h)+\ooo(\sqrt{h\tilde\lambda_j(h)\tilde\lambda_k(h)})\,.
\ee
In particular, it follows from \eqref{eq:orthogproj} that for every $h>0$ small enough, 
the family $(v_{1},\dots,v_{n_{0}})$ is a basis of $\Ran\, \Pi_{h}$.
\end{proposition}

\noindent
\bp Since, for some $c>0$, every $j\in\{1,\dots,n_{0}\}$, and every $h>0$ small enough, it holds
 $\tilde\lambda_j(h)=\ooo(e^{-\frac c h})$, the first identity follows directly from \eqref{eq:orthogQM}, \eqref{eq:estprojQM1}, and from the relation
 $$
 \<v_j,v_k\>\ =\ 
 \<\varphi_j,\varphi_k\>
+ \<\varphi_j,v_k-\varphi_k\>
+
 \<v_j-\varphi_j,v_k\>\,.
 $$
To prove the second estimate, observe that 
\begin{equation*}
\begin{split}
\<L_Vv_j,v_k\>&=
\<L_V\Pi_h \varphi_j,\Pi_h \varphi_k\>\\
&=\<L_V \varphi_j,\varphi_k\>+\<L_V(\Pi_h-1) \varphi_j,\varphi_k\>+
\<\Pi_{h}L_V \varphi_j,(\Pi_h-1)\varphi_k\>\\
&=\<L_V \varphi_j,\varphi_k\>+\<(\Pi_h-1) \varphi_j,L_V^{*}\varphi_k\>+
\<\Pi_{h}L_V \varphi_j,(\Pi_h-1)\varphi_k\>\,.
\end{split}
\end{equation*}
 Moreover, thanks to
Lemma~\ref{prop:quad1formadj}, \eqref{eq.Pi-h}, and
 Lemma~\ref{lem:estprojQM}, one has
$$
\vert\<(\Pi_h-1) \varphi_j,L_V^*\varphi_k\>\vert\ \leq\  \Vert (\Pi_h-1) \varphi_j\Vert\Vert L_V^*\varphi_k\Vert
\ =\ \ooo(\sqrt{h\tilde\lambda_j(h)\tilde\lambda_k(h)})
$$
 and
$$
\vert \<\Pi_{h}L_V \varphi_j,(\Pi_h-1)\varphi_k\>\vert\ \leq\  \|\Pi_{h}\| \,\Vert L_V\varphi_j\Vert \Vert (\Pi_h-1)\varphi_k\Vert\ =\ \ooo(\sqrt{h^2\tilde\lambda_j(h)\tilde\lambda_k(h)})
\,.
$$
Gathering these two estimates and using Lemma~\ref{lem:orthogphik}, we obtain \eqref{eq:formquadproj}.
\ep

\

\noindent
We now orthonormalize the basis $(v_1,\dots,v_{n_{0}})$ of $\Ran\,\Pi_{h}$  by a Gram-Schmidt procedure: for all $j\in \{1,\dots,n_{0}\}$, let us define by induction 
\be\label{eq:deftildeej}
\tilde e_j\ =\ v_j-\sum_{k=1}^{j-1}\frac {\<v_j,\tilde e_k\>}{\Vert\tilde e_k\Vert^2}\,\tilde e_k
\quad\text{and then}\quad e_j\ =\ \frac{\tilde e_j}{\Vert \tilde e_j\Vert}\,.
\ee
\begin{lemma}\label{lem:orthGS}
There exists $c>0$ such that for all $j\in\{1,\dots,n_{0}\}$, one has
in the limit $h\to 0$:
$$
\tilde e_j\ =\ v_j+\sum_{k=1}^{j-1}\alpha_{j,k} v_k
$$
with 
$\alpha_{jk}=\ooo(e^{-\frac ch})$. In particular, it holds: 
$$
\forall\,j\,\in\,\{1,\dots,n_{0}\}\,,\ \ \ \Vert\tilde e_j\Vert\ =\ 1+\ooo(e^{-\frac ch}).
$$
\end{lemma}

\noindent
\bp One proceeds by induction on $j$. For $j=1$, one has $\tilde e_1=v_1=\varphi_1=1$ and there is nothing to prove.
Suppose now that the above formula is true for all $\tilde  e_l$ with $1\leq l\leq j<n_{0}$. Then
$\tilde e_{j+1}=v_{j+1}-r_{j+1}$ with 
$$
r_{j+1}\ =\ \sum_{k=1}^{j}\frac {\<v_{j+1},\tilde e_k\>}{\Vert\tilde e_k\Vert^2}\,\tilde e_k\,.
$$
Since by induction, $\Vert\tilde e_k\Vert=1+\ooo(e^{-\frac ch})$ for all $k\in\{1,\dots,j\}$, it follows that
$$
r_{j+1}\ =\ (1+\ooo(e^{-\frac ch}))\sum_{k=1}^{j}\<v_{j+1},\tilde e_k\>\tilde e_k\,.
$$
Moreover, for all $k\in\{1,\dots, j\}$, one also has by induction
$$
\tilde e_k\ =\ v_k+\sum_{l=1}^{k-1}\alpha_{k,l}v_l\ =\ \sum_{l=1}^{k}\beta_{k,l} v_k
$$
with $\beta_{k,l}=\ooo(1)$ for any $l\in\{1,\dots,k\}$ (and actually $\beta_{k,l}=\ooo(e^{-\frac ch})$ when $l<k$), which implies
$$
r_{j+1}=(1+\ooo(e^{-\frac c h}))\sum_{k=1}^{j}\sum_{l,m=1}^{k}\beta_{k,l}\beta_{k,m}\<v_{j+1},v_l\>v_m.
$$
Since, thanks to Proposition~\ref{lem:formquadproj},  it holds $\<v_{j+1},v_l\>=\ooo(e^{-\frac c h})$ for all $l, m\leq k<j+1$, then
$$
r_{j+1}\ =\ \sum_{m=1}^{j}\gamma_{j,m}v_m\,,
$$
where  $\gamma_{j,m}=\ooo(e^{-\frac c h})$ for all $m\in\{1,\dots,j\}$. This proves the first part of the lemma. The second one is obvious.
\ep

\begin{proposition}\label{prop:matapp}
For all $j,k\in\{1,\ldots,n_0\}$, one has in the limit $h\to 0$:
$$
\<L_V e_j,e_k\>\ =\ \delta_{jk}\tilde\lambda_j(h)+\ooo(\sqrt{h\tilde\lambda_j(h)\tilde\lambda_k(h)})\,.
$$
\end{proposition}

\noindent
\bp
Thanks to Lemma~\ref{lem:orthGS}, one has for  all $j,k\in\{1,\ldots,n_0\}$,
$$
\<L_V\tilde e_j,\tilde e_k\>=\<L_Vv_j,v_k\>+\sum_{p=1}^{j-1}\sum_{q=1}^{k-1}\alpha'_{p,q}\<L_Vv_p,v_q\>\,,
$$
where, for all $p,q$, it holds  $\alpha'_{p,q}=\alpha_{j,p}\alpha_{k,q}=\ooo(e^{-\frac c h})$. Combined with Proposition~\ref{lem:formquadproj}, this implies
\be\label{eq:matapp1}
\<L_V\tilde e_j,\tilde e_k\>=\delta_{jk}\tilde\lambda_j(h)+\ooo(\sqrt{h\tilde\lambda_j(h)\tilde\lambda_k(h)})
+\sum_{p=1}^{j-1}\sum_{q=1}^{k-1}\alpha'_{p,q}\<L_Vv_p,v_q\>\,.
\ee
On the other hand, thanks to Proposition~\ref{lem:formquadproj} and 
\eqref{eq:inclambdaj}, one has
in the limit $h\to 0$,
for all $1 \leq p< j$ and $1\leq q <k$,
\begin{align*}
\<L_V v_p,v_q\>\ =\ \delta_{pq}\tilde\lambda_p(h)+\ooo(\sqrt{h\tilde\lambda_p(h)\tilde\lambda_q(h)})
&\ =\ \ooo(\sqrt{\tilde\lambda_p(h)\tilde\lambda_q(h)})\\
&\ =\ \ooo(\sqrt{\tilde\lambda_j(h)\tilde\lambda_k(h)})\,.
\end{align*}
Combined with \eqref{eq:matapp1} and using the fact that $\alpha'_{p,q}=\ooo(e^{-\frac c h})=\ooo(\sqrt h)$, this shows that 
$$
\<L_V\tilde e_j,\tilde e_k\>\ =\ \delta_{jk}\tilde\lambda_j(h)+\ooo(\sqrt{h\tilde\lambda_j(h)\tilde\lambda_k(h)}).
$$
Eventually, since $e_k=(1+\ooo(e^{-\frac c h}))\tilde e_k$ according to Lemma~\ref{lem:orthGS}, we obtain
$$
\<L_Ve_j,e_k\>=(1+\ooo(e^{-\frac c h}))\<L_V\tilde e_j,\tilde e_k\>=\delta_{jk}\tilde\lambda_j(h)+\ooo(\sqrt{h\tilde\lambda_j(h)\tilde\lambda_k(h)})\,,
$$
which completes the proof.\ep
\

\noindent
We are now in position to prove Theorem \ref{th.main}.
We recall that $(e_1,\ldots,e_{n_0})$ is an orthonormal basis of $\Ran\,\Pi_{h}$ and 
that $L_V|_{\Ran\,\Pi_{h}}:\Ran\,\Pi_{h}\rightarrow \Ran\,\Pi_{h}$ has exactly 
$n_0$ eigenvalues $\lambda_1,\dots,\lambda_{n_0}$, with $\lambda_{j}=0$ iff $j=1$,  counted with algebraic multiplicity. Let us denote $\hat e_j=e_{n_0+1-j}$ and let $\mc M$ denote the matrix of $L_V$ in the basis
$(\hat e_1,\ldots,\hat e_{n_0})$. Since this  basis is orthonormal, it holds
$$\mc M\ =\ \big(\<L_V\hat e_k,\hat e_j\>\big)_{j,k\in\{1,\dots,n_{0}\}}\,.$$ 
Moreover, since 
$$L_V(\hat e_{n_0})=L_V(e_1)=0\ \text{ and }\ L_V^*(\hat e_{n_0})=0\,,$$
then $\mathcal M$ has the form 
$$\mathcal M\ =\ \begin{pmatrix} \mathcal M' & 0 \\ 0 &0 \end{pmatrix} \quad
\text{with}\quad \mathcal M'\ :=\  \big (\<L_V\hat e_k,\hat e_j\> \big)_{j,k\in\{1,\ldots n_0-1\}}\,. $$
On the other hand, denoting
$\hat\lambda_j(h):=\tilde\lambda_{n_0+1-j}(h)$ for $j\in\{1,\dots,n_{0}-1\}$, one deduces from Proposition~\ref{prop:matapp} that 
for every $j,k\in\{1,\ldots n_0-1\}$, it holds in the limit $h\to0$,
\begin{equation*}
\begin{split}
\<L_V\hat e_k,\hat e_j\>\ =\ \<L_Ve_{n_0-k},e_{n_0-j}\>
&\ =\ \delta_{jk}\hat\lambda_j(h)+\ooo(\sqrt{h\hat \lambda_{j}(h)\hat\lambda_{k}(h)})\,,
\end{split}
\end{equation*}
that is
\be\label{eq:coeffmat}
\<L_V\hat e_k,\hat e_j\>\ =\ \sqrt{\hat \lambda_{j}(h)\hat\lambda_{k}(h)}\big(\delta_{jk}+\ooo(\sqrt h)\big).
\ee
For all $j\in\{1,\ldots,n_0-1\}$, let us now define 
\begin{align*}
\hat S_j\ :=\ S_{n_0+1-j} \quad \text{and}\quad \nu_j\ :=\ \zeta(\m_{n_{0}+1-j})\ &=\ 
\sum_{\s \in {\bf j}(\m_{n_{0}+1-j})}\frac{|\mu(\s)|}{2\pi}   \,    \frac{ D_{\m_{n_{0}+1-j}}}{D_{\s} }\\
\ &=\ 
e^{\frac{\hat S_j}{h}}\hat\lambda_j(h)\big(1+\ooo(h)\big)\,,
\end{align*}
where 
$\zeta(\m)$, $\m\in\uuu^{(0)}\setminus\{\underline\m\}$, is defined in \eqref{eq.pref},
and
the last estimate follows from \eqref{eq:asympt-tildelambdaj}.
Since the sequence $(S_j)_{j\in\{2,\dots,n_{0}\}}$ is non-increasing, there exists a partition $J_1\sqcup\ldots\sqcup J_p$ of $\{1,\ldots,n_0-1\}$ such that for all $k\in\{1,\ldots,p\}$, there exists $\iota(k)\in\{1,\ldots,n_0-1\}$ such that 
\be\label{eq:inciota}
\forall j\in J_k\,,\ \hat S_j=\hat S_{\iota(k)}\quad\text{and}\quad \forall 1\leq k< k'\leq p\,,\  \hat S_{\iota(k)}<\hat S_{\iota(k')}.
\ee
Hence, we deduce from \eqref{eq:coeffmat} that
$$
\mc M'\ =\ \widehat\Omega \,\big(\jjj+\ooo(\sqrt h)\big)\,\widehat\Omega
$$
with 
$$\jjj\ =\ \diag(\nu_j,\;j=1,\ldots,n_0-1)$$ and 
$$\widehat\Omega
\ =\ \diag(e^{-\frac{\hat S_{j}}{2h}},\; j=1,\ldots,n_0-1)
\ =\ \diag(e^{-\frac{\hat S_{\iota(k)}}{2h}}I_{r_k},\;k=1,\ldots,p)\,,$$ where, for every $k\in\{1,\dots,p\}$, 
$r_k={\rm card} (J_k)$. Factorizing by 
$e^{-\frac{\hat S_{\iota(1)}}h}$, we get
$$
\mc M'\ =\ e^{-\frac {\hat S_{\iota(1)}}h} \Omega \,\big(\jjj+\ooo(\sqrt h)\big)\,\Omega
$$
with 
$$\Omega\ =\ \diag\big(e^{\frac{\hat S_{\iota(1)}-\hat S_{\iota(k)}}{2h}}\,I_{r_k},k=1,\ldots,p\big)\,.$$ Denoting 
$\tau_1=1$ and, for $k\in\{2,\dots,p\}$, $\tau_k=e^{\frac{\hat S_{\iota(k-1)}-\hat S_{\iota(k)}}{2h}}$, we observe that, thanks to 
\eqref{eq:inciota}, $\tau_k$ is exponentially small when $h\rightarrow 0$. Moreover, with this notation, one has
$$
\Omega\ =\ \diag\big(\tau_1 I_{r_1},\tau_1\tau_2 I_{r_2},\ldots,(\Pi_{j=1}^p\tau_j )I_{r_p}\big).
$$
This shows that  $e^{-\frac {\hat S_{\iota(1)}}h} \mc M'$ is a graded matrix in the sense of Definition \ref{defin:gradedmat}. Hence, we can apply Theorem \ref{th:specgradedmat}
and we get that in the limit $h\to 0$,
$$
\sigma(\mc M')\ \subset\ \bigsqcup_{k=1}^p e^{-\frac{\hat S_{\iota(1)}}h}\varepsilon_k^2\big(\sigma(M_k)+\ooo(\sqrt h)\big)\,,
$$
where for every $k\in\{1,\dots,p\}$, $\varepsilon_k=\prod_{l=1}^k\tau_l$ and 
$M_k=\diag(\nu_j,\,j\in J_k)$. 
Moreover, still according to Theorem~\ref{th:specgradedmat},
$\mc M'$ admits in the limit $h\to 0$, for every $k\in\{1,\dots,p\}$ and every eigenvalue
$\lambda$ of $M_{k}$ with multiplicity $r'_{k}$, exactly
$r'_{k}$ eigenvalues counted with multiplicity of order $e^{-\frac{\hat S_{\iota(1)}}h}\varepsilon_k^2\big(\lambda+\ooo(\sqrt h)\big)$.\medskip

\noindent
Going back to the initial parameters, one has, for every $k\in\{1,\dots,p\}$, 
$$e^{-\frac{\hat S_{\iota(1)}}h}\varepsilon_k^2=e^{-\frac{\hat S_{\iota(k)}}h}
\quad\text{and}\quad
\sigma(M_k)=\{\nu_j,\,j\in J_{k}\}\,.
$$  Hence, the eigenvalues of $\mc M'$ satisfy:
$$
\forall\,j\,\in\,\{1,\dots,n_{0}-1\}\,,\ \ \ \lambda_{n_0+1-j}(h)\ =\ e^{-\frac{\hat S_j}{h}}\big(\nu_j+\ooo(\sqrt h)\big)\,,
$$
which is exactly the announced result.

\subsection{Proof of Theorem~\ref{co.main}}

As in the preceding subsection, 
 we denote for shortness
$$\lp\cdot,\cdot \rp=\lp\cdot,\cdot \rp_{L^{2}(m_{h})}\,,\ \ \ 
\|\cdot \|=\|\cdot \|_{L^{2}(m_{h})}\,,\ \ \ 
L_{V,b,\nu}=L_V\,,$$
and we label the local minima $\m_1,\ldots,\m_{n_0}$ of $V$ so that
$(S(\m_{j}))_{j\in\{1,\dots, n_{0}\}}$ is non-increasing (see \eqref{eq.S}):
$$S(\m_1)=+\infty \ \ \text{and, for all $j \in\{2,\dots,n_{0}\}$,}\ \  S(\m_{j+1})\leq S(\m_j)<+\infty\,.$$
Let moreover
$\m^{*}\in \uuu^{(0)}\setminus\{\underline\m\}$ be such that
\begin{equation}
\label{eq.m*'}
S(\m^{*})\ =\ \ S(\m_{2})
\quad\text{and}\quad  \zeta(\m^{*})\ =\ \min_{\m\in S^{-1}(S(\m_{2}))}\zeta(\m)\,,
\end{equation}
where the prefactors $\zeta(\m)$, $\m \in  \uuu^{(0)}\setminus\{\underline\m\}$, are defined
 in \eqref{eq.pref}, and let us  define, 
 for any $h>0$,
 $$
 \lambda(h)\ :=\ \zeta(\m^{*})\,e^{-\frac{S(\m^{*})}h}\,.
 $$
According to the unitary equivalence  (see \eqref{eq.unit})
$$L_{V}\ =\ \frac1h \mathsf{U} \,P_{\phi} \,\mathsf{U}^{*}\,,$$
and
 to the localization of the spectrum of $P_{\phi}$ stated in Proposition~\ref{prop:elemPphi}  
and in Theorem~\ref{th:small-ev-Pphi}, 
it holds  for every $h>0$ small enough, taking $\epsilon_{0}$ as in the statement of Theorem~\ref{th:small-ev-Pphi},
\begin{equation}
\label{eq.sep}
\|\,e^{-tL_V}\,-\Pi_{0}\,\|\ \leq\ 
\|\,e^{-tL_V}\,\Pi_{h}-\Pi_{0}\,\|+
\|\,e^{-tL_V}\,(\Id-\Pi_{h})\|\,,
\end{equation}
where, as in the preceding subsection, 
$$
\Pi_{h}\ : =\ \frac 1{2i\pi}\int_{z\in \partial D(0,\frac{\epsilon_{0}}{2})}(z-L_V)^{-1}dz.
$$
Moreover, it follows from Proposition \ref{prop:elemPphi} that $\sigma(P_\phi)\subset \Gamma_{\Lambda_0}\subset \tilde \Gamma_{\Lambda_0}$ with 
$\tilde \Gamma_{\Lambda_0}=\{z\in\C,\;\vert \Im(z)\vert\leq \Lambda_0(\Re(z)+1)\}$. Hence, 
 for every  $t>0$,
 the operator $e^{-tL_V}(I-\Pi_{h})$
can be written as the complex integral
$$
e^{-tL_V}(\Id-\Pi_{h})\ =\ -\int_{\Gamma_{0}\cup\Gamma_{\pm}} e^{-tz}(z-L_V)^{-1}dz\,,
$$
where$$
\Gamma_{0}\ =\ \left\{\frac{\epsilon_{0}}{2}+i \Lambda_{0} x\,,\ x\in[-\frac{\epsilon_{0}}{2}-\frac1h,\frac{\epsilon_{0}}{2}+\frac1h]\right\}
$$
and
$$
\Gamma_{\pm}\ =\ \left\{x\pm i \Lambda_{0}(x+\frac1h)\,,\ x\in[\frac{\epsilon_{0}}{2},+\infty)\right\}
\,.
$$
From the resolvent estimates proven in
 Theorem~\ref{th:small-ev-Pphi}, it holds
$(z-L_{V})^{-1}=\ooo(1)$
uniformly on $\Gamma_{0}$, and then, for every  $t>0$,  
$$
\int_{\Gamma_{0}} e^{-tz}(z-L_{V})^{-1}dz\ =\ e^{-t\frac{\epsilon_{0}}{2}}\ooo(\frac 1h)\,.
$$
Using 
in addition the resolvent estimates proven in  Proposition~\ref{prop:elemPphi},
it holds $\|(z-L_{V})^{-1}\|\leq \frac 1{\Re z}\leq \frac2{\epsilon_{0}}$ on $\Gamma_{\pm}$, and then 
$$
\int_{\Gamma_{\pm}} e^{-tz}(z-L_{V})^{-1}dz\ =\ \ooo(1)\int_{\frac{\epsilon_{0}}{2}}^{+\infty}
e^{-tx}dx
\ =\ \frac{e^{-t\frac{\epsilon_{0}}{2}}}t\ooo(1)\,.
$$
It follows that for every $t>0$, it holds
$$
\|\,e^{-tL_{V}}\,(\Id-\Pi_{h})\|\ =\ e^{-t\frac{\epsilon_{0}}{2}}\,\ooo\big(\frac 1t + \frac 1h\big)\,. 
$$
Moreover,  $e^{-tL_{V}}\,(\Id-\Pi_{h})=\ooo(1)$
since $\Pi_{h}=\ooo(1)$ (see \eqref{eq.Pi-h}) and $e^{-tL_{V}}=\ooo(1)$ (by maximal accretivity of $L_V$).
Hence, there exists $C>0$ such that for every $t\geq 0$ and $h>0$ small enough, it holds
$$\|\,e^{-tL_{V}}\,(I-\Pi_{h})\|\ \leq\ C\min\{1, \frac{e^{-t\frac{\epsilon_{0}}{2}}}{h}\}
\ \leq\ 2Ce^{-\lambda(h)t}\,.
$$
Thus, according to \eqref{eq.sep}, it just remains to show
that 
\begin{equation}
\label{eq.exp-L}
\exists \,C\,>\,0\,,\ \|\,e^{-tL_{V}}\,\Pi_{h}-\Pi_{0}\,\|\ \leq\ C\,e^{-(\lambda(h)-C\sqrt h)t}\,.
\end{equation}

\noindent
To this end, let us first recall from 
Proposition~\ref{prop:elemPphi} that 
the spectral projector $\Pi_{\{0\}}$ associated with
the eigenvalue
$0$ of
 $L_{V}$  has rank $1$ and is actually the orthogonal projector
$\Pi_{0}$ on $\sspan\{1\}$ according to the relations
$$\sspan\{1\}\ =\ \Im \Pi_{\{0\}}\ =\ \Im \Pi^{*}_{\{0\}}\ =\ (\Ker \Pi_{\{0\}})^{\perp}\,.$$
It follows that
$$
e^{-tL_{V}}\,\Pi_{h}-\Pi_{0}\ =\  e^{-tL_{V}}\,\big(\Pi_{h}-\Pi_{\{0\}}\big)\,.
$$
Since moreover $\Pi_{h}-\Pi_{\{0\}}=\ooo(1)$ (thanks to the resolvent estimate of Theorem \ref{th:small-ev-Pphi}), 
it suffices to show that
$$
\exists \,C\,>\,0\,,\ \|\,e^{-tL_{V}}\,\big(\Pi_{h}-\Pi_{\{0\}}\big)|_{\Ran (\Pi_{h}-\Pi_{\{0\}})}\,\|\ \leq\ C\,e^{-(\lambda(h)-C\sqrt h)t}\,.
$$

\noindent
Using  the notation of the preceding subsection,
this means proving that
the matrix $\mc M'$ of $L_{V}$ in the orthonormal basis
$(\hat e_1,\ldots,\hat e_{n_0-1})$ of $\Ran(\Pi_{h}-\Pi_{0})$
satisfies
$$
\exists \,C\,>\,0\,,\ \|\,e^{-t\mc M'}\,\|\ \leq\ C\,e^{-(\lambda(h)-C\sqrt h)t}\,.
$$
Let us now consider a subset 
$ \mathcal V^{(0)}$ (in general non unique) of $\uuu^{(0)}\setminus\{0\}$
such that 
$$\m\in \mathcal V^{(0)}\mapsto (\zeta(\m),S(\m))\in \{(\zeta(\m),S(\m)), \m\in \uuu^{(0)}\setminus\{0\} \}
\ \ \ \text{is a bijection.}$$
Then,
for any $K>0$ and for every $h>0$ small enough,
the closed disks of the complex plane
$$
D_{\m,K}\ :=\ D\big(\zeta(\m)e^{-\frac {S(\m)}h},
K\sqrt h e^{-\frac {S(\m)}h}\big)\,,\ \ \m\in\vvv^{(0)}\,,
$$
are included in $\{\Re z>0\}$ and two by two disjoint.
Moreover, according to 
Theorem~\ref{th.main}, $K>0$ can be chosen large enough 
so that when $h>0$ is small enough,
the $n_{0}-1$  non zero small eigenvalues 
of $L_{V}$ are included in 
$$
\cup_{\m\in\vvv^{(0)}} D\big(\zeta(\m)e^{-\frac {S(\m)}h},
\frac K2\sqrt h e^{-\frac {S(\m)}h}\big)\,.
$$
In particular, for every $t\geq 0$ and for every $h>0$ small enough,
 it holds
\begin{align}
\nonumber
e^{-t\mc M'}&\  =\  \sum_{\m\in \mathcal V^{(0)}} \frac 1{2i\pi}\int_{z\in \partial D_{\m,K}}e^{-tz}(z-\mc M')^{-1}dz.
\end{align}
Using now the specific form of $\mathcal M'$ exhibited in the preceding section
and Theorem~\ref{th:specgradedmat},
it holds  for every $\m\in \vvv^{(0)}$, in the limit $h\to 0$,
\begin{align*}
\frac 1{2i\pi}\int_{z\in \partial D_{\m,K}}e^{-tz}(z-\mc M')^{-1}dz
 \ =\ \ooo\Big(
e^{-t \zeta(\m)e^{-\frac {S(\m)}h}(1-K\sqrt h)}\Big)\,.
\end{align*}
Indeed, the  resolvent estimate of Theorem~\ref{th:specgradedmat} implies
\begin{align}
\label{eq.exp-L'}
\forall z\in\pa D_{\m,K}\,,\  \Vert (\mc M'-z)^{-1}\Vert  =  \ooo\big(\dist(z,\sigma(\mc M'))^{-1}\big)
 =  \ooo(\frac{1}{\sqrt h}e^{\frac{S(\m)}{h}}).
\end{align}
The relation \eqref{eq.exp-L} follows easily, which concludes the first part of Theorem~\ref{co.main}.\medskip

\noindent
Finally, let us assume that the element $\m^{*}$
satisfying \eqref{eq.m*'} is unique.
In this case, $\m^{*}$ necessarily belongs to $\vvv^{(0)}$ and the associated eigenvalue
$\lambda(\m^{*},h)$  (see \eqref{eq.vp}) is then real and simple for every $h>0$ small enough.
In particular, it holds
$$
\frac 1{2i\pi}\int_{z\in \partial D_{\m^{*},K}}e^{-tz}(z-\mc M')^{-1}dz
\ =\ e^{-t\lambda(\m^{*},h)}\Pi_{\{\lambda(\m^{*},h)\}},
$$
where $\Pi_{\{\lambda(\m^{*},h)\}}$ is the spectral projector (whose rank is one)
$$
\Pi_{\{\lambda(\m^{*},h)\}}\ =\ 
\frac 1{2i\pi}\int_{z\in \partial D_{\m^{*},K}}(z-\mc M')^{-1}dz.
$$
Moreover, the resolvent estimate
\eqref{eq.exp-L'} 
 shows that
$\Pi_{\{\lambda(\m^{*},h)\}}=\ooo(1)$.
Since in addition, it holds in this case 
 (see \eqref{eq.vp})
\begin{align*}
\forall \m\in\vvv^{(0)}\setminus\{\m^{*}\}\,,\ \ \ \lambda(\m^{*},h)
&\ =\ 
\zeta(\m^{*})e^{-\frac {S(\m^{*})}h}(1+\ooo(\sqrt h))\\
&\ \geq\  
\zeta(\m)e^{-\frac {S(\m)}h}(1-K\sqrt h)
\end{align*}
for every $K>0$ and for every $h>0$ small enough,
 we obtain that in the limit $h\to0$,
 $$
 e^{-t\mc M'}\ =\  \ooo(e^{-t\lambda(\m^{*},h)})\,,
 $$
 and thus the relation
 \eqref{eq.exp-L} remains valid if ones replaces
 $\lambda(h)-C\sqrt h$ there by $\lambda(\m^{*},h)$.
This concludes the proof of Theorem~\ref{co.main}.

\appendix
\section{Some results in linear algebra}

The aim of this appendix is to give some handy tools of linear algebra adapted to the setting of non-reversible metastable problems 
considered in this paper. 
Let us start with some notations.\medskip

Given any matrix $M\in\mmm_d(\C)$ and $\lambda\in\sigma(M)$, we denote by $m(\lambda)$ the multiplicity of $\lambda$,
$m(\lambda)=\dim\Ker (M-\lambda)^d$. We recall that 
for every $r>0$ small enough,
\be
\label{eq:defdimproj}
m(\lambda)\ =\ \rank\big(\Pi_{D(\lambda,r)}(M)\big)\ =:\ n(D(\lambda,r);M)\,,
\ee
 where
$$
\Pi_{D(\lambda,r)}(M)\ =\ \frac 1{2i\pi}\int_{\partial D(\lambda,r)}(z-M)^{-1}dz\,.
$$

\noindent
We denote by $\Dr_0(E)$ the set of complex matrices on a vector space $E$ which are diagonalizable and invertible.  

\noindent
Given two subsets $A,B\subset \C$, we say that $A\subset B+\ooo(h)$ if there exists $C>0$ such that $A\subset B+B(0,Ch)$.
\begin{definition}\label{defin:gradedmat}
Let $\Er=(E_j)_{j=1,\ldots,p}$ be a sequence of finite dimensional vector spaces $E_j$ of dimension $r_j>0$, let 
$E=\oplus_{j=1,\ldots, p}E_j$  and let $\tau=(\tau_2,\ldots,\tau_p)\in (\R_+^*)^{p-1}$. Suppose that  $(h,\tau) \mapsto\mmm_h(\tau)$ is a  map from 
$(0,1]\times(\R_+^*)^{p-1}$ into the set of complex matrices on $E$.

We say that $\mmm_h(\tau)$ is an $(\Er,\tau,h)$-graded matrix if there exists  
$\mmm'\in \Dr_0(E)$ independent of $(h,\tau)$ such that $\mmm_h(\tau)=\Omega(\tau) (\mmm'+\ooo(h))\Omega(\tau)$ 
with $\Omega(\tau)$ and $\mmm'$ such that 
\begin{itemize}
\item[--] $\mmm'=\diag(M_j,\,j=1,\ldots,p)$ with $M_j\in \Dr_0(E_j)$,
\item[--] $\Omega(\tau)=\diag(\varepsilon_j(\tau)I_{r_j},\,j=1,\ldots,p)$ with
$\varepsilon_1(\tau)=1$ and $\varepsilon_j(\tau)=(\prod_{k=2}^j\tau_k) $ for all $j\geq 2$.
\end{itemize}
Throughout, we denote by $\Gr(\Er,\tau,h)$ the set of $(\Er,\tau,h)$-graded matrices.
\end{definition}
\begin{lemma}\label{lem:Malphataugrad} Suppose that $\mmm_h(\tau)$ is a   family of $(\Er,\tau,h)$-graded matrices and that $p\geq 2$. Then,
one has 
\be\label{eq:taugradmatrix}
\mmm_h(\tau)=
\left(
\begin{array}{cc}
J(h)&\tau_2B_h^+(\tau')^{*}\\
\tau_2B^-_h(\tau')&\tau_2^2\nnn_h(\tau')
\end{array}
\right)
\ee
where
\begin{itemize}
\item[--] $J(h)=M_1+\ooo(h)$ with $M_1\in\Dr_0(E_1)$,
\item[--] $\nnn_h(\tau')\in\Gr(\Er',\tau',h)$ with $\tau'=(\tau_3,\ldots,\tau_p)$ and $\Er'=(E_j)_{j=2,\ldots,p}$,
\item[--] $B^\pm_h(\tau')\in \Mr(E_1,\oplus_{j=2}^p E_j)$ satisfies
$$
B^\pm_h(\tau')^*=(b^\pm_{2}(h)^*,\tau_3b^\pm_{3}(h)^*,\tau_3\tau_4b^\pm_{4}(h)^*,\ldots,\tau_3\ldots\tau_pb^\pm_{p}(h)^*)
$$
with $b^\pm_{j}(h):E_1\rightarrow E_j$ independent of $\tau$ and $b_j(h)=\ooo(h)$.
\end{itemize}
Moreover, the matrix $\nnn_h(\tau')-B^-_h(\tau')J(h)^{-1}B^+_h(\tau')^*$ belongs to $\Gr(\Er',\tau',h)$.
\end{lemma}

\noindent
\bp
Assume that $\mmm_h(\tau)=\Omega(\tau) (\mmm'+\ooo(h))\Omega(\tau)$ with $\Omega(\tau)$  and $\mmm'$ as in Definition \ref{defin:gradedmat}.
 First observe that
$$
\Omega(\tau)=
\left(
\begin{array}{cc}
I_{r_1}&0\\
0&\tau_2\Omega'(\tau')
\end{array}
\right)
$$
with 
$$
\Omega'(\tau')=
\left(
\begin{array}{ccccc}
I_{r_{2}}&0&\hdots&\hdots&0\\
0&\tau_3I_{r_{3}}&0&\hdots&0\\
\vdots&0&\ddots&\ddots&\vdots\\
\vdots&\ddots&\ddots&\ddots&0\\
0&\hdots&\hdots&0&\tau_3\tau_4\ldots\tau_pI_{r_p}
\end{array}
\right).
$$
On the other hand, we can write 
$$
\mmm'+\ooo(h)=
\left(
\begin{array}{cc}
J(h)&B'^{+}(h)^*\\
B'^{-}(h)&\nnn'(h)
\end{array}
\right)
$$
where 
$J(h)=M_1+\ooo(h)$ with  $M_1\in \Dr_0(E_1) $, $B'^{\pm}(h)=\ooo(h)$, and $\nnn'(h)=\nnn'_0+\ooo(h)$ with 
$\nnn'_0=\diag(M_j,\,j=2,\ldots,p)$.
Therefore,
$$
\Omega(\tau)(\mmm'+\ooo(h))\Omega(\tau)
=\left(
\begin{array}{cc}
J(h)&\tau_2 B'^+(h)^*\Omega'(\tau')\\
\tau_2 \Omega'(\tau')B'^-(h)&\tau_2^2\Omega'(\tau')\nnn'(h)\Omega'(\tau')
\end{array}
\right)
$$
has exactly the  form \eqref{eq:taugradmatrix} with $B^\pm_h(\tau')= \Omega'(\tau')B'^\pm(h)$ and 
$\nnn_h(\tau')=\Omega'(\tau')\nnn'(h)\Omega'(\tau')$.
By construction, $\nnn_h(\tau')$ belongs to $\Gr(\Er',\tau',h)$ and $B^\pm_h(\tau')$ has the required form.
\ep

\begin{lemma}\label{lem:invertmatdiag}
Let $M$ be a complex diagonalizable matrix. Then there exists $C>0$ such that 
$$
\forall\lambda\notin\sigma(M),\;\Vert (M-\lambda)^{-1}\Vert\leq C\dist(\lambda,\sigma(M))^{-1}\,.
$$
\end{lemma}

\noindent
\bp
Let $P$ be an invertible matrix such that $PMP^{-1}=D$ is diagonal. Then 
$$
\Vert (M-\lambda)^{-1}\Vert=\Vert P(D-\lambda)^{-1}P^{-1}\Vert\leq C\Vert (D-\lambda)^{-1}\Vert=C\dist(\lambda,\sigma(M))^{-1}.
$$
\ep 

The following theorem gives precise informations on the spectrum of graded matrices as introduced above. The proof is based on standard arguments, namely on
the Schur complement method and complex analysis. The use of these two tools permits to work by induction and to decompose the base vector space in order to isolate eigenspaces corresponding to eigenvalues of the same order and to see the remainder of the matrix as a perturbation.
Similar arguments were used in \cite{Mi19} in a self-adjoint framework.
We believe that this result could be useful in other contexts where the computation of 
clouds of eigenvalues cannot be carried out by standard self-adjoint arguments.
\begin{theorem}\label{th:specgradedmat}
Suppose that $\mmm_h(\tau)$ is $(\Er,\tau,h)$-graded. Then, there exists $\tilde \tau_0, h_0>0$ such that for all $0<\tau_j\leq\tilde \tau_0$ and all $h\in(0,h_0]$, one has
$$
\sigma(\mmm_h(\tau))\subset \bigsqcup_{j=1}^p\varepsilon_j(\tau)^2(\sigma(M_j)+\ooo(h)).
$$
Moreover, for any eigenvalue $\lambda$ of $M_j$ with multiplicity $m_j(\lambda)$, there exists $K>0$ such that, denoting $D_j:=\{z\in\C,\;\vert z-\varepsilon_j(\tau)^2\lambda_j\vert <K\varepsilon_j(\tau)^2 h\}$,
one has 
\be\label{eq:dimprojspect}
n(D_j;\mmm_h(\tau))=m_j(\lambda),
\ee
where $n(D_j;\mmm_h(\tau))$ is defined by \eqref{eq:defdimproj}.
Moreover, there exists $C>0$ such that
$$\Vert (\mmm_h(\tau)-z)^{-1}\Vert \leq C\dist(z,\sigma(\mmm_h(\tau)))^{-1}$$
 for all $z\in\C\setminus\cup_{j=1}^p\cup_{\lambda\in \sigma(M_j)}
B( \varepsilon_j(\tau)^2\lambda, \varepsilon_j(\tau)^2Kh)$.

\end{theorem}

\noindent
\bp
We prove the theorem by induction on $p$. Throughout the proof the notation $\ooo(\cdot)$ is uniform with respect to the parameters $h$ and $\tau$.
For $p=1$, one has $\mmm_h(\tau)=M_1+\ooo(h)$ with $M_1\in\mmm_{r_1}(\R)$ independent of $h$, diagonalizable and invertible. 
Let us denote $\lambda_j^1$, $j=1,\ldots,n_1$ its eigenvalues and $m_j=m(\lambda_j^1)$ the corresponding multiplicities.
The function $z\mapsto (\mmm_h-z)^{-1}$ is meromorphic on $\C$ with poles in $\sigma(\mmm_h)$. Moreover, 
Lemma~\ref{lem:invertmatdiag} and the identity 
$$
\mmm_h-z=(M_1-z)(\Id+(M_1-z)^{-1}\ooo(h)),\;\forall z\notin\sigma(M_1)
$$
show that for any $C>0$ large enough,
$(\mmm_h-z)$ is invertible on $\C\setminus\cup_{j=1}^{n_1}D(\lambda_j^1,Ch)$ with 
$\|(M_1-z)^{-1}\|=\ooo(\frac{1}{Ch})$
and
\be\label{eq:resolvneum}
(\mmm_h-z)^{-1}=(\Id+(M_1-z)^{-1}\ooo(h))^{-1}(M_1-z)^{-1}.
\ee
Hence, for every $C>0$ large enough, the associated spectral projector writes
$$
\Pi_{D(\lambda_j^1,Ch)}(\mmm_h)=\frac 1{2i\pi}\int_{\partial D(\lambda_j^1,Ch)}(\Id+\mathcal O(\frac 1C))^{-1}(z-M_1)^{-1}dz.
$$
This implies that for $C>0$ large enough,
$$\rank(\Pi_{D(\lambda_j^1,Ch)}(\mmm_h))=\rank(\Pi_{D(\lambda_j^1,Ch)}(M_1))=m_j,$$
which is exactly \eqref{eq:dimprojspect}.
As a consequence
$$
\sum_{j=1}^{n_1}\rank(\Pi_{D(\lambda_j^1,Ch)}(\mmm_h))=\sum_{j=1}^{n_1}m_j=r_1
$$
is maximal and hence $\sigma(\mmm_h)\subset \cup_{j=1}^{n_1}D(\lambda_j^1,Ch)$.
Eventually, \eqref{eq:resolvneum} shows that for any $z\in\C\setminus  \cup_{j=1}^{n_1}D(\lambda_j^1,Ch)$, one has 
$$
\Vert (\mmm_h-z)^{-1}\Vert\leq C'\Vert (M_1-z)^{-1}\Vert
$$
for some constant $C'>0$. Using Lemma \ref{lem:invertmatdiag} we get
$$
\Vert (\mmm_h-z)^{-1}\Vert\leq C'\dist(z,\sigma(M_1))^{-1}\leq C''\dist(z,\sigma(\mmm_h))^{-1}
$$
for all $z\in\C\setminus  \cup_{j=1}^{n_1}D(\lambda_j^1,2Ch)$. This completes the initialization step.\medskip

\noindent
Suppose now that $p\geq 2$ and 
let $\mmm_h(\tau)\in\Gr(\Er,\tau,h)$. We have
$$
\mmm_h(\tau)=\left(
\begin{array}{cc}
J(h)&\tau_2B^+_h(\tau')^*\\
\tau_2B^-_h(\tau')&\tau_2^2\nnn_h(\tau')
\end{array}
\right)
$$
with $J(h), B^\pm_h(\tau')$ and $\nnn_h(\tau')$ as in Lemma \ref{lem:Malphataugrad}.
In order to lighten the notation, we will drop the variables $\tau, \tau'$ in the proof below.
For $\lambda\in\C$, let
\begin{equation}\label{eq:schur1}
\ppp(\lambda):=\mmm_h(\tau)-\lambda
=
\left(
\begin{array}{cc}
J(h)-\lambda &\tau_2 B_h^{+,*}\\
\tau_2 B_h^-&\tau_2^2\nnn_h-\lambda
\end{array}
\right).
\end{equation}  
This is  an holomorphic function, and since it is non trivial, its inverse is well defined 
excepted for a finite number of values of $\lambda$ which are exactly the  spectral values of $\mmm_h$.\medskip 

\noindent
We first study the part of the spectrum of $\mmm_h$  which is of largest modulus.
Let $\lambda_{n}^1$, $n=1,\ldots,n_1$, denote the  eigenvalues of the matrix $M_1$. 
Since $J(h)=M_1+\ooo(h)$ and  $M_1\in \Dr_0(E_1)$, then the initialization step shows that  there exists $C>0$ such that 
$\sigma(J(h))\subset\cup_{n=1}^{n_1}D(\lambda_n^1,Ch)$. 
Moreover, since $M_1$ is invertible, there exists $c_1,d_1>0$ and $h_0>0$ such that for all $n=1,\ldots,n_1$, one has
$\lambda_n^1\in K(c_1,d_1)$ where $K(c_1,d_1)=\{z\in\C,\;c_1\leq \vert z\vert\leq d_1\}$.
Let $n\in\{1,\ldots,n_1\}$ be fixed and consider 
$D_n=D_n(h)=\{z\in\C,\;\vert z-\lambda^1_{n}\vert\leq Mh\}$ for some $M>C>0$ and
$\tilde D_n=\{z\in\C,\;\vert z-\lambda^1_{n}\vert\leq 2Mh\}$.
Observe that for $h>0$ small enough, the disks $\tilde D_n$ are disjoint. By definition, one has $\nnn_h(\tau')=\ooo(1)$ and since
 $\vert \lambda\vert \geq c_1-\ooo(h)\geq  c_1/ 2$, this implies that for $\tau_2>0$ small enough with respect to $c_{1}$ and  $\lambda\in \tilde D_n$, the matrix  $\tau_2^2\nnn_h(\tau')-\lambda$ is  invertible, and 
$(\tau_2^2\nnn_h(\tau')-\lambda)^{-1}=\ooo(1)$. 
Moreover, it follows from the initialization step that  for $\lambda\in\tilde D_n\setminus D_n$, $J(h)-\lambda$ is invertible and 
$$\Vert(J(h)-\lambda)^{-1}\Vert=\ooo(\dist(\lambda,\sigma(J(h))^{-1})=\ooo(h^{-1}).$$
 Combined with the fact that $B^\pm_h=\ooo(h)$, this implies that for 
$h>0$ small enough and $\lambda\in\tilde D_n\setminus D_n$, $J(h)-\lambda-\tau_2^2 B_h^{+,*}(\tau_2^2\nnn_h-\lambda)^{-1}B^{-}_h$ is invertible with 
\be\label{eq:invertshur1}
\begin{split}
\big(J(h)&-\lambda-\tau_2^2 B_h^{+,*}(\tau_2^2\nnn_h-\lambda)^{-1}B_h^-\big)^{-1}\\
&=(J(h)-\lambda)^{-1}\Big(I-\tau_2^2 B_h^{+,*}(\tau_2^2\nnn_h-\lambda)^{-1}B_h^-(J(h)-\lambda)^{-1}\Big)^{-1}\\
&=(J(h)-\lambda)^{-1}(I+\ooo(h)).
\end{split}
\ee
Hence, the standard Schur complement procedure shows that for $\lambda\in\tilde D_n\setminus D_n$, 
$\ppp(\lambda)$ is invertible with  inverse $\eee(\lambda)$  given by
\begin{equation}\label{eq:schur2}
\eee(\lambda)=
\left(
\begin{array}{cc}
E(\lambda)&-\tau_2 E(\lambda)B_h^{+,*}(\tau_2^2\nnn_h-\lambda)^{-1}\\
-\tau_2(\tau_2^2\nnn_h-\lambda)^{-1}B_h^-E(\lambda)&E_0(\lambda)
\end{array}
\right)
\end{equation}
with 
\begin{equation*}
E(\lambda)=\Big(J(h)-\lambda-\tau_2^2 B_h^{+,*}(\tau_2^2\nnn_h-\lambda)^{-1}B_h^-\Big)^{-1}
\end{equation*}
and 
\begin{equation*}
E_0(\lambda)=(\tau_2^2\nnn_h-\lambda)^{-1}+\tau_2^2(\tau_2^2 \nnn_h-\lambda)^{-1}B_h^-E(\lambda)B_h^{+,*}(\tau_2^2\nnn_h-\lambda)^{-1}\,.
\end{equation*}
Let us now consider  the spectral projector $\Pi_{D_n}(\mmm_h)$. Then, 
$$\rank(\Pi_{D_n}(\mmm_h))\  \geq\  \rank(\tilde\Pi_n)\,,$$ 
where we defined
$$
\tilde\Pi_n=\left(\begin{array}{cc}
\Id&0\\
0&0
\end{array}\right)\Pi_{D_n}(\mmm_h)\left(\begin{array}{cc}
\Id&0\\
0&0
\end{array}\right).
$$
On the other hand, an elementary computation shows that  
$$
\tilde \Pi_n=\frac 1{2i\pi}\int_{\partial D_n}
\left(
\begin{array}{cc}
E(\lambda)&0\\
0&0
\end{array}
\right) d\lambda=\left(\begin{array}{cc}
E_n&0\\
0&0
\end{array}
\right) 
$$
with 
\begin{align*}
E_n&\ =\ \frac 1{2i\pi}\int_{\partial D_n}\Big(J(h)-\lambda-\tau_2^2 B_h^{+,*}(\tau_2^2\nnn_h-\lambda)^{-1}B_h^-\Big)^{-1}d\lambda
\\ 
&\ =\ 
\frac 1{2i\pi}\int_{\partial D_n}
(J(h)-\lambda)^{-1}(I+\ooo(h)) d\lambda
\,,
\end{align*}
where the last equality follows from \eqref{eq:invertshur1}.
 It follows that for $h>0$ small enough,
the 
rank
of $E_n$ is bounded from below by the  multiplicity $m(\lambda_n^1)$ of $\lambda_n^1$ and hence 
\be\label{eq:minorrank}
\rank(\Pi_{D_n}(\mmm_h))\geq m(\lambda_n^1)
\ee
for all $n=1,\ldots,n_1$.\medskip

\noindent
Let us now study the part of the spectrum of order smaller than  $\tau_2^2$. 
Thanks to the last part of Lemma \ref{lem:Malphataugrad}, the matrix $\zzz_h(\tau'):=\nnn_h-B_h^-J(h)^{-1}B_h^{+,*}$ is classical $(\Er',\tau')$-graded. 
Hence, it follows from the induction hypothesis that uniformly with respect to $h$, one has
\be\label{eq:locspec2}
\sigma(\zzz_h(\tau'))\subset\bigsqcup_{j=2}^p\tilde \varepsilon_j^{2}(\sigma(M_j)+\ooo(h))
\ee
with $\tilde \varepsilon_j=\tau_2^{-1}\varepsilon_j=\prod_{l=3}^{j}\tau_l$ for $j\geq 3$ and $\tilde \varepsilon_2=1$. 
One also knows  that for all $j=2,\ldots,p$ 	and all $\lambda\in\sigma(M_j)$, one has
$$
\rank\Pi_{D_j}(\zzz_h)=m_j(\lambda)
$$
where $D_j=D(\lambda\tilde\varepsilon_j^2,Kh\tilde\varepsilon_j^2)$ for some $K>0$.
Moreover,
one has for all $z\notin \cup_{j=2}^p\cup_{\lambda\in \sigma(M_j)} D(\lambda\tilde\varepsilon_j^2,Kh\tilde\varepsilon_j^2)$ the resolvent estimate
\be\label{eq:resolvest2}
(\zzz_h(\tau')-z)^{-1}=\ooo(\dist(z,\sigma(\zzz_h(\tau'))^{-1}).
\ee
For $j=2,\ldots,p$, let  $\lambda_{1}^j,\ldots, \lambda_{n_j}^j$ denote the eigenvalues of the  matrix 
$M_j\in\ddd_0$. As above,  there exists $c_j,d_j>0$ such that
  $\lambda_{n}^j\in  K(c_j,d_j)$ for all $n=1,\ldots,n_j$.
Suppose now that $j\in\{2,\ldots,p\}$ and $n\in\{1,\ldots,n_j\}$ are fixed and consider, for $M>K$, 
$$D'_{j,n}=\{z\in\C,\;\vert z-\varepsilon_j^2\lambda^j_{n}\vert\leq Mh\varepsilon_j^2\}
= \tau_2^{-2}\{z'\in\C,\;\vert z'-\tilde\varepsilon_j^2\lambda^j_{n}\vert\leq Mh\tilde\varepsilon_j^2\}\,.
$$
Since $M_1$ is invertible, $J(h)-\lambda$ is invertible
and $(J(h)-\lambda)^{-1}=\mathcal O(1)$
 for $\lambda$ in   $D'_{j,n}$ and $h,\tau_{2}$ small enough.
 Moreover, for any
 $\lambda\in \pa D'_{j,n}$, it holds, noting $\lambda'=\tau_2^{-2}\lambda$,
\begin{align*}
\tau_2^2 \nnn_h-\lambda-\tau_2^2 B_h^-(J(h)-\lambda)^{-1}B_h^{+,*}
&=\tau_2^2 (\nnn_h-\lambda'-B_h^-(J(h)-\lambda)^{-1}B_h^{+,*})\\
&=\tau_2^2 (\zzz_h-\lambda'-B_h^-\big((J(h)-\lambda)^{-1} - J(h)^{-1}\big)B_h^{+,*})\\
&=\tau_2^2 (\zzz_h-\lambda')(I+\mathcal O(h^{2}|\lambda|\|(\zzz_h-\lambda')^{-1})\|).
\end{align*}
Hence, according to the relations \eqref{eq:locspec2}, \eqref{eq:resolvest2}, and to $\varepsilon_j=\tau_2\tilde \varepsilon_j$,
it holds
\begin{align}
\nonumber
\tau_2^2 \nnn_h-\lambda-\tau_2^2 B_h^-(J(h)-\lambda)^{-1}B_h^{+,*}
&=\tau_2^2 (\zzz_h-\lambda')(I+\mathcal O(h^{2}\varepsilon_{j}^{2}\|(\zzz_h-\lambda')^{-1})\|)\\
\nonumber
&=\tau_2^2 (\zzz_h-\lambda')(I+\mathcal O(h\frac{\varepsilon_{j}^{2}}{\tilde\varepsilon_{j}^{2}}))\\
&=\tau_2^2 (\zzz_h-\lambda')(I+\mathcal O(h\tau_2^2)).
\label{eq.Schur}
\end{align}
The latter operator is then invertible around $\pa D'_{j,n}$
for $h,\tau_{2}$ small enough, and
the Schur complement formula then permits to
write the inverse of $\ppp(\lambda)$ as
\begin{equation}\label{eq:schur5}
\eee(\lambda)=
\left(
\begin{array}{cc}
E_0(\lambda)&-\tau_2 (J(h)-\lambda)^{-1}B_h^{+,*}E(\lambda)\\
-\tau_2 E(\lambda)B_h^{-}(J(h)-\lambda)^{-1}&E(\lambda)
\end{array}
\right)
\end{equation}
with 
$$
E(\lambda)=\Big(\tau_2^2 \nnn_h-\lambda-\tau_2^2 B_h^-(J(h)-\lambda)^{-1}B_h^{+,*}\Big)^{-1}
$$
and 
$$
E_0(\lambda)=(J(h)-\lambda)^{-1}+\tau_2^2(J(h)-\lambda)^{-1}B_h^{+,*}E(\lambda)B_h^-(J(h)-\lambda)^{-1}.
$$
As above, 
let us consider the
corresponding projector
$
\Pi_{ D'_{j,n}}(\mmm_h)$.
From $\lambda=\tau_2^2 \lambda'$, we get
$$
\Pi_{ D'_{j,n}}(\mmm_h)=\frac {\tau_2^2}{2i\pi}\int_{\partial \hat D'_{j,n}}\eee(\tau_2^2 \lambda')d\lambda'
$$
with $\hat D'_{j,n}=\{z'\in\C,\;\vert z'-\tilde\varepsilon_j^2\lambda^j_{n}\vert\leq Mh\tilde \varepsilon_j^2\}$.
It follows moreover from \eqref{eq.Schur}
that for every $\lambda'\in \pa \hat D'_{j,n}$ and $h,\tau_{2}$ small enough, 
%
%
\be\label{eq:resolvest3}
E(\tau_2^2\lambda')
=\tau_2^{-2}(\zzz_h-\lambda')^{-1}(I+\ooo(h)),
\ee
and the same argument as above shows that 
$\rank(\Pi_{ D'_{j,n}}(\mmm_h))\geq \rank(E'_n)$ 
with
\begin{equation*}
\begin{split}
E'_n=\frac {\tau_2^2}{2i\pi}\int_{\partial \hat D'_{j,n}}E(\tau_2^2 z)dz
=\frac 1{2i\pi}\int_{\partial \hat D'_{j,n}}(  \zzz_h-z )^{-1}(I+\ooo(h))^{-1}dz.
\end{split}
\end{equation*}
By the induction hypothesis, this shows that for $h$ small enough, the 
rank 
of $E'_n$  is exactly the multiplicity of $\lambda_n^j$ and hence
$$
\rank (\Pi_{D_{j,n}'}(\mmm_h))\geq m(\lambda_n^j)
$$
for all $j=2,\ldots,p$ and $n=1,\ldots,n_j$. Combined with \eqref{eq:minorrank}, this shows that 
for all $j=1,\ldots,p$ and $n=1,\ldots,n_j$, one has 
$$
\rank (\Pi_{D_{j,n}}(\mmm_h))\geq m(\lambda_n^j)
$$
with $D_{j,n}=\varepsilon^{2}_jD(\lambda_n^j,Mh)$.
Since $\sum_{j,n} m(\lambda_n^j)$ is equal to the total dimension of the space, this implies that 
\be\label{eq:minorrankproj4}
\rank (\Pi_{D_{j,n}}(\mmm_h))= m(\lambda_n^j)
\ee
which proves the localization of the spectrum and \eqref{eq:dimprojspect}.\medskip

\noindent
It remains to prove the resolvent estimate. Suppose that $\lambda\in\C$ is such that $\lambda\notin \cup_{j=1}^p\cup_{\mu\in \sigma(M_j)}
D( \varepsilon^{2}_j(\tau)\mu, \varepsilon^{2}_j(\tau)Kh)$. We suppose first that $\vert \lambda\vert\geq c_0$ for $c_0>0$ such that 
$\vert \lambda_{n}^1\vert\geq  2c_0$ for all $n=1,\ldots,n_1$. Then $\ppp(\lambda)=\mmm_h(\tau)-\lambda$ is invertible with inverse $\eee(\lambda)$ given by 
\eqref{eq:schur2}. Using \eqref{eq:invertshur1} it is clear that 
$E(\lambda)=\ooo(h^{-1})= \ooo(\dist(\lambda,\sigma(\mmm_h(\tau))^{-1})$. On the other hand, since 
$(\tau_{2}^{2}\nnn_h-\lambda)^{-1}=\ooo(1)$ and $B_h^{\pm}=\ooo(h)$ we have also 
$E_0(\lambda)=\ooo(1)$ and then $\eee(\lambda)=\ooo(\dist(\lambda,\sigma(\mmm_h(\tau)))^{-1})$.\medskip

\noindent
Suppose now that $\vert \lambda\vert\leq c_0$. Then $\ppp(\lambda)=\mmm_h(\tau)-\lambda$ is invertible with inverse $\eee(\lambda)$ given by 
\eqref{eq:schur5}. Setting $\lambda'=\tau_2^{-2}\lambda$ one deduces from \eqref{eq:resolvest3}
and from \eqref{eq:locspec2},\eqref{eq:resolvest2}  that 
$$
E(\lambda)=\ooo(\tau_2^{-2}\dist (\lambda',\sigma(\zzz_h))^{-1})=\ooo(\dist(\lambda,\sigma(\mmm_h(\tau))^{-1}).
$$
This completes the proof. 
 \ep


\bibliographystyle{amsplain}
\bibliography{ref_nonrevsmol}

\end{document}